\documentclass[10pt,a4paper]{article}
\usepackage[left=1in,right=1in,top=1in,bottom=1in]{geometry}

\usepackage[utf8]{inputenc}
\usepackage{pgfplots}
\pgfplotsset{compat=1.16} 
\usepackage{tikz}
\usetikzlibrary{automata,topaths,arrows,shapes}
\usetikzlibrary{backgrounds}
\usepackage{xcolor}
\usepackage{tabularx}
\usepackage{graphicx}
\usepackage{subfig}
\usepackage{textcomp}
\usepackage{gensymb}
\usepackage{amsmath}
\usepackage{amssymb}
\usepackage{multicol}
\usepackage{multirow}
\usepackage{varwidth}
\usepackage{ulem, hyperref}
\usepackage{authblk}
\usepackage{parskip}
\usepackage{adjustbox}

  \newcommand{\Bree}[1]{{\color{violet}\textbf{Bree:} #1}}

\title{
A yeast cell cycle pulse generator model shows consistency with multiple oscillatory and checkpoint mutant datasets}
\author[1]{Julian Fox\thanks{Corresponding author: julianfox8@gmail.com}}
\author[1]{Breschine Cummins}
\author[2]{Robert C. Moseley}
\author[3]{Marcio Gameiro}
\author[2]{Steven B. Haase}

\affil[1]{Department of Mathematical Sciences, Montana State University, Bozeman, MT, USA}
\affil[2]{Department of Biology, Duke University, Durham, NC, USA}
\affil[3]{Department of Mathematics, Rutgers University, New Brunswick, NJ, USA}

\date{\today}

\begin{document}

\maketitle

\begin{abstract}
% For years uncovering the regulatory mechanism driving progression of the cell cycle has been under investigation. Elucidating the relationship between cyclin dependent kinases (CDKs) and the transcription factor (TF) network has shed light on the regulatory mechanisms present that drive cell-cycle progression. The goal of this work was to highlight the oscillatory nature present within the TF network with and without the presence of CDKs, specifically Clb2. To do so we applied a systems biology level approach to investigate the dynamics present in a simplified hypothesized network model across wild-type and mutant datasets. It is shown that the simplified network model is capable of exhibiting the dynamics of interest which indicates a coupled interaction between the TF network and CDKs that regulates progression of the cell cycle. Additionally, our results introduce additions to the simplified network model which are hypothesized to improve the validity of the network model.   

The regulatory mechanisms driving progression of the yeast cell cycle appears to be comprised of an interacting network of transcription factors (TFs), cyclin-dependent kinases (CDK) and ubiquitin ligases. %The controlling regulatory network must be capable of driving periodic behavior. 
From a systems perspective the controlling regulatory network must produce robust periodic behavior during proliferative phases, but have the capability to halt the cycle when unfavorable conditions trigger a checkpoint arrest. How the individual components of the network contribute to these dynamical phenotypes remains an open question. Here we evaluate the capability of a simplified network model hypothesized to contain key elements of the regulation of cell-cycle progression to reproduce observed transcriptomic behaviors. We match time-series data from both cycling and checkpoint arrested cells to the predictions of a relatively simple cell-cycle network model using an asynchronous multi-level Boolean approach. We show that this single network model, despite its simplicity, is capable of exhibiting dynamical behavior similar to the datasets in most cases, and where it does not, we identified hypotheses that suggest missing components of the network. 

\end{abstract}

\section{Introduction}

The study of the cell cycle of simple organisms like the budding yeast \textit{Saccharomyces cerevisiae} can provide insight into more complex systems such as the mammalian cell cycle due to the similarity of the cellular machinery across eukaryotic organisms \cite{bertoli:2013,RN1839, RN1928}. Cell proliferation is a fundamental property of cells and understanding the regulation of the cell cycle is critical for explaining and controlling many biological processes from development to cancer progression. 

Models for the mechanisms controlling cell cycle progression have evolved over time. Early biochemical studies in marine invertebrates identified cyclins and cyclin dependent kinases (CDKs) as key regulators of cell-cycle oscillations along with the anaphase promoting complex (a ubiquitin ligase) \cite{RN772, RN361, RN495}. By forming a negative feedback loop, it was hypothesized that this simple cell-cycle network was capable of producing periodic behavior \cite{RN775, RN242}. Genetic studies in both \textit{S. cerevisiae} and \textit{S. pombe }also identified cyclins, CDKs, and APC complexes as critical regulators of cell cycle progression, indicating that machinery was highly conserved across phyla and that these components are the basic components of the cell-cycle oscillator \cite{RN58, RN869, RN385, RN485}. Cell-cycle models became more complex over time and mounting evidence indicated that oscillations of cell cycle events could continue when the CDK/APC motif was prevented from oscillating \cite{RN1728}. The advent of genome technologies identified large programs of dynamic gene expression associated with the cell cycle, and new models emerged suggesting the importance of transcriptional networks on producing cell-cycle oscillations \cite{RN1285, RN1452, RN1378, RN1282}.

Although ODE models of the yeast cell-cycle that are focused primarily on biochemical interactions have been remarkably predictive of mutant phenotypes \cite{RN1829}, there is a compelling argument that  transcription factors (TFs), cyclin-dependent kinases (CDKs), and ubiquitin ligases all play key roles in regulating cell-cycle progression \cite{RN1935}. Multiple studies have demonstrated that temporally ordered, high-amplitude transcript dynamics were present in budding yeast with non-oscillating levels of CDK \cite{RN1869, RN1939, RN1935, RN1452, RN1712} indicating that the CDK/APC oscillator identified in early embryonic systems may not be the core motif driving periodic behavior during the yeast cell cycle. 
 
Recent work suggested that a pulse-generating TF network containing an oscillatory mechanism was responsible for a transcriptional pulse that was thought to drive global phase-specific transcription \cite{RN1939}. The pulse generator seems to operate in a start-stop manner, where the network is first quiescent and then, after receiving a start signal, creates a wave of sequential transcription through the network which is hypothesized to be driven by the interaction of CDKs with a TF network \cite{orlando:2008,bristow:2014,cho:2019}. The transcriptional pulse driving progression through the cell cycle operates consistently, meaning that the gene products express in a stereotyped order \cite{cho:2019, powers:2017}, and the timing and robustness of this periodic transcription was affected when certain CDKs were knocked out or up-regulated, creating a weaker and less robust cycle \cite{RN1712}. It is important to note that periodic transcription, while weaker and less robust, was still present in some of the mutants. Therefore, the dynamics expressed by any hypothesized network model should exhibit oscillations under both wild-type (WT) and CDK mutant phenotypes.
 
Synchronous and autonomous Boolean models of yeast cell-cycle networks containing TFs and CDK have produced robust oscillatory behavior \cite{RN1452, RN1550}, but more sophisticated dynamical models that match observed dynamics (including transcriptional dynamics) across wild-type and mutant strains have not been reported.  Moreover, these models have not addressed the observation that certain environmental perturbations (e.g. DNA damage or spindle assembly defects) can reversibly arrest the cell cycle until damage is repaired.  
 
Here we describe a previously developed modeling framework, Dynamic Signatures Generated by Regulatory Networks (DSGRN)~\cite{cummins:2018,cummins:2016,cummins:2018,DSGRN_repo}, which is a type of asynchronous multi-level Boolean approach~\cite{crawford-kahrl_joint_2022} via the mechanism of switching systems~\cite{Thomas1991, veflingstad07,deJong2002,edwards00,Ironi2011}. The DSGRN software package can exhaustively compute all the dynamics a genetic regulatory network (GRN) can produce, allowing for a comprehensive description of potential network behaviors. It has been used to model and predict genetic network behavior in similar biological systems, such as the epithelial-to-mesenchymal transition in cancer \cite{xin:2020} and the Rb-E2F mechanism of the mammalian cell cycle~\cite{gedeon:2018}.

We investigated a simple integrated yeast cell-cycle network model that contains both TFs and CDKs, (Fig~\ref{fig:wavepool} (Right)). We observed that this network can match observed transcriptional dynamics from a variety of conditions where network components are mutated or checkpoint arrests are induced~\cite{orlando:2008,bristow:2014,cho:2017,cho:2019} as well as wild-type data, indicating that it may encode important features of cell-cycle regulation.

\section{Network modeling approach}\label{sec:phenotypes}

In this section, we discuss the evidence-based network model that we check for consistency with multiple wild-type and mutant datasets. Each of these datasets exhibits a cellular phenotype, namely, the cell cycle is either progressing or arrested. We distinguish between these cellular phenotypes and \textit{transcriptional phenotypes}
based on wild type and mutant microarray and RNAseq time series datasets from~\cite{orlando:2008,bristow:2014,cho:2017,cho:2019}. 
Transcriptional phenotypes consist of observed cycling or steady equilibrium behavior seen in the time series. 
We also define \textit{dynamical phenotypes} based on DSGRN predictions of network model dynamics that allow us to determine model consistency with the data. Dynamical phenotypes are graphical structures that capture stable or unstable cycling behavior as well as equilibria. 
%These are qualitative properties that track time series expression trajectories as discrete steps between qualitative states, i.e. high, low, and various grades of intermediate expression. To check for model consistency with data, the time series datasets are discretized in a unique way to construct target dynamical phenotypes (Section~\ref{sec:timeseries}). 
In order to be declared consistent with the observed data, a network model must be able to reproduce cyclic patterns in WT and mutant data and needs to support the arrest of cycling behavior during a triggered checkpoint.

\subsection{Wavepool model}

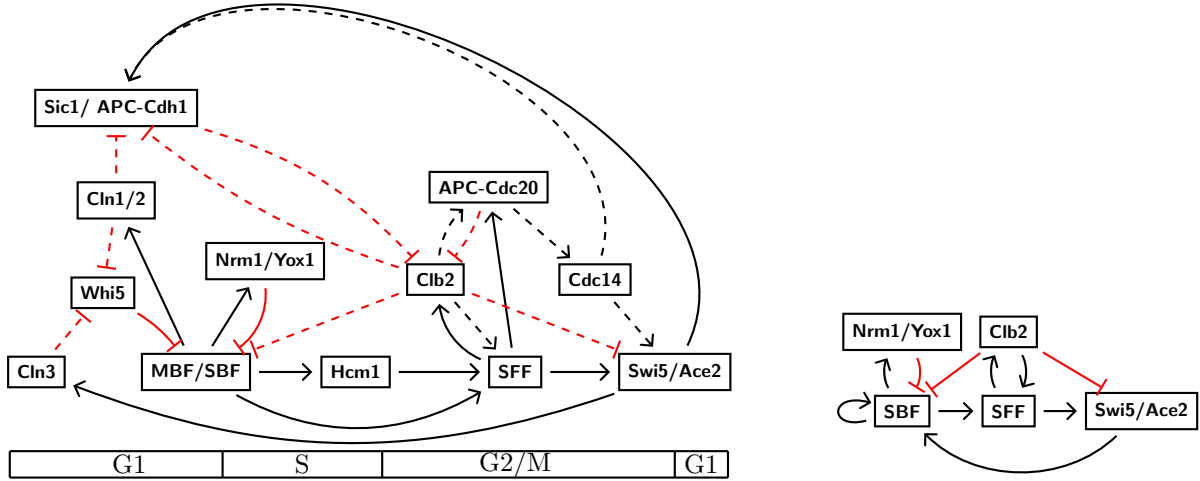
\begin{figure}[]
\centering
\begin{tabularx}{\textwidth}{cX}
\begin{tikzpicture}[main node/.style={rectangle,fill=white!20,thick,draw,font=\sffamily\scriptsize\bfseries},scale=0.7]

\draw[thick] (-.5,-.5)--(13,-.5);
\draw[thick] (-.5,-1)--(13,-1);
\draw[thick] (-.5,-.5)--(-.5,-1);
\draw[thick] (12,-.5)--(12,-1);
\draw[thick] (3.5,-.5)--(3.5,-1);
\draw[thick] (6.5,-.5)--(6.5,-1);
\draw[thick] (13,-.5)--(13,-1);

\node (b1) at (1.75,-.75) {G1};
\node (b2) at (5,-.75) {S};
\node (b3) at (9,-.75) {G2/M};
\node (n4) at (12.5,-.75) {G1};

\node[main node] (n1) at (0,1) {Cln3};
\node[main node] (n2) at (3,1) {MBF/SBF};
\node[main node] (n3) at (6,1) {Hcm1};
\node[main node] (n4) at (9,1) {SFF};
\node[main node] (n5) at (12, 1) {Swi5/Ace2};
\node[main node] (n6) at (1.25, 2.5) {Whi5};
\node[main node] (n7) at (4.3, 3.1) {Nrm1/Yox1};
\node[main node] (n8) at (7.5,2.75) {Clb2};
\node[main node] (n9) at (10.5,2.75) {Cdc14};
\node[main node] (n10) at (1.5,4.25) {Cln1/2};
\node[main node] (n11) at (1.5,6) {Sic1/ APC-Cdh1};
\node[main node] (n12) at (8.5,4.5) {APC-Cdc20};

\path[->,>=angle 90,thick]
(n1) edge[-|, dashed, red, shorten <= 3pt, shorten >= 3pt] node[] {} (n6)
(n2) edge[shorten <= 3pt, shorten >= 3pt] node[] {} (n3)
(n2) edge[shorten <= 3pt, shorten >= 3pt ,bend right] node[] {} (n4)
(n2) edge[shorten <= 3pt, shorten >= 3pt] node[] {} (n7)
(n2) edge[shorten <= 3pt, shorten >= 3pt] node[] {} (n10)
(n3) edge[shorten <= 3pt, shorten >= 3pt] node[] {} (n4)
(n4) edge[shorten <= 3pt, shorten >= 3pt] node[] {} (n5)
(n4) edge[shorten <= 3pt, shorten >= 3pt,bend left] node[] {} (n8)
(n4) edge[shorten <= 3pt, shorten >= 3pt] node[] {} (n12)
(n5) edge[shorten <= 3pt, shorten >= 3pt, bend right = 90] node[] {} (n11)
(n5) edge[shorten <= 3pt, shorten >= 3pt, bend left= 20] node[] {} (n1)
(n6) edge[-|, red, shorten <= 3pt, shorten >= 3pt, bend left = 10] node[] {} (n2)
(n7) edge[-|, red, shorten <= 3pt, shorten >= 3pt, bend left] node[] {} (n2)
(n8) edge[-|, dashed, red, shorten <= 3pt, shorten >= 3pt] node[] {} (n2)
(n8) edge[dashed, shorten <= 3pt, shorten >= 3pt] node[] {} (n4)
(n8) edge[-|,dashed, red, shorten <= 3pt, shorten >= 3pt, bend left = 10] node[] {} (n11)
(n8) edge[-|, dashed, red, shorten <= 3pt, shorten >= 3pt] node[] {} (n5)
(n8) edge[dashed, shorten <= 3pt, shorten >= 3pt,bend left = 20] node[] {} (n12)
(n9) edge[dashed, shorten <= 3pt, shorten >= 3pt, bend right = 85] node[] {} (n11)
(n9) edge[dashed, shorten <= 3pt, shorten >= 3pt] node[] {} (n5)
(n10) edge[-|, dashed, red, shorten <= 3pt, shorten >= 3pt] node[] {} (n11)
(n10) edge[-|, dashed, red, shorten <= 3pt, shorten >= 3pt] node[] {} (n6)
(n11) edge[-|, dashed, red, shorten <= 3pt, shorten >= 3pt, bend left = 15] node[] {} (n8)
(n12) edge[-|, dashed, red, shorten <= 3pt, shorten >= 3pt, bend left = 10] node[] {} (n8)
(n12) edge[dashed, shorten <= 3pt, shorten >= 3pt] node[] {} (n9)

;
\end{tikzpicture}

& 

\begin{tikzpicture}[main node/.style={rectangle,fill=white!20,thick,draw,font=\sffamily\scriptsize\bfseries},scale=0.7]
\node[main node] (n1) at (1,1) {SBF};
\node[main node] (n2) at (1,2.5) {Nrm1/Yox1};
\node[main node] (n3) at (3,1) {SFF};
\node[main node] (n4) at (3,2.5) {Clb2};
\node[main node] (n5) at (5.5,1) {Swi5/Ace2};

\path[->,>=angle 90,thick]
(n1) edge[shorten <= 3pt, shorten >= 3pt, bend left] node[] {} (n2)
(n1) edge[shorten <= 3pt, shorten >= 3pt,loop left] node[] {} (n1)
(n1) edge[shorten <= 3pt, shorten >= 3pt] node[] {} (n3)
(n2) edge[-|,red, shorten <= 3pt, shorten >= 3pt, bend left] node[] {} (n1)
(n3) edge[shorten <= 3pt, shorten >= 3pt, bend left] node[] {} (n4)
(n3) edge[shorten <= 3pt, shorten >= 3pt] node[] {} (n5)
(n4) edge[-|, red, shorten <= 3pt, shorten >= 3pt] node[] {} (n1)
(n4) edge[shorten <= 3pt, shorten >= 3pt, bend left] node[] {} (n3)
(n4) edge[-|, red, shorten <= 3pt, shorten >= 3pt] node[] {} (n5)
(n5) edge[shorten <= 3pt, shorten >= 3pt, bend left = 45] node[] {} (n1)
;                       
\end{tikzpicture}
\end{tabularx}

\caption{(Left) The wavepool model, or network oscillator, corresponding to Fig 7A in \cite{cho:2017}, including the pulse generator subnetwork composed of Nrm1/Yox1, MBF/SBF, SFF, Hcm1, Swi5, and Sic1. (Right) The mini wavepool, a gene regulatory network simplified from the wavepool model that attempts to capture the effect of Clb2 on pulse generator activity. We call the subnetwork without Clb2 the mini pulse generator, consisting of nodes SBF, SFF, Swi5/Ace2, and Nrm1/Yox1. Black arrows indicate activation, red blunt arrows indicate repression, dashed lines represent post-transcriptional interactions and solid lines represent transcriptional interactions. Notice that dashed lines are not included on the right, because we do not model post-transcriptional modification differently than transcriptional regulation in this paper. }
\label{fig:wavepool}
\end{figure}

A network oscillator model~\cite{cho:2017} representing a collection of regulatory interactions hypothesized to be capable of exhibiting multiple cell cycle behaviors is visualized as a GRN in Fig~\ref{fig:wavepool} (Left). We will call this network oscillator the wavepool model and the transcriptional subnetwork involving the nodes Nrm1/Yox1, MBF/SBF, SFF, Hcm1, Swi5, and Sic1 will be referred to as the pulse generator \cite{cho:2017}. The edges in the GRN reflect different regulatory mechanisms. A TF regulates another TF or CDK through transcriptional control of gene products. A CDK regulates another CDK or TF only when bound in a cyclin/CDK complex (post-transcriptional control). Once assembled, the complex is able to phosphorylate the target protein, which can have either an activating or inhibiting effect depending on the target. In Fig~\ref{fig:wavepool}, black arrows indicate activation, red arrows indicate inhibition, solid arrows indicate transcriptional control, and dashed arrows post-transcriptional control.

Each node within the network in Fig \ref{fig:wavepool} (Left) may not pertain to a specific TF or cyclin/CDK; in some cases the node refers to a complex of TFs or cyclins and CDKs. By convention, the name of a complex of proteins has all letters capitalized, the name of an individual protein has only the first letter capitalized, and the name of a gene or mutant gene has all letters lowercase and italicized. 

While it is understood that the CDKs and TFs identified in Fig \ref{fig:wavepool} (Left) are important for progression of the cell cycle, some research has focused on investigating the importance of smaller subnetworks \cite{orlando:2008,cho:2019,bristow:2014}; for example, concentrating on the impact that Clb2 activity has on its targets and the progression of the transcriptional program. We continue along these lines by computationally investigating a small network derived from the wavepool model. We will refer to this network as the mini wavepool, seen in Fig \ref{fig:wavepool} (Right). 
%Our intention is to explore the most minimal network that still contains the necessary interactions for the pulse generator to operate as well as being able to exhibit checkpoint activity. With this in mind, t
The mini wavepool was chosen to focus on the interactions between the single CDK, Clb2, and a cyclic arrangement of TFs that we refer to as the mini pulse generator. Clb2 is the B-type cyclin of choice because it is thought to play a role in the activation of multiple checkpoints and targets TFs \cite{bristow:2014} and because it acts broadly across cell cycle phases.

The size of the mini wavepool model was partially chosen for tractability of computation and many nodes and edges from the original hypothesized pulse generator are absorbed or disregarded. 
%Because we are interested in a minimal control mechanism of multiple checkpoints, w
% We chose to remove most of the post-transcriptional regulation acting on the mini pulse generator primarily due to the fact that during the preparation of this manuscript, DSGRN did not explicitly model post-transcriptional regulation (although this is a newly available feature~\cite{Cummins21}). For this reason, we removed the nodes Cdc14, APC-Cdc20, and Sic1/APC-Cdh1, and the cyclins Cln 1/2/3 have been absorbed into shortened paths. 
The self-edge on SBF comes from the length three path from MBF/SBF to itself through Cln1/2 and Whi5 in Fig~\ref{fig:wavepool}. Since there is a double repression, this indirect self-regulation is activating.  %The resulting self-edge is represented as post-transcriptional since the cyclins Cln1/2 are  contained within the collapsed path. 
Likewise, the edge from Swi5 to SBF in the mini wavepool is a collapsed path through double repression in Cln3 and Whi5. The removal of Sic1 is justified by evidence that Sic1 may not be necessary for cell-cycle oscillations \cite{RN1046}. Post-transcriptional regulation is not explicitly modeled as such; this is a recent advance in DSGRN~\cite{Cummins21}.

%Of our modifications, the least justifiable removal is the TF Sic1, acknowledged as part of the larger pulse generator~\cite{cho:2017}. The removal is performed primarily for computational reasons, but we also note that the effect of Sic1 on the pulse generator in Fig~\ref{fig:wavepool} is completely mediated through Clb2. There is also evidence that Sic1 is not necessary for cell-cycle oscillations \cite{RN1046}.

\subsection{Transcriptional Phenotypes}

We describe seven different transcriptional phenotypes expressed in seven datasets. The transcriptional phenotypes are called WT, Clb2 OFF, Clb2 ON, Clb2 intermediate-low (INT-L), Clb2 intermediate-high (INT-H), spindle assembly checkpoint (SAC), and DNA replication checkpoint (DRC).

The wild-type (WT) microarray dataset comes from \cite{orlando:2008} where genome-wide transcription was analyzed for \textit{S. cerevisiae}. The WT dataset acts a baseline for the dynamics present within the cell cycle. Oscillations were seen not only within the transcription factor network but within CDKs as well. 

 We identify the Clb2 ON mutant with a \textit{cdc20$\Delta$} knockout mutant~\cite{bristow:2014}. Cdc20 acts to promote progression through the cell cycle in M phase from metaphase to anaphase. This knockout leads to cells arrested at the metaphase-to-anaphase border in mitosis containing high levels of non-oscillating Clb2 protein. However, transcriptional oscillations in the pulse generator persist. We identify the Clb2 OFF mutant with a \textit{clb$\Delta$} mutant in which \textit{clb1-6} are knocked out~\cite{orlando:2008}. These cyclins (Clb1-6)  act as the S-phase and mitotic cyclins; therefore this mutant was arrested at the G1/S border. It was seen that these mutants were not able to replicate DNA, separate spindle pole bodies, undergo isotropic bud growth, or complete nuclear division, indicating that these cells do not have active Clb-CDK complexes \cite{orlando:2008} and that therefore these mutants were void of any nontrivial Clb2 activity. Again, the data show that transcriptional oscillations in the pulse generator persist despite the cell cycle arrest.

The two mutant phenotypes with intermediate expression are identified with the \textit{cdc14-3} and \textit{cdc15-2} mutants from \cite{cho:2017}. 
%Clb2 INT-H represents the \textit{cdc14-3} mutant and the Clb2 INT-L represents the \textit{cdc15-2} mutant. 
Instead of knockouts these mutants are temperature sensitive, thus when incubated at 37\degree{C} there is a drop off in the activity of these genes. Cdc14 and Cdc15 are key players in a cell's exit from mitosis \cite{powers:2017,baumer:2000}. 
There are two APC complexes responsible for lowering Clb activity, APC/Cdc20 and APC/Cdh1. Cdc15 and Cdc14 are important for the formation of the APC/Cdh1 complex, therefore the \textit{cdc14-3} and \textit{cdc15-2} mutants are thought to only have APC/Cdc20 present. It was observed that the temperature sensitive \textit{cdc14-3} and \textit{cdc15-2} mutants had moderate levels of non-oscillating Clb2 activity \cite{cho:2017,baumer:2000}.  Cdc14 acts to lower Clb2 activity directly through the phosphorylation of CDKs, while Cdc15 simply acts to re-initiate Cdc14. Thus, in the \textit{cdc14-3} mutant there is no degradation of Clb2 via Cdc14 activity and within the \textit{cdc15-2} mutant there is the initial degradation of Clb2 through Cdc14, but Cdc14 is never re-initiated by Cdc15. This indicates that the \textit{cdc14-3} mutant should exhibit higher levels of Clb2 activity than the \textit{cdc15-2} mutant, and this reasoning is supported by the experimental data. We name the corresponding phenotypes Clb2 INT-H and Clb2 INT-L accordingly. In these two mutant datasets, pulse generator transcriptional oscillations are present. 

We identify the SAC phenotype with a \textit{cse4} mutant that contains a mutant allele of the kinetochore protein, which acts to disrupt the spindle assembly \cite{bristow:2014} due to improper orientation of the kinetochore. If this occurs, signals are sent from the SAC mechanism to initiate the formation of the mitotic checkpoint complex (MCC), which binds to APC and inhibits the formation of the APC/Cdc20 complex through competition \cite{Yang:2015,Tavormina:1998}. In this mutant, pulse generator transcriptional oscillations are silenced and the resulting steady state serves as the transcriptional SAC phenotype.
% \textcolor{red}{We identify the SAC phenotype with a \textit{cse4} mutant that contains a mutant allele of the kinetochore protein, which acts to disrupt the spindle assembly \cite{bristow:2014} due to improper orientation of the kinetochore. If this occurs, signals are sent from the SAC mechanism to checkpoint proteins such as Bub1, Bub2, Mad1, and Mad2. These signals initiate the formation of the mitotic checkpoint complex (MCC), which binds to APC and inhibits the formation of the APC/Cdc20 complex through competition \cite{Yang:2015,Tavormina:1998}. In addition, it has been seen that Cdc20 is targeted for degradation via ubiquitination from APC which is crucial for maintaining mitotic arrest during SAC \cite{Wang:2017,Nilsson:2008}. Thus, SAC arrest is hypothesized to be driven by two mechanisms, one being the activation of the MCC and the other being ubiquitination of Cdc20 by APC.   In this mutant, pulse generator transcriptional oscillations are silenced and the resulting steady state serves as the SAC phenotype.} 

The DNA replication checkpoint (DRC) mutant is defined as the \textit{cdc8} mutant from~\cite{bristow:2014} that contains a temperature sensitive allele of \textit{cdc8} that disrupts the thymidylate kinase responsible for DNA synthesis \cite{Sclafani:1984}. Given a disruption to the progression of the DNA replication fork, the thymidylate kinase cascade regulates timing of S-phase events, such as spindle elongation, allowing the cell to deal with insults to the DNA replication fork by arresting the cell cycle. The pulse generator is non-oscillatory in this dataset, with the resulting steady state defining a transcriptional DRC phenotype.

To investigate the mini wavepool it is necessary to ascertain qualitative transcript levels for all nodes in the network. As is clear from the presence of complexes and boxed TFs in Fig~\ref{fig:wavepool}, multiple proxies are available for each node in the mini wavepool. The determination of qualitative transcript levels for the mini wavepool involved the choice of proxies for SBF and SFF complexes, which we chose to be Swi4 and Ndd1 respectively. The node Nrm1/Yox1 naturally has two proxies, Nrm1 and Yox1. 
% While Yhp1 is also a TF of interest, it exhibits a large double peak that is inconsistent with the expression of its paralog Yox1 and the other pulse generator TFs. \Bree{[I suggested swapping out Yhp1 for Yox1 in Fig~\ref{fig:wavepool}.]}
% \Bree{[put picture in appendix/supp]} 
Both Swi5 and Ace2 were explored as proxies to represent the Swi5/Ace2 node as pictured in Fig~\ref{fig:wavepool}. We evaluate each of the transcriptional phenotypes for each of the proxy sets listed in Table \ref{fig:proxy_decisions}. 

\begin{table}[]
\centering
\begin{tabular}{|llllll|}
\hline
\multicolumn{6}{|c|}{Mini wavepool proxy decisions}                        \\ \hline
                                & \multicolumn{5}{c|}{Mini wavepool nodes} \\ \cline{2-6} 
\multicolumn{1}{|l|}{Proxies}   & SBF    & Nrm1/Yox1   & SFF   & Clb2  & Swi5  \\ \hline
\multicolumn{1}{|l|}{Swi5-Nrm1} & Swi4   & Nrm1   & Ndd1  & Clb2  & Swi5  \\
\multicolumn{1}{|l|}{Swi5-Yox1} & Swi4   & Yox1   & Ndd1  & Clb2  & Swi5  \\
\multicolumn{1}{|l|}{Ace2-Nrm1} & Swi4   & Nrm1   & Ndd1  & Clb2  & Ace2  \\
\multicolumn{1}{|l|}{Ace2-Yox1} & Swi4   & Yox1   & Ndd1  & Clb2  & Ace2  \\ \hline
\end{tabular}
\caption{Proxy choices for the mini wavepool.}
\label{fig:proxy_decisions}
\end{table}

% We next describe corresponding computational phenotypes that are expressed as stable dynamical behaviors arising from the DSGRN discrete approach. We provide a brief overview of DSGRN in order to explain these ``DSGRN dynamical phenotypes''. A more detailed explanation of DSGRN is provided in Section~\ref{sec:methods}.

\subsection{Dynamical Phenotypes}

In the traditional ODE modeling paradigm, a
mechanistic model is constructed, often using Hill functions in the genetic network setting, and numerous parameters are either fit to data or drawn from the literature. Frequently, a sensitivity analysis is performed to check the variability of model output to small perturbations in parameters. This approach, while valuable, has limitations. First, the area of parameter space that can be explored is minuscule, and second, the large number of parameters can lead to overfitting. We offer an alternative framework in which large regions of parameter space are excluded as unable to produce the desired dynamical behavior using a multi-level Boolean modeling approach. After identifying a reduced parameter space, parameterized Hill models can be created to replicate the desired behavior.

The modeling framework Dynamical Signatures Generated by Regulatory Networks (DSGRN)~\cite{cummins:2018,cummins:2016,cummins:2018,DSGRN_repo} is based on an ODE system with switching functions (Section~\ref{sec:dsgrn}) that leads to a fundamentally different approach. First, DSGRN provides a searchable database  of coarse dynamics over the entirety of parameter space (Sections~\ref{sec:pg}-\ref{sec:stg}). This is possible because DSGRN divides parameter space into a finite number of regions, each called a DSGRN parameter node. The coarse dynamics exhibited by a network model are constant for all real-valued parameter sets sampled from  the same DSGRN parameter node. Second, DSGRN uses only coarse information from a time series dataset, which avoids overfitting. 

DSGRN output contains information on the number and type of fixed points (equilibria) and oscillations that the network can exhibit at a given DSGRN parameter node. The fixed points identified by DSGRN are not quantitative, they only indicate whether the associated gene product is predicted to have high, low, or intermediate concentrations. Likewise, cycling behavior is not modeled by a smooth trajectory, but rather a sequence of maxima and minima for each gene's expression level.

Model consistency in cycling behavior is determined by what we call a pattern match (Sections~\ref{sec:timeseries}-\ref{sec:patternmatch}) between the sequence of maxima and minima predicted by a network model and the observed sequence of maxima and minima in the data. Model consistency for a fixed point involves a judgment call on whether gene expression is best described as high, low, or intermediate at the time when steady behavior is observed. For both types of model consistency, the proportion of DSGRN parameter nodes that exhibit the observed behavior indicates how robustly the network model recapitulates the data.

Each DSGRN parameter can be decomposed into a collection of independent DSGRN factor parameters, one per node in a GRN. In the mini wavepool, we will distinguish between the DSGRN factor parameter for Clb2 (the Clb2 parameter, Section~\ref{sec:methods:clb2}) and the collection of remaining DSGRN factor parameters, which we will call the mini pulse generator parameter. In particular, we explore the behavior of the mini wavepool as the Clb2 factor parameter changes, but the mini pulse generator parameter remains fixed, mimicking a control mechanism of Clb2 on the mini pulse generator.

 We devise three different dynamical phenotypes based on DSGRN output that are intended to represent important features of the experimental data and associated transcriptional phenotypes; see the correspondence in Table~\ref{tab:dsgrn_pheno}. 
%  If the mini wavepool model is consistent with the WT phenotype, then the mini wavepool model captures the  pulse generator behavior of the undisrupted cell cycle. If the model is consistent with all dynamical phenotypes, then it is a complete hypothesis for explaining the seven datasets.  If the WT phenotype, but not all mutant phenotypes are supported by the model, then it is an indication that the model is incomplete and requires additional regulatory molecules.
 Dynamical phenotype I (WT cycling) is a pattern match to the wild-type data within a stable DSGRN cycle. This is roughly analogous to ensuring similar phase relationships between expression levels that are robust to small perturbations. We identify the DSGRN parameter nodes at which consistency with the WT cycling in the data occurs.
 
 \begin{table}[!h]
\centering
\begin{tabular}{|l|l|l|}
\hline
Dynamical phenotype        & Transcriptional phenotypes & Datasets \\ \hline
I:~~~~~~~~WT cycling       & WT  &  WT      \\ \hline
II:~~~Clb2 mutant cycling &  \begin{tabular}[c]{@{}l@{}} Clb2 ON, Clb2 OFF,\\ Clb2 INT-H \& Clb2 INT-L\end{tabular}  &    \begin{tabular}[c]{@{}l@{}}\textit{cdc20$\Delta$}, \textit{clb$\Delta$},\\ \textit{cdc14-3} \& \textit{cdc15-2} \end{tabular}   \\ \hline
III:~~~~checkpoint arrest &  SAC \& DRC &  \textit{cse4} \& \textit{cdc8}       \\ \hline
\end{tabular}
\caption{The mapping between dynamical and transcriptional phenotypes and their associated datasets~\cite{orlando:2008,bristow:2014,cho:2017}.}
\label{tab:dsgrn_pheno}
\end{table}

% \Bree{Mention that this is co-existence of WT+mutant at the same pulse generator DSGRN parameter. Introduce the idea of splitting DSGRN parameter into a mini pulse-generator part and a Clb2 part, see Methods Sections~\ref{sec:pg}-\ref{sec:methods:clb2}.}

Dynamical phenotype II captures Clb2 mutant cycling in which only the expression levels in the mini pulse generator are stably oscillating. We identify those DSGRN parameter nodes that exhibit pattern matches to Clb2 mutant time series when Clb2 is fixed at an appropriate high, low, or intermediate level, depending on the dataset. In the mini wavepool model in Fig~\ref{fig:wavepool} (Right), Clb2 can attain four discrete states, low (0), intermediate-low (1), intermediate-high (2), and high (3), which is dictated by the number of regulatory targets of Clb2 (see Section~\ref{sec:dsgrn}). We can fix Clb2 at any state $i$ by restricting ourselves to an appropriate collection of DSGRN parameter nodes (see Section~\ref{sec:methods:clb2}). For the Clb2 ON transcriptional phenotype, we fix the Clb2 state at 3, for Clb2 INT-H the fixed state is 2, for Clb2 INT-L, the state is 1, and for Clb2 OFF, the state is 0. For example, a pattern match between the mini wavepool model and the mutant dataset \textit{cdc20}$\Delta$ is only considered successful if it occurs at a DSGRN parameter node where the Clb2 state is fixed at 3; similarly for the other datasets listed in Table~\ref{tab:dsgrn_pheno}, row 2. We say that a DSGRN parameter is ``phenotype-permissible'' if it obeys the appropriate constraint for a given dataset. 

Suppose two DSGRN parameter nodes differ only in Clb2, are identical in the mini pulse generator and there is WT cycling at one parameter node and there is mutant cycling at the other. This dynamical connection is made more rigorous in Section~\ref{sec:methods:clb2}.  In dynamical phenotype II, we seek such parameter pairs, and the implication of the existence of such a pair is that the mini pulse generator is predicted to operate independently of changes to Clb2. 
% On the other hand, if there is mutant cycling at a DSGRN parameter node, but there is no other DSGRN parameter node with WT cycling at the same fixed mini pulse generator parameter, then the model implication is that a regulatory event impacted the mini pulse generator. This regulation could have occurred either through Clb2 alone, or from an unmodeled regulator.

% We seek a specific relationship between DSGRN parameter nodes that exhibit WT cycling and those that exhibit Clb2 mutant cycling. When two DSGRN parameter nodes differ only in Clb2 and are identical in the mini pulse generator, we call this pair Clb2-divergent.   When a Clb2-divergent parameter pair exists between WT and mutant cycling, the parameters for the mini pulse generator undergo no changes from WT when Clb2 levels are fixed, indicating that the mini pulse generator can operate independently of changes to Clb2. On the other hand, if the model does exhibit mutant cycling but no Clb2-divergent parameter pairs exist, then the model implication is that a regulatory event impacted the mini pulse generator. This regulation could have occurred either through Clb2 alone, or from an unmodeled regulator.

\begin{figure}[]
    \centering
    \subfloat[\centering \textit{cdc8} ts with Swi5 and Nrm1]{{\includegraphics[width=0.32\textwidth]{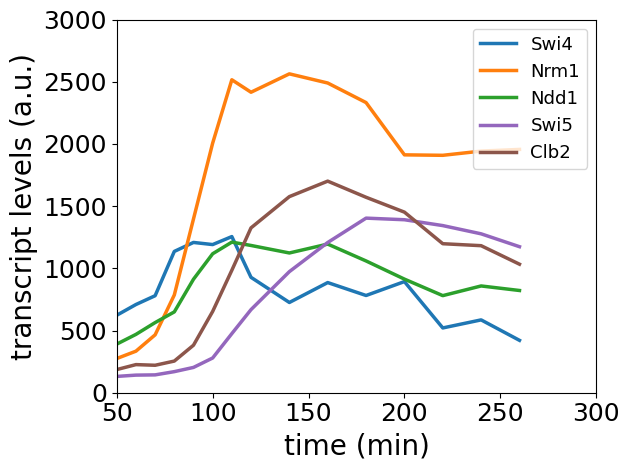} }}\hfill
    \subfloat[\centering WT ts with Swi5 and Nrm1]{{\includegraphics[width=0.32\textwidth]{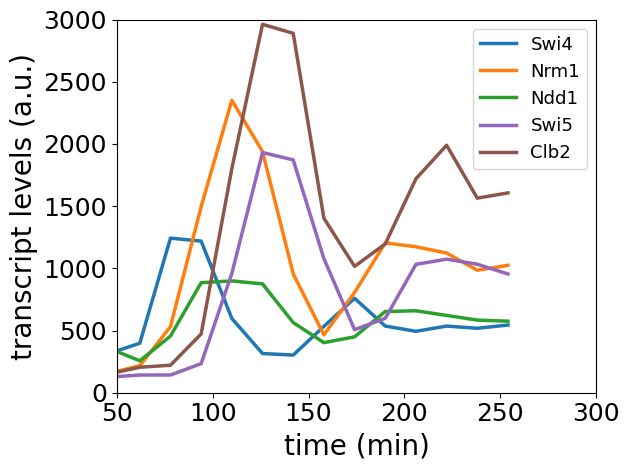} }}\hfill
    \subfloat[\centering \textit{cse4} ts with Swi5 and Nrm1]{{\includegraphics[width=0.32\textwidth]{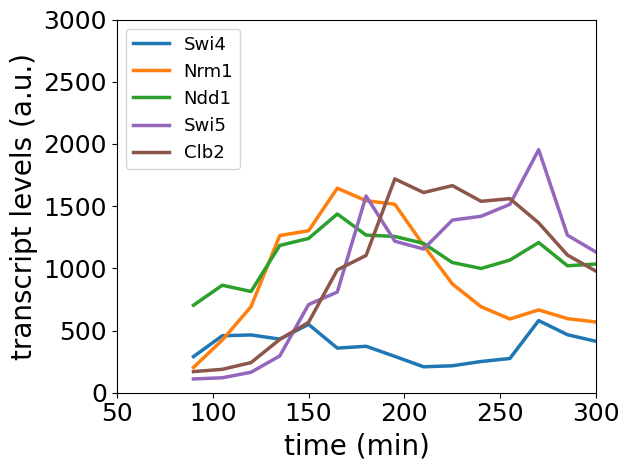} }} \par
    \subfloat[\centering \textit{cdc8} ts with Ace2 and Yox1]{{\includegraphics[width=0.32\textwidth]{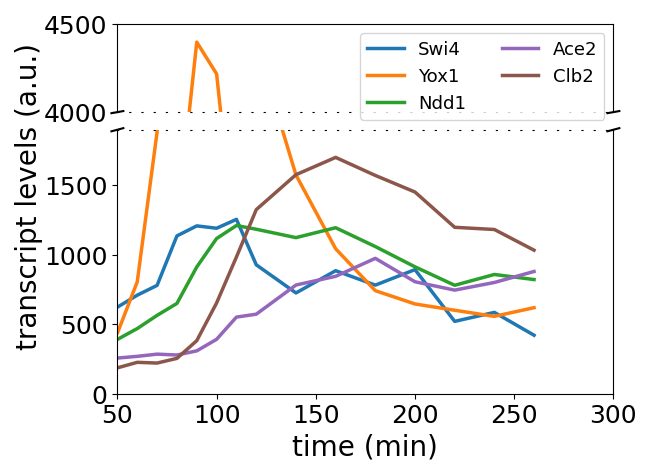} }}\hfill
    \subfloat[\centering WT ts with Ace2 and Yox1]{{\includegraphics[width=0.32\textwidth]{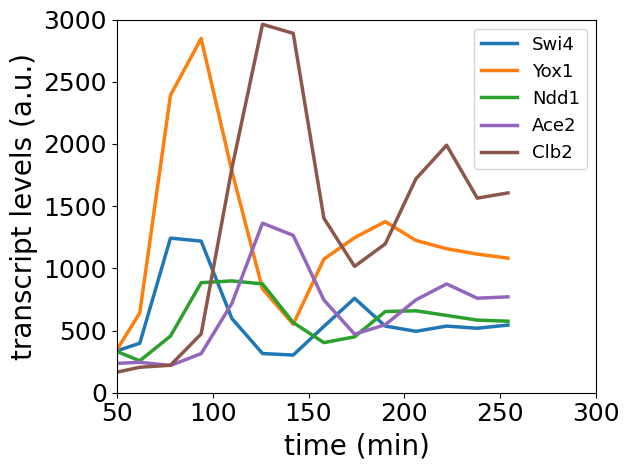} }}\hfill    
    \subfloat[\centering \textit{cse4} ts with Ace2 and Yox1]{{\includegraphics[width=0.32\textwidth]{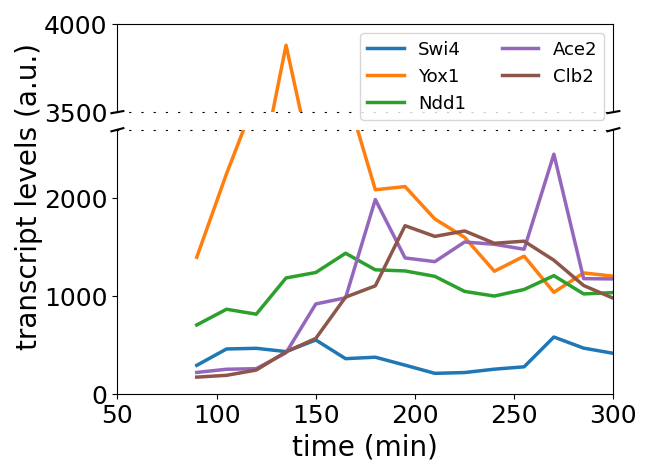} }} 
    \caption{(a) and (d) The transcript levels for proxy sets in the \textit{cdc8} mutant exhibiting DRC arrest. (b) and (e) The transcript levels for proxy sets in the WT data. (c) and (f) The transcript levels for proxy sets in the \textit{cse4} mutant exhibiting SAC arrest. The datasets have been jointly normalized and WT data is in the middle to facilitate comparison with both checkpoint datasets.}
    \label{fig:ts_poset}%
\end{figure}

In dynamical phenotype III (checkpoint arrest), we identified a DSGRN fixed point (FP) surrogate of the SAC equilibrium and of the DRC equilibrium. We examined the data in \cite{bristow:2014} for the \textit{cse4} mutant (SAC), and the \textit{cdc8} mutant (DRC). Qualitatively high, intermediate, and low values at equilibrium were assigned for each network node based on a comparison to the wild-type transcriptomics data. We remark that all seven datasets were normalized together to permit such comparisons (see Data Availability). The activity of each node was determined by comparing the transcript level at the end of the cycle in the \textit{cse4} mutant data or the \textit{cdc8} mutant data to the transcript level of the same node in the WT data, see Fig~\ref{fig:ts_poset} and Table~\ref{tab:qualitative_fps}. Data points from the end of the time series were used as the best representation of the resulting steady state. 

\begin{table}[h!]
\centering
\begin{tabular}{|lllllll|}
\hline
Checkpoint & Proxies & Swi4 & Nrm1/Yox1 & Ndd1 & Swi5/Ace2 & Clb2 \\ \hline
SAC & All proxies & low & low & high & high & intermediate/high \\ \hline
\multirow{2}{*}{DRC} & Yox1 proxies& low & low & high & high & intermediate/high
\\ & Nrm1 Proxies& low & high & high & high & intermediate/high
\\ \hline
\end{tabular}
\caption{The fixed points representing the SAC and DRC based on the comparison between the WT dataset and the \textit{cse4} and \textit{cdc8} datasets, respectively  (see Fig~\ref{fig:ts_poset}).} 
\label{tab:qualitative_fps}
\end{table}

The transcript levels to determine the SAC FP were based upon a maximum level of transcript seen across the SAC mutant compared to the maximum transcript level seen in the WT time series. Low represents any transcripts judged to be substantially below 50\% maximum transcript level and high represents any transcript judged to be substantially above 50\% maximum transcript level. Intermediate levels are judged to be in-between. The SAC FP is given by Swi4 = low, Nrm1/Yox1 = low, Ndd1 = high, Clb2 = not low, and Swi5/Ace2 = high (see Table~\ref{tab:qualitative_fps} first row). Notice that Clb2 expression levels are largely indeterminate since they have not achieved steady state by the end of the time series. 
%This leads to a collection of multiple possible SAC FPs that differ in Clb2, rather than just one FP.

Similar analysis was done for the DRC mutant (Table~\ref{tab:qualitative_fps} second row). For the Yox1 proxies, the SAC and DRC FPs are identical, indicating that the mini wavepool has an insufficient diversity of nodes to distinguish between these two checkpoints. For the Nrm1 proxies, the fixed points differ at the Nrm1 value, where the SAC exhibits low Nrm1 activity and the DRC exhibits high Nrm1 activity.

As in mutant cycling, we checked 
%both for the existence of SAC/DRC FPs  over all DSGRN parameter nodes and 
for the existence of DSGRN parameter node pairs  where one showed WT cycling and the other a SAC/DRC FP with the same mini pulse generator parameter. Unlike in mutant cycling, the arrest at a checkpoint in phenotype III indicates different dynamical behavior in the mini pulse generator as compared to WT. Therefore, the existence of such a pair indicates that regulation through Clb2 alone is sufficient to control entry into a checkpoint within the mini wavepool model without additional external regulation at other network nodes. It is important to stress that this does not exclude the existence of other regulators in a larger network in the cell; it merely indicates that the mini wavepool model as constructed can replicate the mutant phenotype of interest.

% \begin{table}[]
% \centering
% \begin{tabular}{|lll|}
% \hline
% \multicolumn{3}{|l|}{Checkpoint FPs: FP(Swi4, Nrm1/Yox1, Ndd1, Clb2, Swi5/Ace2)}                                                                                                                                                                                                \\ \hline
% \multicolumn{1}{|l|}{SAC FPs:} & \multicolumn{2}{l|}{\begin{tabular}[c]{@{}l@{}}All proxies: \\  FP(0, 0, 2, 1, 1) \\  FP(0, 0, 2, 2, 1) \\  FP(0, 0, 2, 3, 1)\end{tabular}}                                                                                                     \\ \hline
% \multicolumn{1}{|l|}{DRC FPs:} & \begin{tabular}[c]{@{}l@{}}Yox1 proxies:\\  FP(0, 0, 2, 1, 1) \\  FP(0, 0, 2, 2, 1) \\  FP(0, 0, 2, 3, 1)\end{tabular} & \begin{tabular}[c]{@{}l@{}}Nrm1 Proxies:\\  FP(0, 1, 2, 1, 1) \\  FP(0, 1, 2, 2, 1) \\  FP(0, 1, 2, 3, 1)\end{tabular} \\ \hline
% \end{tabular}
% \caption{The order of the nodes for the FPs is Swi4, which has activity levels in the range 0-3, Nrm1/Yox1, which has levels 0-1, Ndd1, which has levels 0-2, Clb2, which has levels 0-3, and Swi5/Ace2, which has levels 0-1. See Section~\ref{sec:dsgrn} for an explanation of the differing state discretizations for each of the nodes. \Bree{We should probably move this to Methods or an appendix.}} 
% \label{fig:sac_fps}
% \end{table}

\section{Results}\label{sec:results}

We assess the consistency of the predictions of the mini wavepool model with the datasets by checking for the existence of the DSGRN dynamical phenotypes at a collection of DSGRN parameters. The number of DSGRN parameter nodes for the mini wavepool is quite large: 275,466,240 total.
We first divide DSGRN parameters into five distinct groups as described below (see Section~\ref{sec:methods:clb2} for technical details). This division of DSGRN parameter space is based solely on the Clb2 DSGRN factor parameter, and not on the mini pulse generator parameter. 

A full 60\% of DSGRN parameter space allows for changing levels of Clb2 and therefore has the capacity for exhibiting the wild-type cycling phenotype. The remaining 40\% of DSGRN parameter space is composed of parameters that have a fixed value of Clb2, with 10\% each at high, low, and intermediate-high/low. We enforce that the mini wavepool can only be consistent with the Clb2 mutant datasets at the appropriately assigned parameters for Clb2 ON, OFF, INT-H, and INT-L (Clb2 mutant cycling phenotypes). The checkpoint datasets do not provide any \textit{a priori} parameter constraints  and consistency is assessed over 100\% of DSGRN parameter space (checkpoint arrest phenotypes). The subsets of DSGRN parameter nodes over which matches to given datasets are searched are called are called \textit{phenotype-permissible} parameters. In particular, WT cycling phenotype-permissible parameters comprise 60\% of parameter space, mutant cycling phenotype-permissible parameters comprise 10\% of parameter space each, and checkpoint phenotype-permissible parameters are unconstrained.

Our results are reported as proportions over subsets of DSGRN parameter space. Dynamical phenotype I results are the percentages of WT cycling matches over 60\% of parameter space. From this set of DSGRN parameters, the number of mini pulse generator (MPG) parameters was extracted and used as the normalization term for dynamical phenotypes II and III. This was justified because we sought to discover regions of parameter space where a change in Clb2 was the only impact on the network behavior; therefore, we looked for MPG parameters that showed both WT cycling and matched at least one other phenotype.
% After checking for consistency with each mutant dataset independently, we examine at which mini pulse generator parameters multiple mutant phenotypes can be observed along with WT behavior.

\subsection{Dynamical Phenotype I: Wild-Type Cycling}\label{sec:pheno_I}

The purpose of dynamical phenotype I is to show that the mini wavepool model can match the normal function of the cell cycle, namely the experimental WT dataset, by identifying the percent of phenotype-permissible DSGRN parameter nodes that can be pattern matched in a stable cycle. Table \ref{tab:pheno1_wt} shows that the proxy groups have WT pattern matches ranging from 3.7\% up to 4.3\% of phenotype-permissible DSGRN parameters. This demonstrates that every proxy group can recapitulate the oscillating sequence of maxima and minima seen in the WT data, and therefore the mini wavepool model cannot be excluded as a reasonable network model for controlling the yeast cell cycle. 

\begin{table}[h]
\centering
\begin{tabular}{|llcc|}
\hline
\multicolumn{4}{|c|}{DSGRN Phenotype I}                             \\ \hline
proxies                         & matches & \% of matches & \multicolumn{1}{l|}{MPG parameters} \\ \cline{2-4} 
\multicolumn{1}{|l|}{Swi5-Nrm1} & 7073137 & 4.3\%  &  1541616                           \\
\multicolumn{1}{|l|}{Ace2-Nrm1} & 7073137 & 4.3\%  &   1541616                        \\
\multicolumn{1}{|l|}{Swi5-Yox1} & 6071292 & 3.7\%  &   1347895                          \\
\multicolumn{1}{|l|}{Ace2-Yox1} & 6071292 & 3.7\%  &   1347895                          \\ \hline
\end{tabular}
\caption{Dynamical phenotype I: wild-type cycling results. The raw number of matches together with proportions of phenotype-permissible parameters are shown in columns 1 and 2. Column 3 shows the number of distinct mini pulse generator (MPG) parameters within the WT cycling matches. These numbers are the normalization constants in Fig~\ref{fig:phenoII-III}.}
\label{tab:pheno1_wt}
\end{table}

\begin{figure}
    \centering
    \begin{tabular}{cc}
\includegraphics[width=3in]{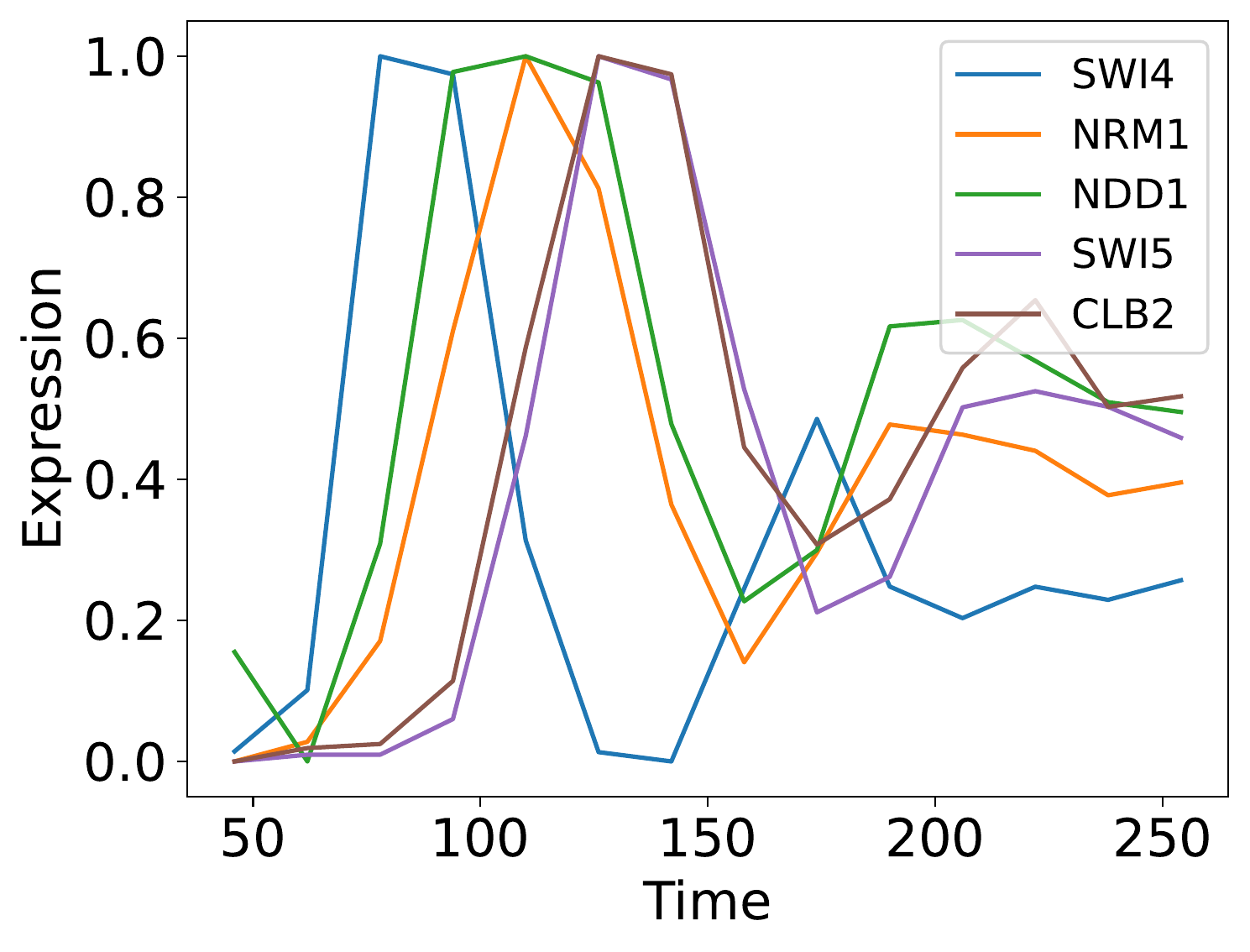} & \includegraphics[width=3in]{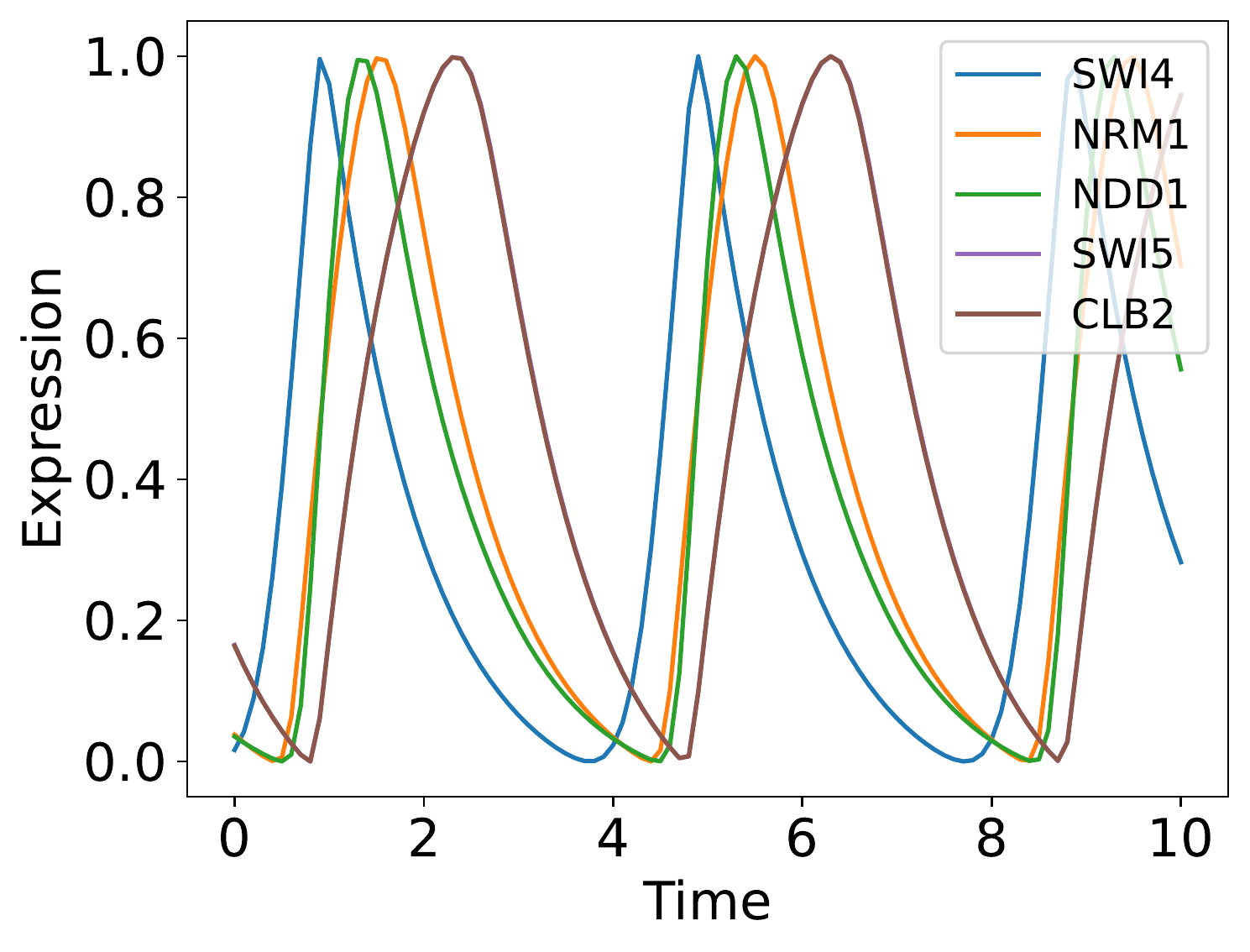} \\ (a) & (b)
\end{tabular}
    \caption{\textbf{Hill model simulation for wild type cells.} (a) Wild type data normalized between zero and one. (b) Hill model simulation using a DSGRN MPG parameter predicted to exhibit wild type oscillations, also normalized. The curves for Swi5 and Clb2 are (nearly) coincident. No attempt was made to match amplitude or period, only the ordering of extrema. }
    \label{fig:wt_hill_model}
\end{figure}

In addition, matching WT cycling greatly reduced parameter space. We parameterized Hill models from this reduced region (see Section~\ref{sec:hill}) to demonstrate the utility of DSGRN in locating useful regions of parameter space. The result of one such simulation is shown in Fig~\ref{fig:wt_hill_model}(b) with the comparable WT data shown in Fig~\ref{fig:wt_hill_model}(a). Both sets of curves are normalized between 0 and 1 to emphasize the order of maxima and minima in the time series, which is the behavior that DSGRN predicts. We see that in both cases, the peak of Swi4 precedes the peaks of Nrm1 and Ndd1, which themselves precede the peaks of Swi5 and Clb2. Swi5 and Clb2 have nearly indistinguishable peaks in both panels. The order of the Ndd1 and Nrm1 peaks are indistinguishable in the WT data in panel (a), but are ordered with Ndd1 first in the simulation in (b). Since Ndd1 has a faster rise in the WT data, the simulation ordering shows consistency with the WT data. We regard the recapitulation of this order of extrema as a successful model; however, we made no attempt to match either amplitude or period of the WT data.

We determined the number of unique mini pulse generator parameters in the set of WT cycling pattern matches, shown in Table \ref{tab:pheno1_wt}, column 3 for each proxy. These numbers are the normalization factors for the dynamical phenotypes II and III results shown in Fig~\ref{fig:phenoII-III}.

\subsection{Dynamical Phenotype II: Clb2 mutant cycling}\label{sec:pheno_II}

\begin{figure}[!h]
    \centering
    \subfloat[\centering Swi5-Nrm1]{{\includegraphics[width=0.475\textwidth]{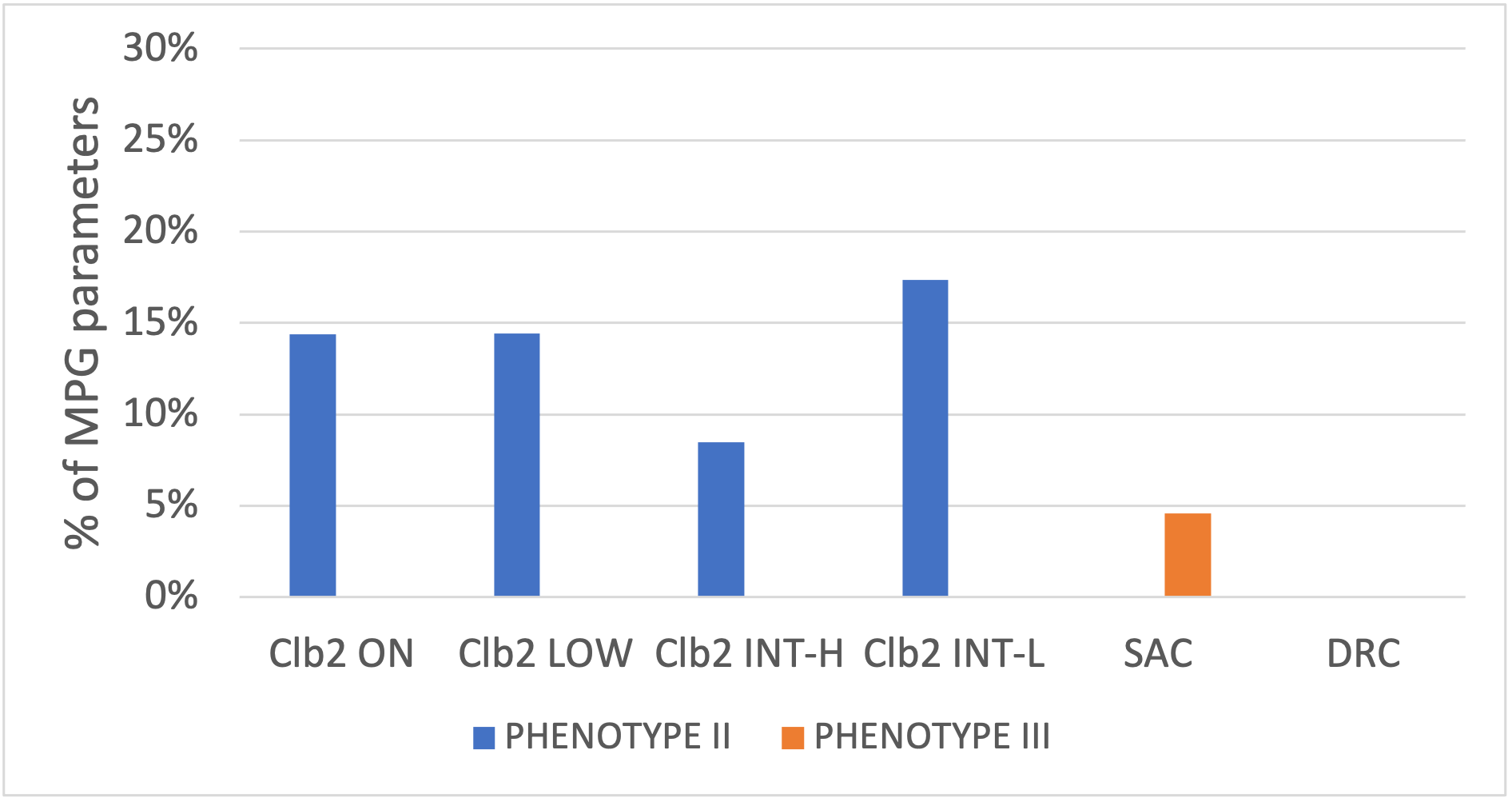} }}\hfil
    \subfloat[\centering Swi5-Yox1]{{\includegraphics[width=0.475\textwidth]{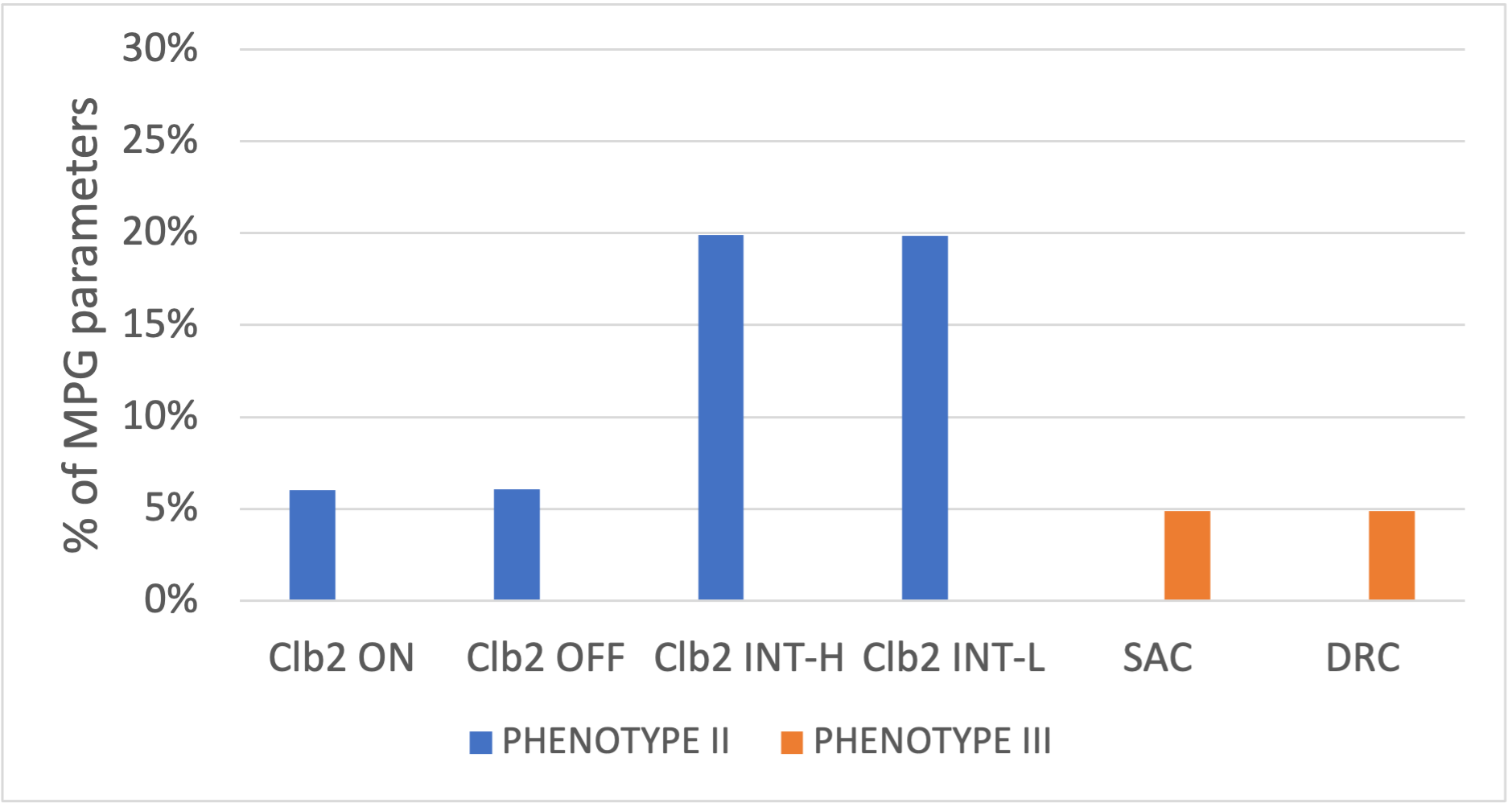} }}\par
    \subfloat[\centering Ace2-Nrm1]{{\includegraphics[width=0.475\textwidth]{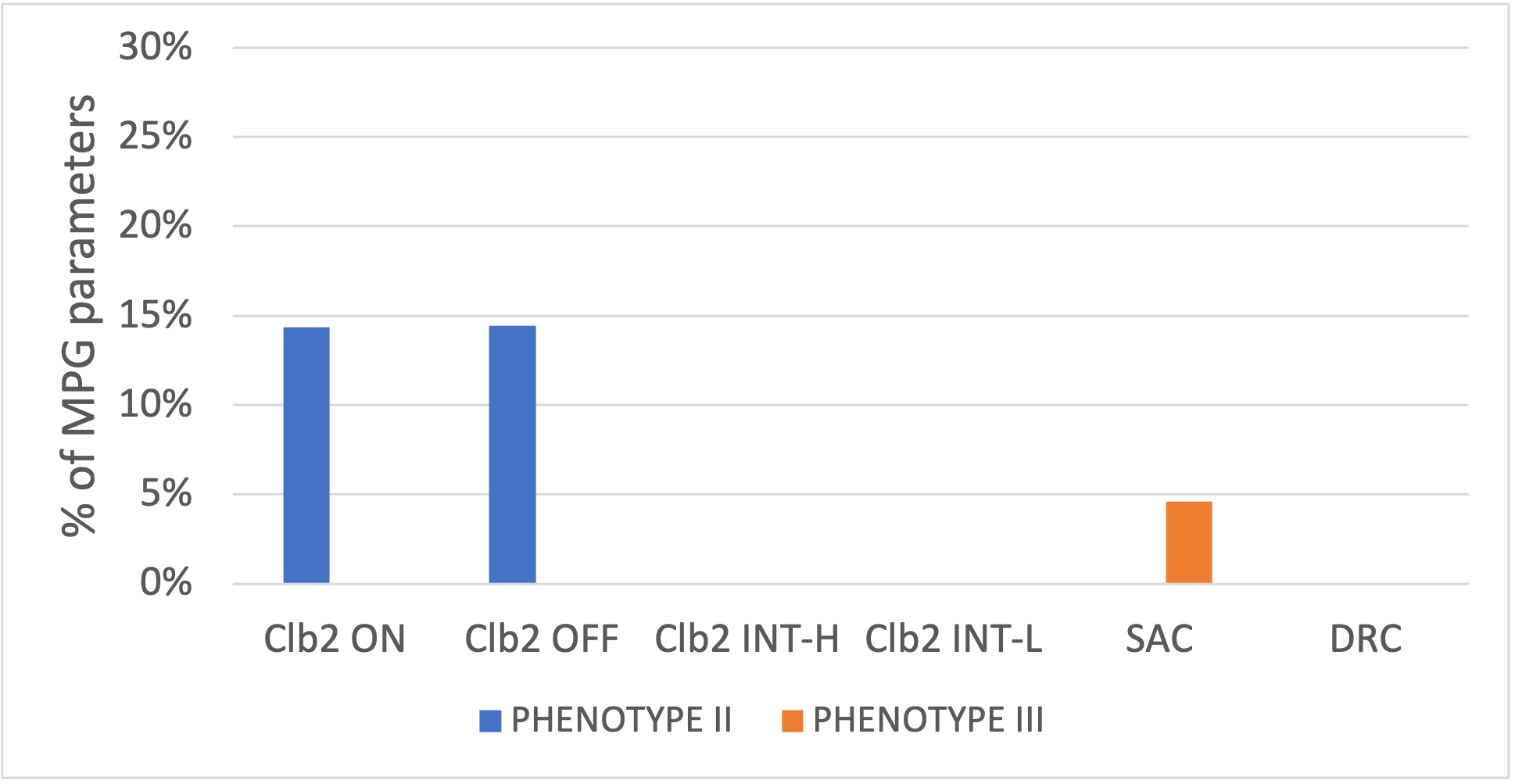} }} \hfil
    \subfloat[\centering Ace2-Yox1]{{\includegraphics[width=0.475\textwidth]{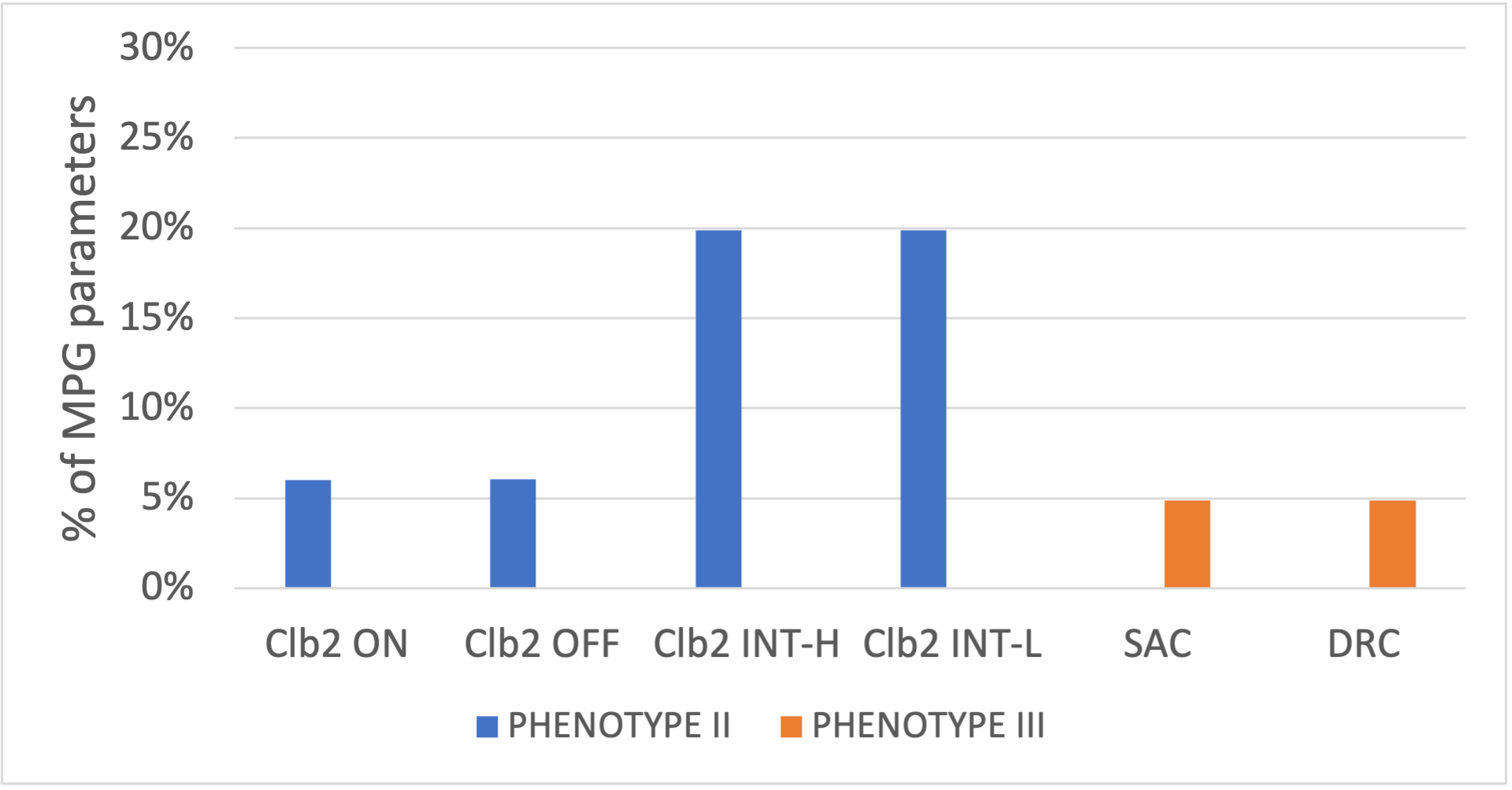} }}
    \caption{The percentage of WT cycling MPG parameters that also are predicted to show behavior consistent with a dataset associated to dynamical phenotype II or III. All percentages are nonzero except for the Clb2 INT-H and INT-L transcriptional phenotypes for the Ace2 proxy groups and the DRC transcriptional phenotype for the Nrm1 proxy groups.}
    \label{fig:phenoII-III}%
\end{figure}

We checked for the existence of dynamical phenotype II (mutant cycling) at Clb2 ON, Clb2 OFF, Clb2 INT-H, and Clb2 INT-L phenotype-permissible parameters; i.e., we checked for the co-existence of mutant and WT cycles at a single MPG parameter, which would indicate that a change in Clb2 alone does not disrupt mini pulse generator oscillations in the mini wavepool network. 
%We found that there existed pattern matches in all proxies for each mutant at phenotype-permissible parameters except for the Ace2-Nrm1 proxies for the intermediate mutant cycling phenotypes. All other proxy groups for the network model could be pattern-matched to all the mutant data at some .
% To complete dynamical phenotype II, we checked for  
% When pattern matched phenotype-permissible DSGRN parameters existed, the corresponding mini pulse generator parameter could also exhibit WT cycling at another Clb2 factor parameter; 
The percentages of MPG parameters where this occurs are shown in Figure~\ref{fig:phenoII-III}. The nonzero results indicate that a perturbation of Clb2 allows the mini wavepool to transition from WT cycling to every type of mutant cycling in three of the four proxy groups.

\begin{figure}
\centering
\begin{tabular}{cc}
\includegraphics[width=2.75in]{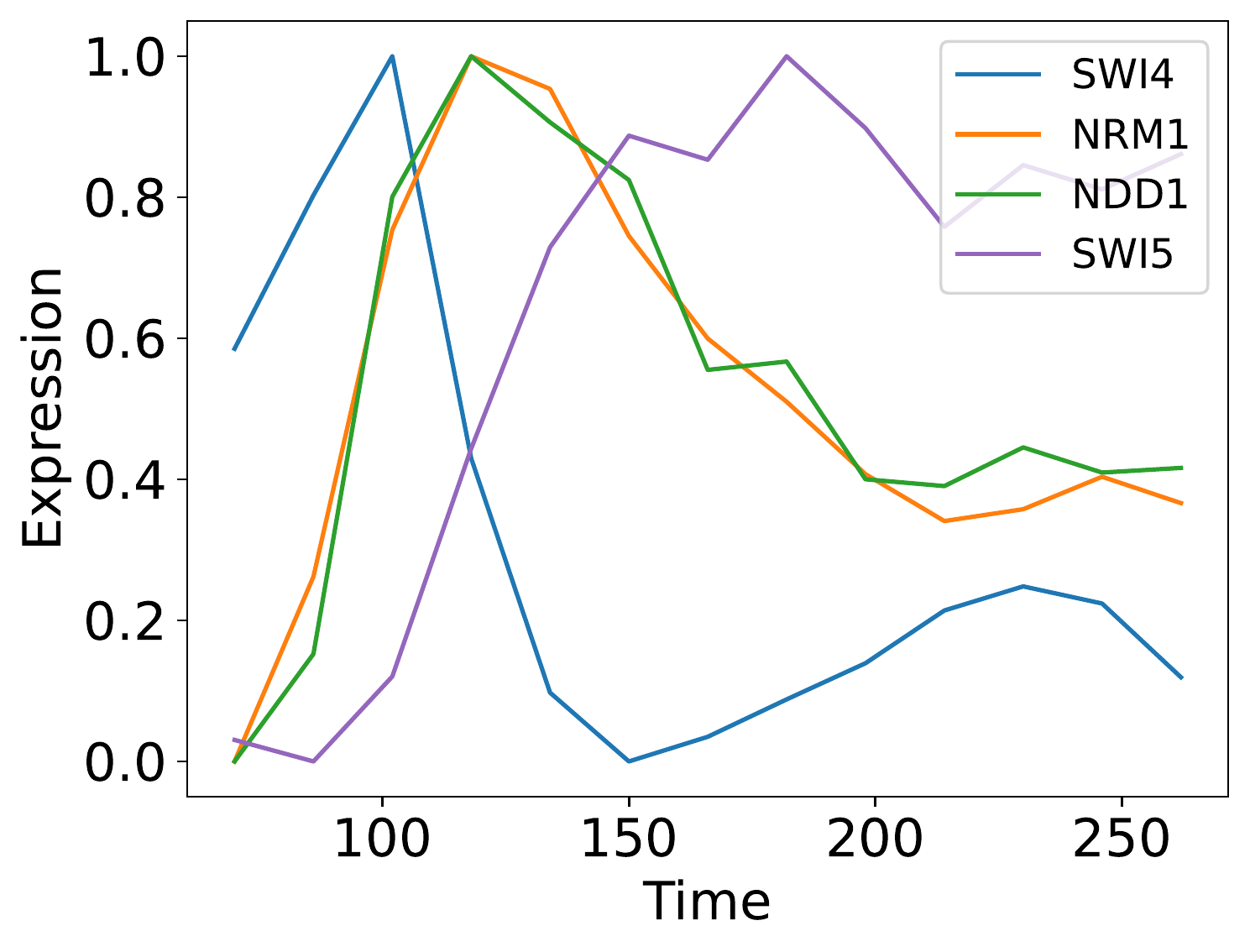} & \includegraphics[width=2.75in]{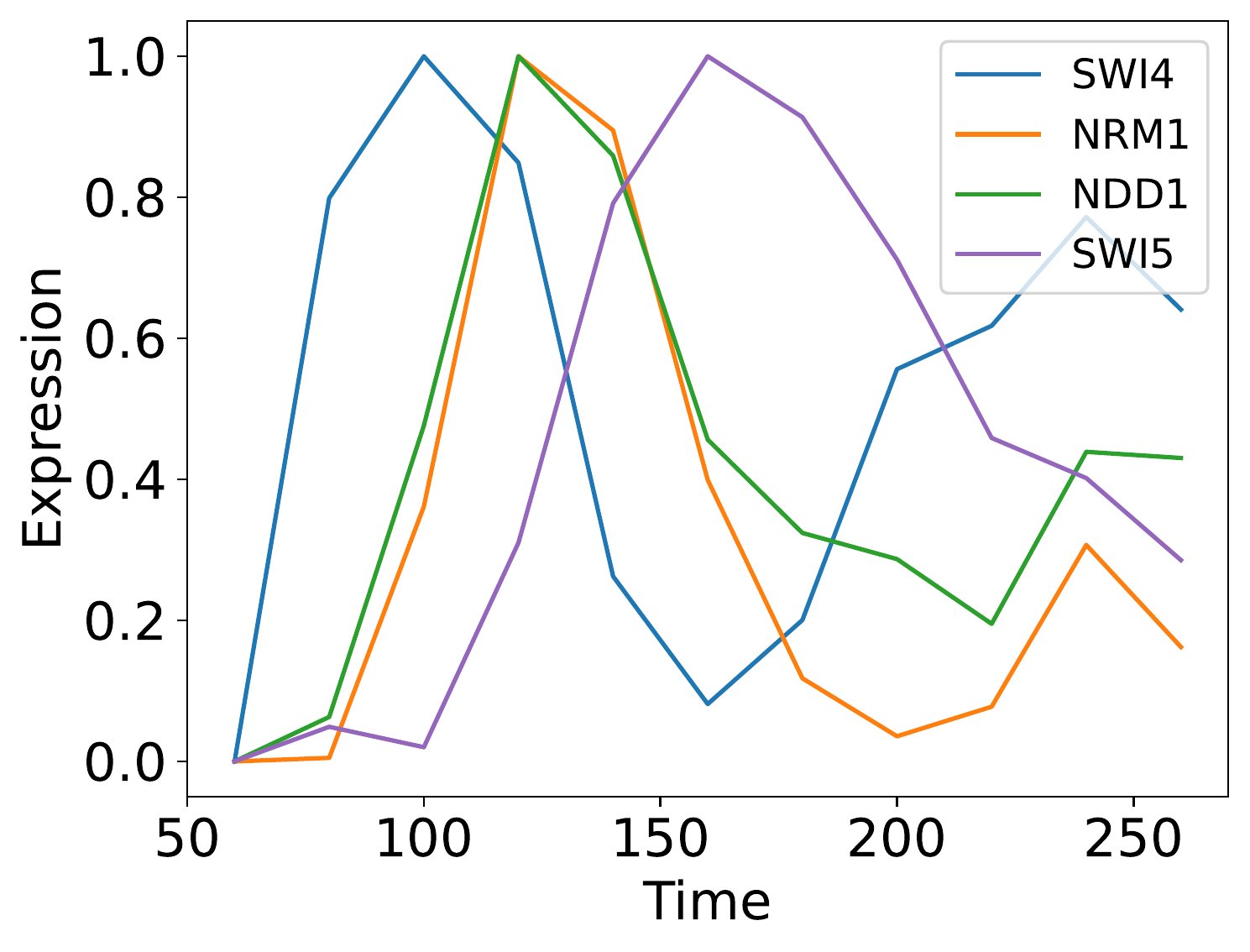} \\
(a) Clb2 ON & (b) Clb2 OFF \\
\includegraphics[width=2.75in]{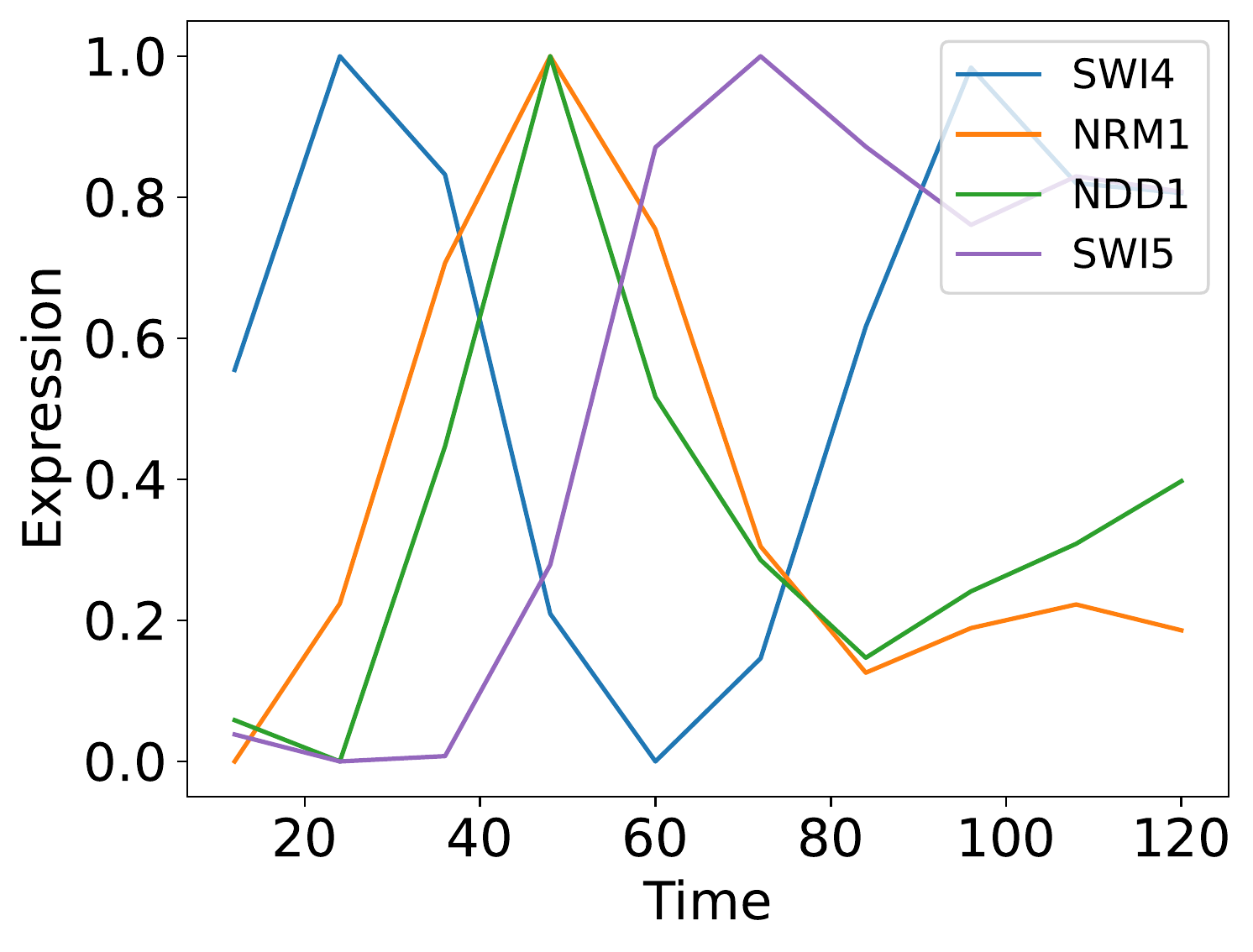} & \includegraphics[width=2.75in]{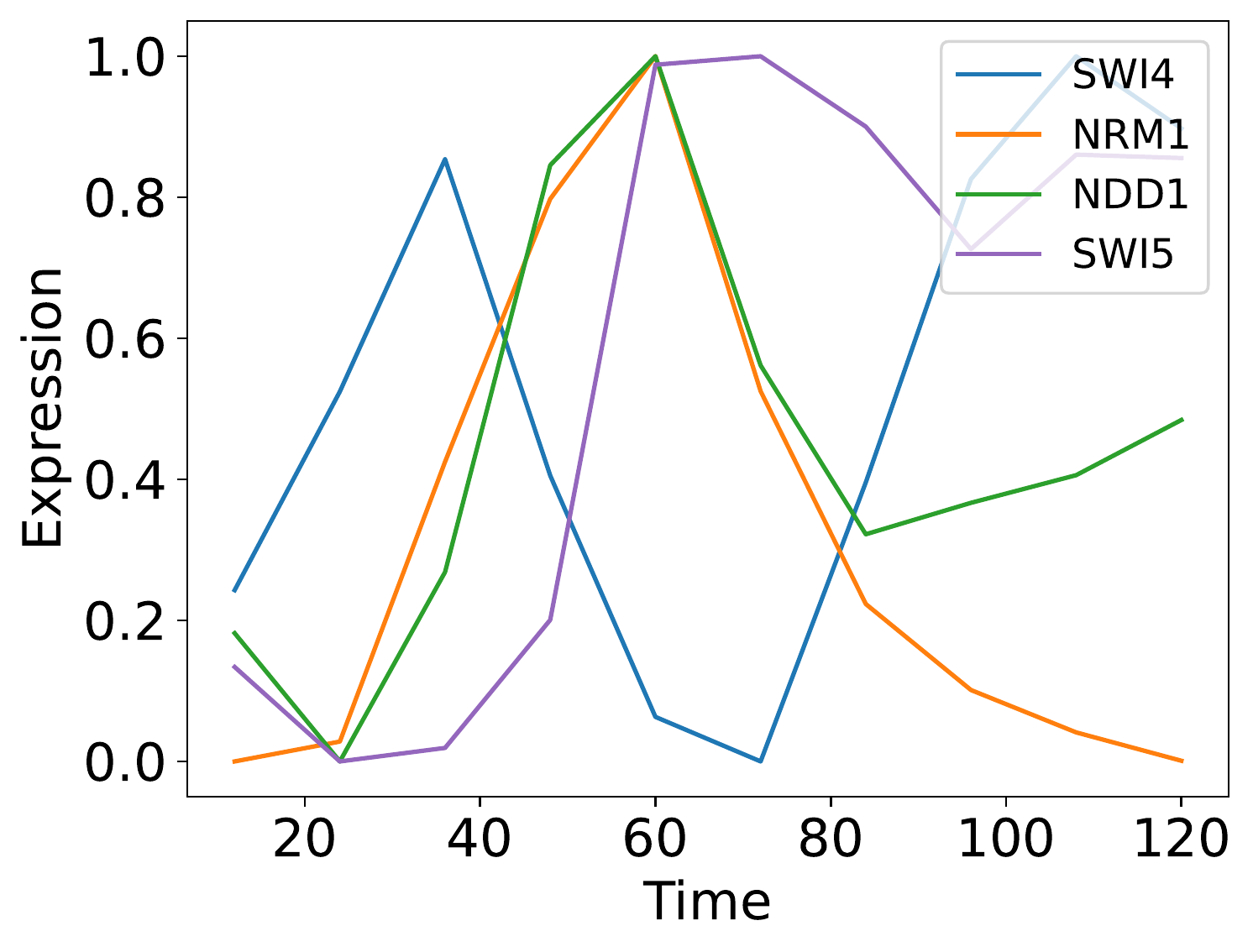} \\
(c) Clb2 INT-H & (d) Clb2 INT-L \\
 \includegraphics[width=2.75in]{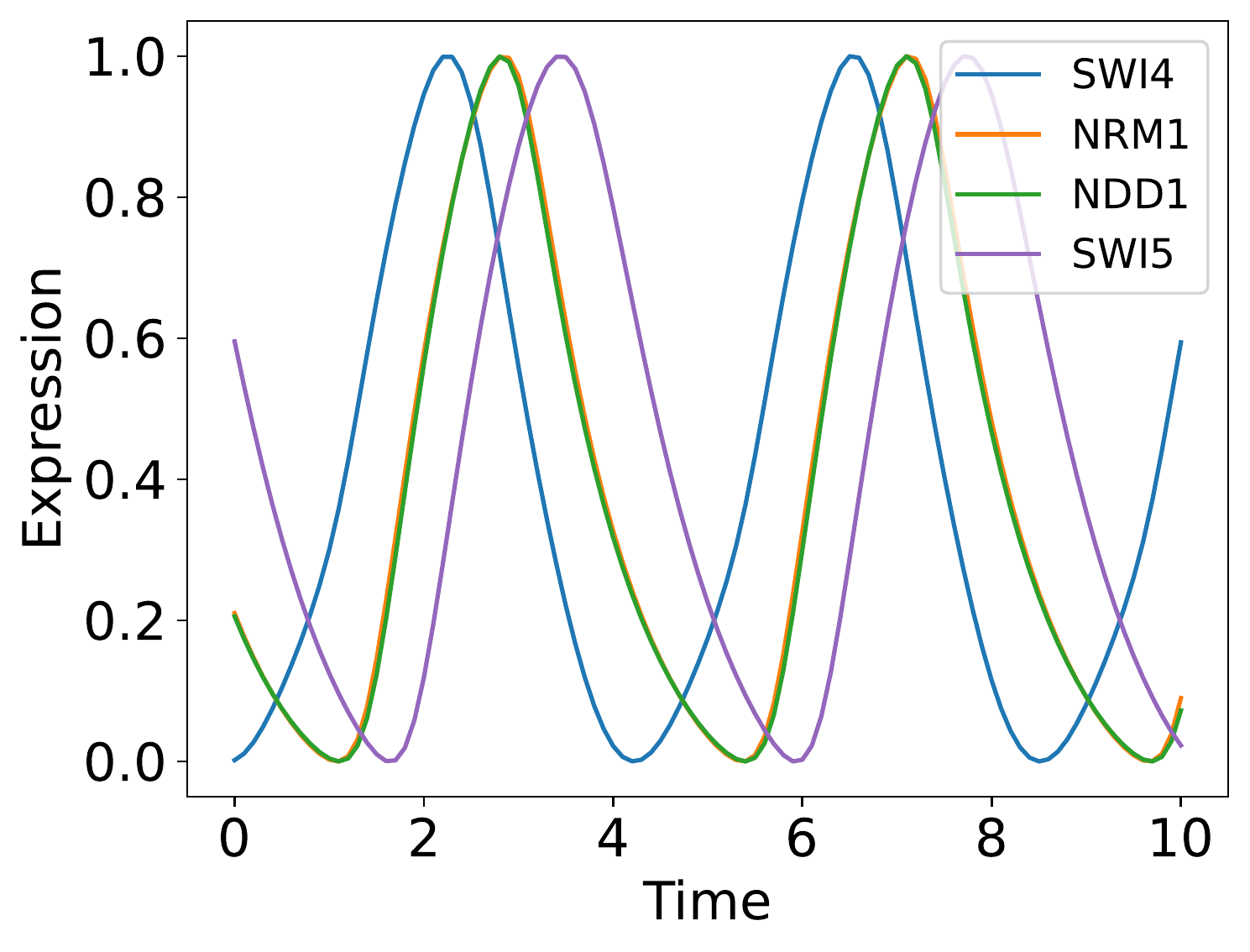} & \includegraphics[width=2.75in]{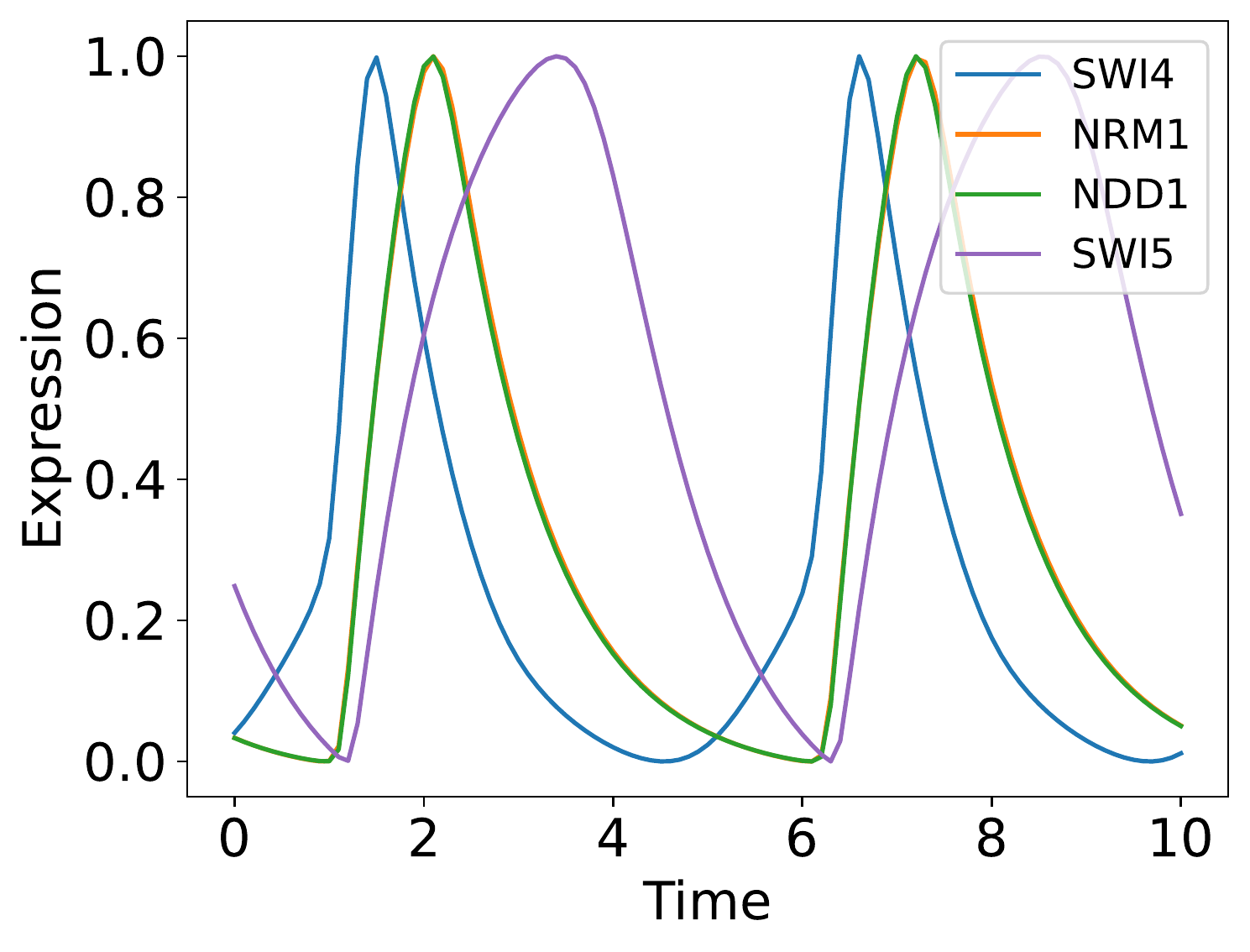} \\
 (e) Hill model 1 & (f) Hill model 2
\end{tabular}
\caption{\textbf{Hill modeling for mutant cycling in the mini pulse generator.} (a)-(d) The four mutant datasets associated to phenotype II, normalized between 0 and 1. (e)-(f) Hill models parameterized to match the order of extrema in the datasets. Notice that all four datasets show the Swi4 peak preceding the Nrm1 and Ndd1 peaks (which are coincident), followed by the Swi5 peak. Therefore, all four datasets can be roughly matched by the same parameterization of the mini pulse generator, and two examples are given in (e)-(f) for very different MPG parameters. In both of these choices, the Nrm1 and Ndd1 curves are nearly coincident. As for the WT simulation, we made no attempt to match amplitude or period. }\label{fig:phenoII_hill_model}
\end{figure}

The four Clb2 mutant datasets are shown in panels (a)-(d) of Fig~\ref{fig:phenoII_hill_model} normalized between 0 and 1. It can be seen that, given noise, they all roughly have the same order of maxima and minima. Therefore, we can demonstrate the ability of DSGRN to reduce parameter space for a Hill model. We chose to sample mini pulse generator parameters for the Swi5-Nrm1 proxy group that are predicted by DSGRN to be able to exhibit all five types of WT and mutant cycling behavior for varying Clb2 factor parameters (see Section~\ref{sec:coexist}). Fig~\ref{fig:phenoII_hill_model} (e)-(f) show simulations for two different MPG parameters that match the expected order of maxima and minima. We remark that Clb2 expression also oscillated in these simulations (not shown). This case is not excluded by DSGRN, since DSGRN will not predict oscillations that are too small to impact downstream targets.

It may seem inconsistent that there are Hill models that can recapitulate all four mutant cycling datasets, and yet the percent of pattern matches for the Clb2 ON, OFF, INT-L, and INT-H in Fig~\ref{fig:phenoII-III} are not identical within each proxy group. This occurs because the DSGRN pattern matching methodology is sensitive to noise.

\subsection{Dynamical Phenotype III: Checkpoint Arrest}\label{sec:pheno_III}

% \begin{table}[h!]
% \centering
% \begin{tabular}{|lll|}
% \hline
% \multicolumn{3}{|c|}{DSGRN Phenotype III}                  \\ \hline
% Proxies                         & SAC & DRC  \\ \cline{2-3} 
% \multicolumn{1}{|l|}{Swi5-Nrm1} & 4.59\%     & 0.00\%        \\
% \multicolumn{1}{|l|}{Ace2-Nrm1} & 4.59\%     & 0.00\%          \\
% \multicolumn{1}{|l|}{Swi5-Yox1} & 4.88\%     & 4.88\%          \\
% \multicolumn{1}{|l|}{Ace2-Yox1} & 4.88\%     & 4.88\%         \\ \hline
% \end{tabular}
% \caption{dynamical phenotype III results for the SAC FP (identically the Yox1 proxy DRC FP) in Table~\ref{fig:sac_fps}. The DRC FP for the Nrm1 proxies was not recoverable anywhere in DSGRN parameter space.}
% \label{tab:pheno_III}
% \end{table}

In DSGRN phenotype III, we identified MPG parameters that were predicted to exhibit a WT pattern match at one Clb2 factor parameter and SAC or DRC FPs  (Table~\ref{tab:qualitative_fps}) at others. In Fig~\ref{fig:phenoII-III}, the ``SAC" bars correspond to the percentages of MPG parameters exhibiting WT cycling and the SAC FP, which is identically the DRC FP for the Yox1 proxies seen by looking at the ``DRC" bars. The fact that these percentages are nonzero indicates that the Yox1 proxies are consistent with both mutant datasets, which must be true given that the dynamical phenotype representation of the two datasets is identical. 
% The existence of mini pulse generator DSGRN parameters exhibiting both WT cycling and the SAC/DRC FP indicates that regulation through Clb2 alone is sufficient to initiate the SAC for this model with a Yox1 proxy group, even if other regulators may in fact also exist.

On the other hand, the Nrm1 proxies are only consistent with the \textit{cse4} mutant dataset representing the SAC phenotype. The DRC FP for the Nrm1 proxies was absent not only in MPG parameters with WT cycling, but also absent across all of DSGRN parameter space. This indicates that the mini wavepool topology contains insufficient regulatory interactions to recapitulate the DRC phenotype when the Nrm1 time trace is considered. 

\begin{figure}[h]
    \centering
    \includegraphics[width=3.5in]{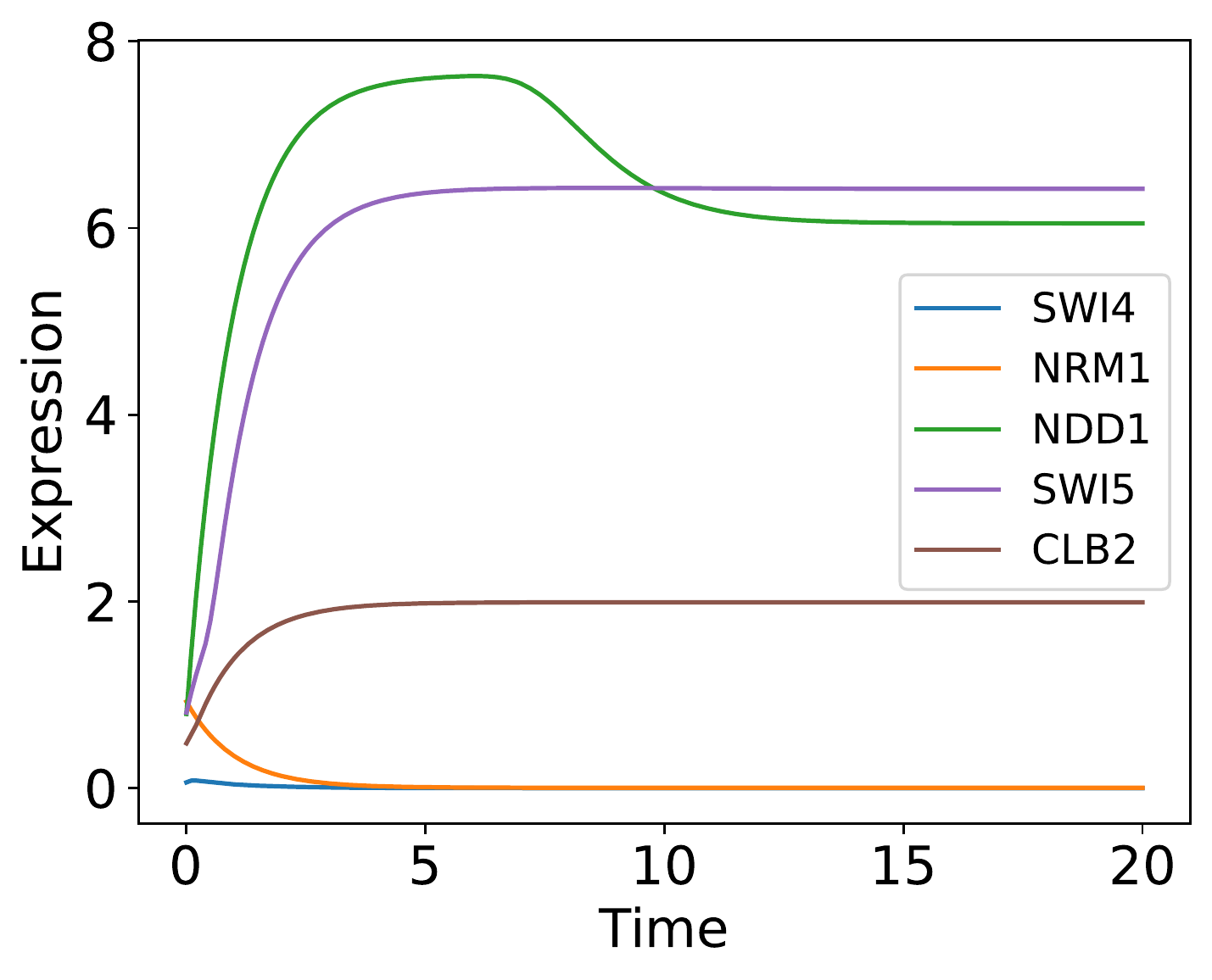}
    \caption{Hill model of the SAC, where Swi4 and Nrm1 are low, Ndd1 and Swi5 are high, and Clb2 is not low.}
    \label{fig:hill_model_sac}
\end{figure}

An example of a simulation from a Hill model parameterized with an MPG predicted to show both WT cycling and a SAC FP is shown in Fig~\ref{fig:hill_model_sac} for the Swi5-Nrm1 proxy group. It is easily seen by examination that the simulation traces are consistent with the SAC FP in Table~\ref{tab:qualitative_fps}.

\subsection{Co-existing Phenotypes}\label{sec:coexist}

\begin{figure}[!h]
    \centering
    \subfloat[\centering Swi5-Nrm1]{{\includegraphics[width=0.45\textwidth]{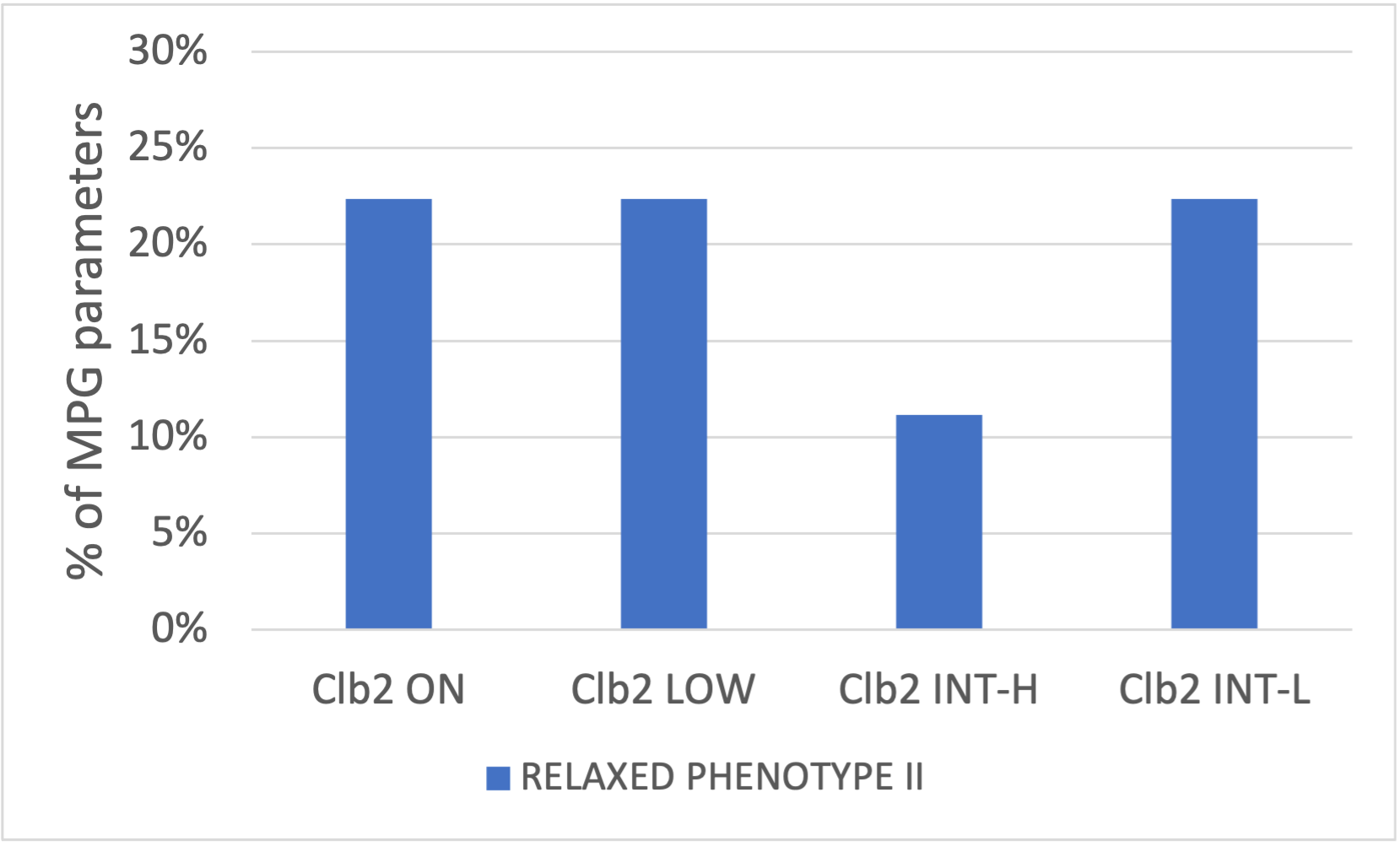} }}\hfil
    \subfloat[\centering Swi5-Yox1]{{\includegraphics[width=0.45\textwidth]{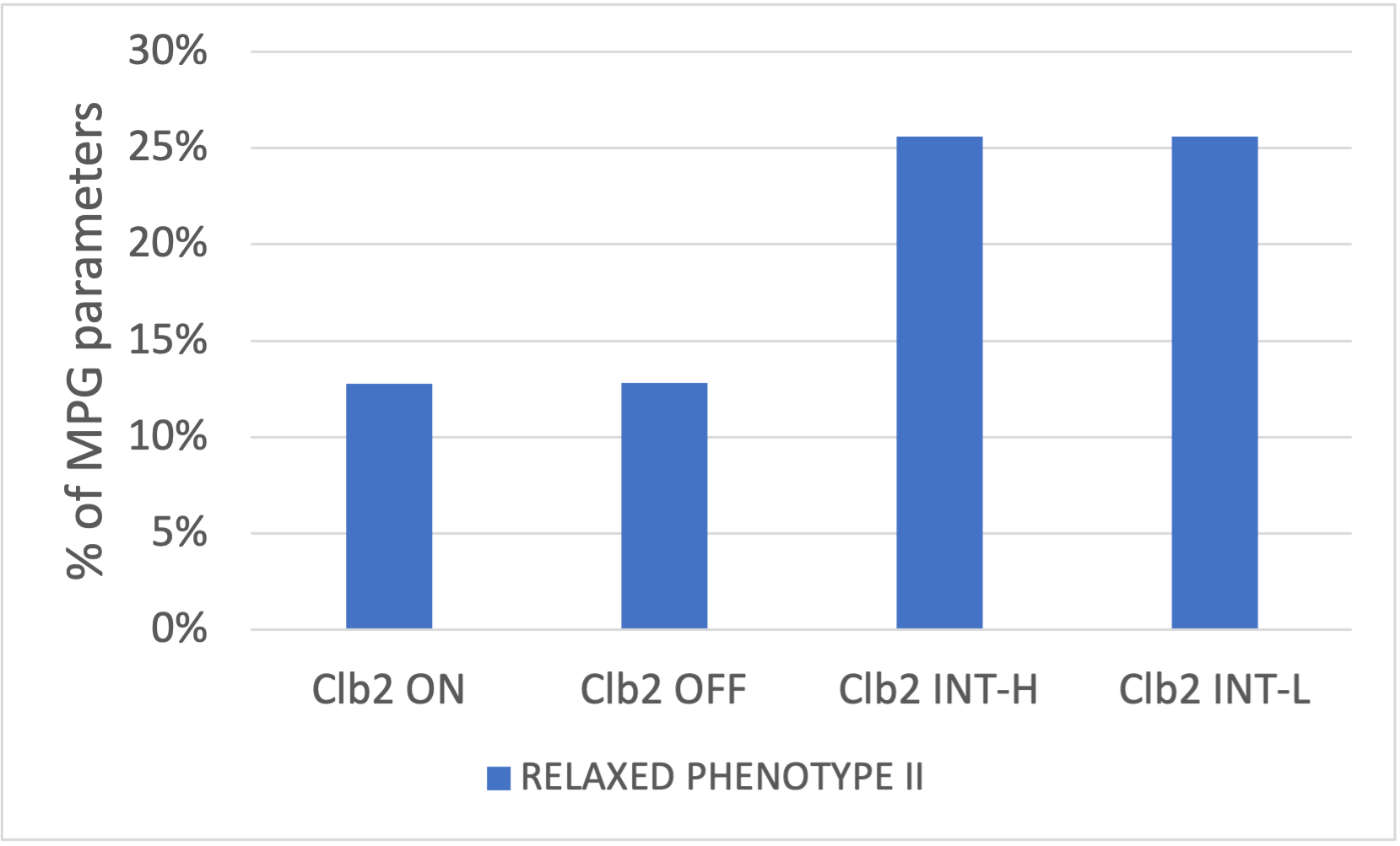} }}\par
    \subfloat[\centering Ace2-Nrm1]{{\includegraphics[width=0.45\textwidth]{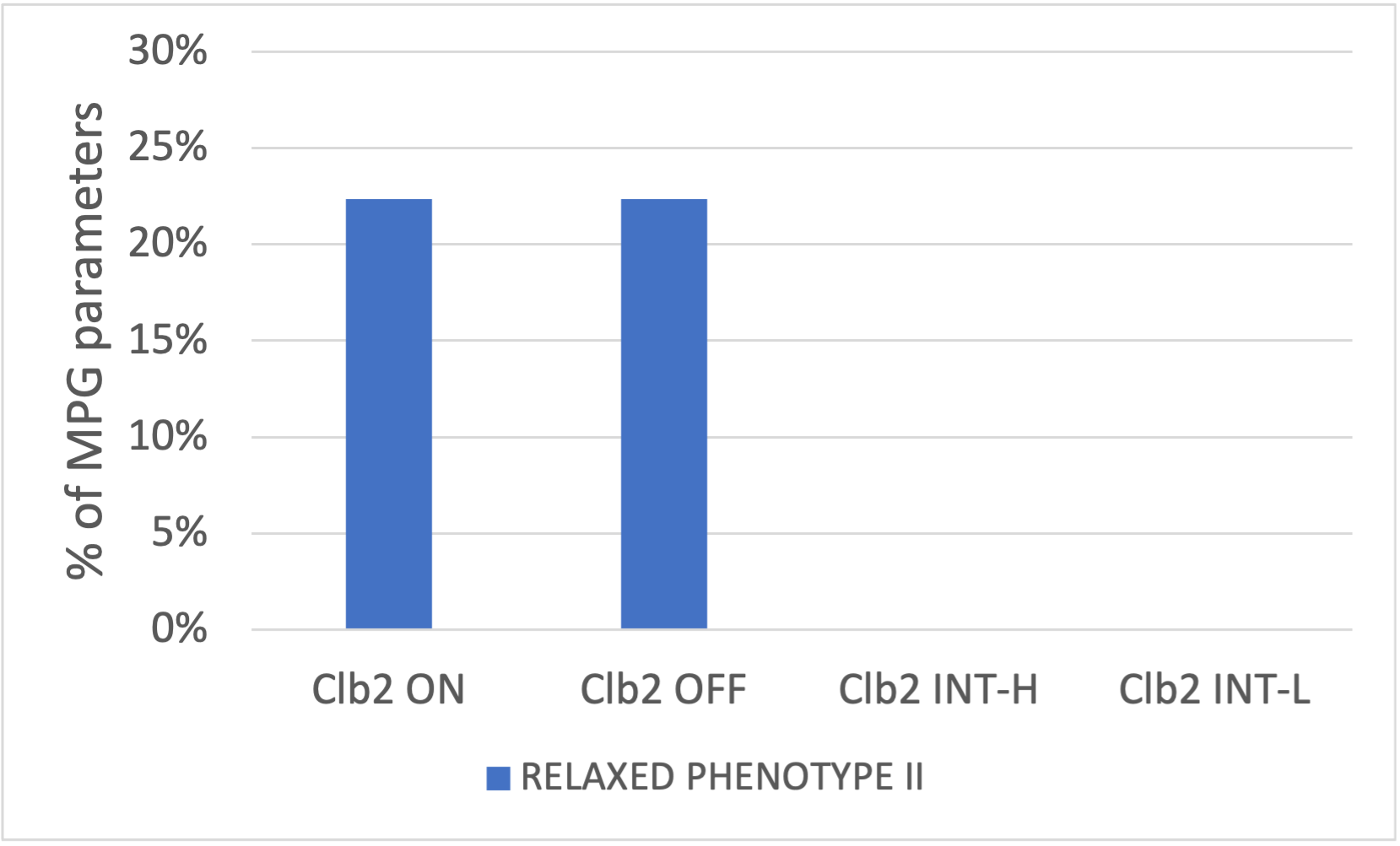} }} \hfil
    \subfloat[\centering Ace2-Yox1]{{\includegraphics[width=0.45\textwidth]{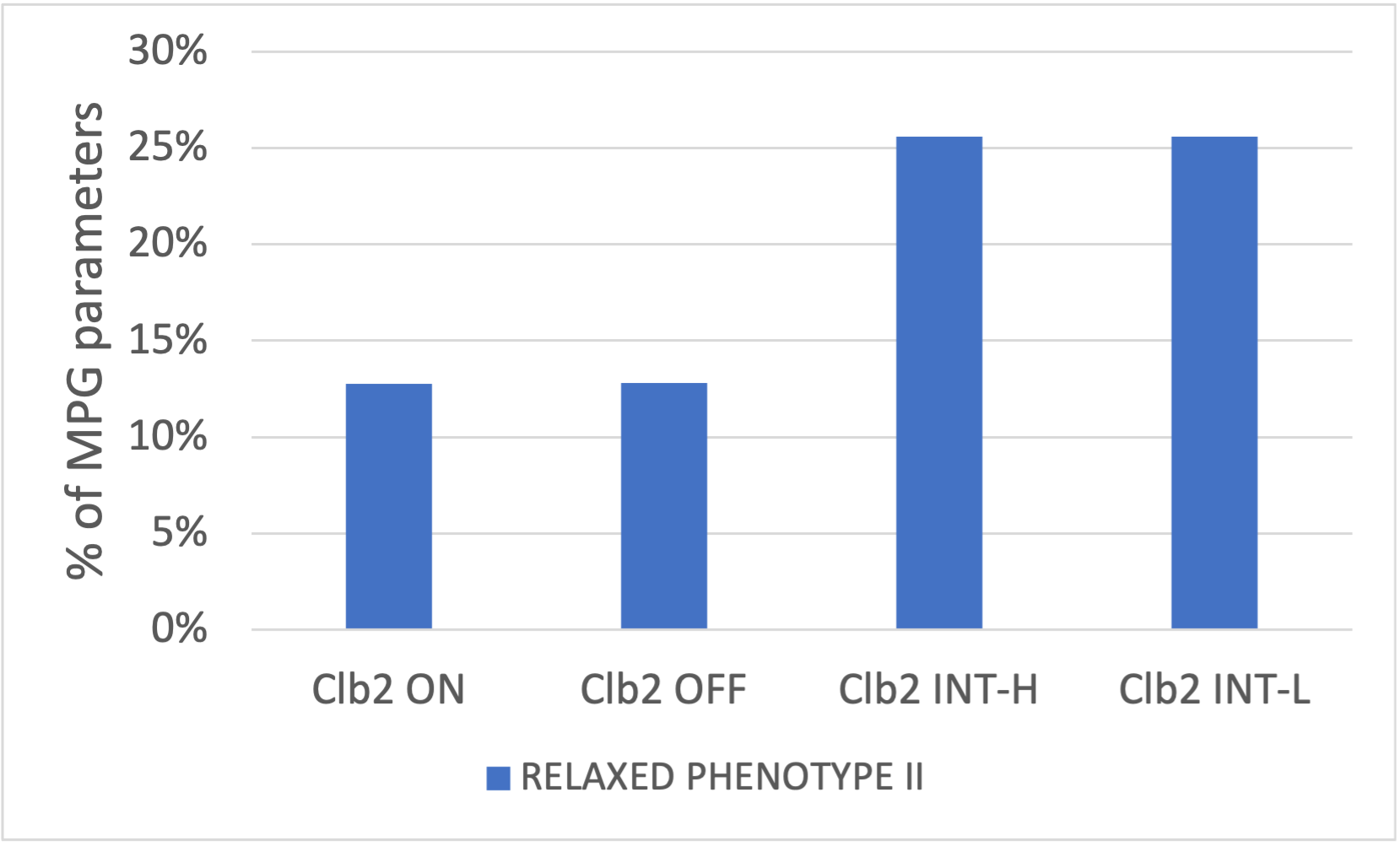} }}
    \caption{The percentage of WT MPG parameters that also are predicted to show behavior
consistent with a dataset associated to dynamical phenotype II after relaxing the phenotype-permissible restriction. Similar to the results from dynamical phenotype II with the phenotype-permissible restriction, all percentages are nonzero except
for the Clb2 INT-H and INT-L transcriptional phenotypes for the Ace2 proxy groups.}
    \label{fig:phenoII_co}%
\end{figure}

We looked for mini pulse generator parameters that could support all four mutant cycling phenotypes in addition to WT cycling at phenotype-permissible parameters. There were none such, indicating that either an internal change in or an external regulation of the mini pulse generator occurs between mutant phenotypes. 
This result is dependent on the modeling choice of using only phenotype-permissible parameters associated to the Clb2 mutant cycling datasets. If this restriction is relaxed, the percentage of mini pulse generator parameters at which there is both WT cycling and a single mutant cycling dataset increases (compare Fig~\ref{fig:phenoII_co} to Fig~\ref{fig:phenoII-III}). More precisely, relaxing the phenotype-permissible restriction means that mutant cycling phenotypes are searched over all values of the Clb2 factor parameter, not just the 10\% of DSGRN parameter nodes originally specified.
% Moreover, we are able to find thousands of mini pulse generator parameters that support WT cycling and all four types of mutant cycling. \Julian{not sure if this can be said for every dataset}\Bree{we know it must be true for the Swi5-Nrm1 proxy set. Do you know for the other proxy sets?}

However, relaxing the phenotype-permissible modeling restriction still does not allow a mini pulse generator parameter that can show matches to every dataset in all three dynamical phenotypes, including checkpoint arrest. The closest achiever was the Swi5-Nrm1 proxy set for which the mini wavepool model was able to exhibit the WT cycling phenotype, the four mutant cycling phenotypes, and the SAC phenotype at 3575 mini pulse generator parameters. 

\section{Discussion}

We have demonstrated that a small network model (mini wavepool, Fig~\ref{fig:wavepool} (Right)) approximating the pulse generation capacity of the yeast cell cycle transcriptional oscillator (Fig~\ref{fig:wavepool} (Left)) is capable of matching multiple datasets with different experimental perturbations. This process provides further evidence for the validity of a pulse generator mechanism with CDK control in the yeast cell cycle and highlights the flexibility of dynamical behavior that a single network is capable of producing. We conclude that the mini wavepool network is \textit{controllable}, in the sense that a large number of dynamical phenotypes are accessible via parameter perturbation, and that the region of wild-type behavior is relatively small within parameter space. In contrast, a network might instead be \textit{robust}, wherein wild-type behavior is prevalent across parameter space, making it difficult to push network behavior into a different dynamical regime.

We devised three computational phenotypes using the software DSGRN~\cite{cummins:2016,DSGRN_repo} that associate to seven different datasets (Table~\ref{tab:dsgrn_pheno}). Dynamical phenotype I corresponds to wild-type cell cycle behavior under standard conditions. Dynamical phenotype II encompasses mutant cycling---oscillations of the mini pulse generator under fixed values of the CDK Clb2 induced by various knockout experiments. In dynamical phenotype III, checkpoint behavior was evaluated by computationally seeking qualitative equilibrium values determined by each of two datasets, one representing the spindle assembly checkpoint and one representing the DNA replication checkpoint. 

We interpret model consistency with the data in the following way. If the mini wavepool model is consistent with wild-type cycling, then the model faithfully captures the behavior of 
the undisrupted cell cycle. If it is consistent with mutant cycling, then the model is capable of reproducing observed oscillatory behavior in the mini pulse generator regardless of the state of Clb2. If it is consistent with checkpoint behavior, then it provides the hypothesis that regulation through Clb2 is sufficient to induce the checkpoint in the mini pulse generator, even if in reality other regulators are implicated in cell cycle control. If the mini wavepool model is consistent with all of these dynamical phenotypes, then it is a complete hypothesis for explaining the seven datasets, while the non-existence of any mutant phenotype indicates network model incompleteness.

% Checks for consistency occur at each phenotype-permissible DSGRN parameter, of which there are millions, and the percentage of parameters exhibiting a dynamical phenotype indicates how prevalent the behavior is. 

We found that 3.7-4.3\% of potential wild-type DSGRN parameters predicted dynamical behavior consistent with the WT dataset, indicating that the network model is capable of reproducing standard cellular conditions. Of these WT pattern-matched parameters,  up to 20.0\% were consistent with various mutant datasets at various proxy groups, although not all phenotypes were exhibited at every proxy group. 
There were no mini pulse generator parameters that could support WT cycling and all four types of mutant cycling within a single proxy group at phenotype-permissible DSGRN parameters, i.e. fixed Clb2 levels hypothesized to be associated to particular datasets. However, if this modeling choice is not enforced, mini pulse generator parameters were located that support WT and all four mutant cycling phenotypes, as well as (in some cases) SAC arrest.
Relaxing phenotype-permissibility in the mutant cycling phenotype is akin to acknowledging that Clb2 expression in the datasets may not be sufficiently close to constant for our model assumptions to hold. If we allow for this possibility, then there exists a highly narrowed selection of mini pulse generator parameters that can recapitulate the dynamics seen in multiple datasets. One interesting and open problem is whether this selection of parameters is most likely to contain biologically reasonable parameterizations for the mini pulse generator, since only Clb2 parameter modifications are necessary to induce phenotypic changes. This is an area of future research.

% \Bree{Discuss the sensitivity of the method to noise and time series length leading to the discrepancy between a single Hill model for all four mutant cycling datasets, but different percentages of DSGRN pattern matches.}

% % While the restriction of phenotype-permissibility is relaxed we still require that the mini-pulse generator for each DSGRN parameter is pattern matched to a biological dataset.
% \Julian{does removing the restriction on phenotype-permissible parameters still allow for determination of biologically relevant mini-pulse generator remainder parameters?} 
% It is useful to narrow parameter space by determining Clb2-divergent groups exhibiting multiple phenotypes of interest as way to determine biologically relevant parameters. This narrowing of parameter space also acts to validate the network and identifies the affect that Clb2 modulation has on the dynamics of the network. After identifying Clb2-divergent groups in the Swi5-Nrm1 proxy for all cycling phenotypes and additionally for the SAC FP, we can see that Clb2 modulation is sufficient to initiate the SAC FP.

We observed that when using Nrm1 proxy sets, the DRC phenotype was not supported at any DSGRN parameter in the network model. This indicates that the model lacks important regulatory elements that are necessary for the DRC when considering Nrm1 as a proxy in the network. We hypothesize that distinguishing  between SBF and MBF in the mini wavepool may help rectify this issue. Current models like the network seen in \cite{cho:2017} define MBF and SBF as the same node, yet results from \cite{oliveira:2012,travesa:2012} indicate a mechanism for DRC arrest dependent solely on MBF activity. Given replication stress, the protein Rad53p inactivates the MBF co-repressor Nrm1. This activation of MBF induces up-regulation of G1/S genes within S-phase. As a result of this MBF pathway being activated, the DRC is initiated. We suggest that enlarging the mini wavepool to include MBF and Rad53p may permit consistency with the DRC phenotype.
%While the mechanism regulating the DRC may still be unclear, it is apparent that this mechanism has crosstalk with the TF network through MBF and/or Nrm1. Thus, characterizing the difference between SBF and MBF within our hypothetical network model could lead to a more valid model. 

There is another network enhancement that our work suggests. We showed that SAC arrest is supported in the network model by only modifying the parameterization of Clb2, suggesting the presence of a regulatory element controlled by the SAC mechanism that impinges solely on Clb2. Research has shown that one of the physiological goals of the SAC is to inhibit APC/Cdc20 activity which in turn inhibits Clb2 degradation \cite{Nilsson:2008}.
% The wavepool network ideally supports the SAC and DRC, which are known to regulate cell cycle progression across different organisms \cite{Marston:2017,Musacchio:2007,Wang:2017,Nilsson:2008}. In \textit{S. cerevisiae}, these mechanisms are believed to be an interaction with the transcriptional program to arrest the cell at a specific point in the cell-cycle progression \cite{bristow:2014}. 
In \cite{Wang:2017} it was seen that the rapid degradation of Cdc20 is necessary for SAC arrest. Using ordinary nonlinear differential equations and spatial simulations, the authors of \cite{ibrahim:2015} were able to identify a SAC mechanism that acts to fully sequester Cdc20 and inhibit APC activity. It was also suggested within \cite{Wang:2017} that there is potential for indirect regulation of Cdc20 through substrates of APC, such as Clb2 and Clb5. By expanding the mini wavepool to include feedback with APC/Cdc20, we can test the hypothesis that APC and Cdc20 interaction can initiate the Clb2 parameter change required to transition from WT cycling to SAC.
%Together these findings suggest an interaction between the SAC and the transcriptional program which drives the arrest of the cell-cycle.
% It has been seen in the biological research that APC/Cdc20 is critical for initiating the metaphase-to-anaphase transition \cite{Wang:2017,Nilsson:2008,Yang:2015}. 
% Research has also shown that the physiological goal of the SAC is to inhibit APC/Cdc20 activity which in turn inhibits Clb2 degradation \cite{Nilsson:2008}. Thus, without Clb2 degradation the cell does not exit mitosis. We wanted to test whether SAC arrest could occur through control of Clb2, indicating a node not present in the mini wavepool that can facilitate this perturbation to Clb2.  Clb2-divergence from a cycling WT DSGRN parameter to the SAC FP DSGRN parameter suggests that entry into the SAC could occur through additional regulation of Clb2 not characterized by the network. Given that we know Cdc20 is an important factor in reaching SAC arrest, the addition of a Cdc20 node which interacts with Clb2 might be justified. If Cdc20 does not add more validity to the network model, adding in an APC node which interacts with Clb2 could be another option for expanding the model. 

In its original implementation, DSGRN does not differentiate between transcriptional regulation vs protein modifications such as phosphorylation or ubiquitination, although these generalizations are newly available~\cite{Cummins21}, as well as modeling transport across the nuclear membrane~\cite{fox2022modeling}. In this manuscript the original DSGRN package was used, and as a consequence transcriptional interactions and post-transcriptional modifications were not distinguished in the modeling of the mini wavepool in Fig~\ref{fig:wavepool} (Right). This means that we used mRNA expression levels of Clb2 as a substitute for protein levels and activity.  Protein and RNA levels are certainly correlated at least under some conditions \cite{RN1936}.
An obvious next step is to expand the mini wavepool to model post-transcriptional modifications, as well as to incorporate other regulators. We are currently in the midst of an effort to explore the endocycling phenomenon in a larger network using the new methodology from~\cite{Cummins21}.

% \Julian{Do we want to touch on these ideas?:
% Would the introduction of APC possibly allow for bistability between cycling and a FP? There exists a signal to the network which relies on connectivity to Clb2 that allows for initiation of the SAC FP.
% } \Bree{I think we should take out any discussion of multistability. I think it muddies the message.}

\section{Methods}\label{sec:methods}

In this section, we discuss the basic properties of DSGRN~\cite{cummins:2018,cummins:2016,gedeon:2018,DSGRN_repo} that are used when interpreting WT and mutant transcriptional phenotypes as DSGRN dynamical phenotypes. We discuss in detail how the dynamical phenotypes are constructed and interpreted with respect to the data.

DSGRN comprehensively computes coarse features of the dynamical behavior of a genetic regulatory network over a combinatorial representation of parameter space that is finite \cite{cummins:2016}. These coarse features include oscillatory behavior with stereotyped orders of maximum and minimum concentrations of gene products, and the number and type of equilibria. DSGRN uses techniques from ordinary differential equations and graph theory to compute these behaviors.

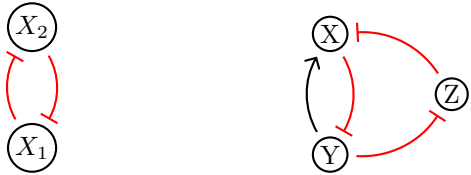
\begin{figure}[h!]
\centering
\begin{tabular}{cc}
\begin{tikzpicture}[main node/.style={circle, draw, thick, inner sep=1pt, minimum size=0pt},scale=0.8]
\node[main node] (n1) at (0,0) {$X_1$};
\node[main node] (n2) at (0,2) {$X_2$};

\path[->,>=angle 90,thick]
(n1) edge[-|,red,shorten <= 3pt, shorten >= 3pt, bend left] node[] {} (n2)
(n2) edge[-|,red, shorten <= 3pt, shorten >= 3pt,bend left] node[] {} (n1)

;                       
\end{tikzpicture}

& 

\hspace{1in}

\begin{tikzpicture}[main node/.style={circle, draw, thick, inner sep=1pt, minimum size=0pt},scale=0.8]
\node[main node] (n1) at (0,2) {X};
\node[main node] (n2) at (0,0) {Y};
\node[main node] (n3) at (2,1) {Z};

\path[->,>=angle 90,thick]
(n1) edge[-|,red,shorten <= 3pt, shorten >= 3pt, bend left] node[] {} (n2)
(n2) edge[shorten <= 3pt, shorten >= 3pt,bend left] node[] {} (n1)
(n2) edge[-|,red,shorten <= 3pt, shorten >= 3pt,bend right] node[] {} (n3)
(n3) edge[-|,red,shorten <= 3pt, shorten >= 3pt,bend right] node[] {} (n1)
;                       
\end{tikzpicture}

\end{tabular}

\caption{Left: The two-node toggle switch where $X_1$ represses $X_2$ and $X_2$ represses $X_1$. Right: An example three-node network, where $X$ represses $Y$, $Y$ activates $X$ and represses $Z$, and $Z$ represses $X$.}
\label{fig:example-networks}
\end{figure}

Four kinds of graphs provide a framework for understanding DSGRN. A gene regulatory network (GRN) is the input to DSGRN and an undirected parameter graph (PG) is the basic structure of DSGRN. The output of DSGRN is a collection of directed state transition graphs (STGs) and their corresponding directed Morse graphs (MGs), one for each node in the parameter graph. As defined earlier a GRN is a system of interacting gene products that inhibit or activate one another. The GRN modeled in this paper is in Fig~\ref{fig:wavepool} (Right), which is built as a directed graph.  As a running example for explaining DSGRN computations in this section we will be talking about the simpler networks in Fig \ref{fig:example-networks}. On the left, node $X_1$ inhibits node $X_2$ and $X_2$ inhibits $X_1$, otherwise known as the ``toggle switch'' \cite{gedeon:2018}. On the right, we see a three-node network that exhibits more complex dynamical behavior than the toggle switch, but remains amenable to manual construction of the various DSGRN graphs.

% \sout{Due to the structure of DSGRN there is a limitation in how the network interactions are defined. In its current state DSGRN does not differentiate between different types of regulatory interactions, such as TFs vs protein modifications such as phosphorylation or ubiquitination. Therefore, the difference between transcriptional interactions and post-transcriptional modifications are ignored for the purpose of investigating the mini wavepool; i.e. activation and repression are the only differentiators between edges.} 

\subsection{Switching systems and DSGRN parameters}\label{sec:dsgrn}

DSGRN is, at its heart, a rigorous connection between Boolean modeling and ordinary differential equation (ODE) modeling. It can equivalently be viewed as an asynchronous multi-level Boolean modeling approach~\cite{crawford-kahrl_joint_2022} or as a dynamical systems approach~\cite{cummins:2016}. In this exposition, we will introduce DSGRN using the formalism of dynamical systems.

A GRN can be modeled using a system of discontinuous ODEs called a switching system~\cite{Thomas1991, veflingstad07,deJong2002,edwards00,Ironi2011}, where each node in the network is modeled by an ODE of the form:

\begin{equation}
    \dot{x}_i = -\gamma_i x_i + \lambda_i(x) , \quad  i = 1,...,n,
\end{equation}
where $x =(x_1,x_2,\dots,x_n)$.
For every node in the network, $x_i$ is the concentration of species $i$, $\gamma_i$ is the decay rate of species $i$, and $\lambda_i(x)$ is a product of sums of step functions,  one for each regulating edge $j$ on node $i$. The step functions $\sigma^\pm_{i,j}(x_i)$ are given by:

\begin{align}
    if \quad j \rightarrow i \quad then \quad \sigma_{i,j}^+(x_j) &= \begin{cases}
    l_{i,j} \quad if \quad x_j < \theta_{i,j} \\
    h_{i,j} \quad if \quad x_j > \theta_{i,j}\\
    \end{cases} \label{eq:positive_interaction}\\
    if \quad j \dashv i \quad then \quad \sigma_{i,j}^-(x_j) &= \begin{cases}
    l_{i,j} \quad if \quad x_j > \theta_{i,j} \\
    h_{i,j} \quad if \quad x_j < \theta_{i,j}\\
    \end{cases} \label{eq:negative_interaction}
    \end{align}
representing activation and inhibition respectively, where $0<l_{i,j} < h_{i,j}$ are low and high constant values. For an activating step function $\sigma_{i,j}^+(x_j)$, an increasing concentration of $x_j$ causes an increasing rate of production of $x_i$ as $x_j$ crosses the threshold value $\theta_{i,j}$. Similarly, a repressing step function $\sigma_{i,j}^-(x_j)$ indicates a decreasing rate of production of $x_i$ as $x_j$ increases across the threshold $\theta_{i,j}$. In the case of the toggle switch we only see a single inhibiting edge per node so both nodes are defined by a single inhibiting step function, i.e. $\lambda_1 = \sigma^-(x_2)$ and $\lambda_2 = \sigma^-(x_1)$.

Suppose that node $j$ regulates both $i$ and $k$. Then we require that $\theta_{k,j} \neq \theta_{i,j}$ so that the effect of $x_j$ on each of its downstream regulatory targets is distinct and therefore totally ordered. For example, the thresholds on the two out-edges from $Y$ in Fig~\ref{fig:example-networks} (Right) are required to be different; however there are no requirements on the relationship between any other pair of thresholds.

A DSGRN parameter is a collection of inequalities governing the relationship of the low, high, and threshold values for each node within the network. Each DSGRN parameter consists of two parts for each node in the network: a logic parameter and a order parameter~\cite{cummins:2016}. A key observation is that the logic and order parameters for a node are independent of all other nodes in the network, and therefore may be chosen independently. 
An order parameter defines the order of the threshold values for a node. For example, the $Y$ node in Fig \ref{fig:example-networks} (Right) has two threshold values due to its two out-edges, $\theta_{X,Y}$ and $\theta_{Z,Y}$. There are then two order parameters for Y: $\theta_{X,Y} < \theta_{Z,Y}$ and $\theta_{X,Y} > \theta_{Z,Y}$. All other nodes in our examples have a single out-edge and are trivially ordered. 

The number of thresholds associated to each node determines the discretization of the corresponding gene product's expression level. Considering the node $Y$ above, the two thresholds mean that $Y$ has three discrete expression levels, 0, 1, and 2 (low, intermediate, and high), where the integer indicates how many thresholds the value of $Y$ exceeds; i.e. $Y < \theta_{X,Y}, \theta_{Z,Y}$ corresponds to 0, and so forth. The consequence is that nodes in a GRN can have, and generally do have, different discretization levels. In the mini wavepool, Swi4 and Clb2 each have 3 out-edges and therefore 4 states (low, intermediate low, intermediate high, and high); Ndd1 has 2 out-edges and therefore 3 states (low, intermediate, and high); and Nrm1/Yox1 and Swi5/Ace2 each have 1 out-edge and therefore the Boolean states 0 and 1 (low and high).

The logic parameter for each node within the network orders the input values to the node with respect to the output thresholds of the node. For example, the $X$ node in Fig~\ref{fig:example-networks} (Right) has two in-edges, one from X and one from Y. Therefore, it has four possible input values: 
$$l_{X,Y}l_{X,Z} < \{ l_{X,Y}h_{X,Z}, \; h_{X,Y}l_{X,Z} \} < h_{X,Y}h_{X,Z}$$ 
that are partially ordered due to the constraint that $0 < l < h$. In particular, notice that the input values $l_{X,Y}h_{X,Z}$ and $h_{X,Y}l_{X,Z}$ cannot have a determined order until real values are assigned to $l$ and $h$. The node $X$ has a single out-edge to node $Y$, giving it one threshold $\theta_{Y,X}$. A logic parameter is the insertion of this threshold into the partial order of inputs. For $X$, there are six possible logic parameters, since there are six possible ways to insert the threshold into the partial order:
\begin{align*}
    & l_{X,Y}l_{X,Z} < \{l_{X,Y}h_{X,Z}, \; h_{X,Y}l_{X,Z}\} < h_{X,Y}h_{X,Z} < \theta_{Y,X} \\
    & l_{X,Y}l_{X,Z} < \{l_{X,Y}h_{X,Z}, \; h_{X,Y}l_{X,Z}\} < \theta_{Y,X} < h_{X,Y}h_{X,Z} \\
    & l_{X,Y}l_{X,Z} < l_{X,Y}h_{X,Z}  < \theta_{Y,X} < h_{X,Y}l_{X,Z} < h_{X,Y}h_{X,Z} \\
    & l_{X,Y}l_{X,Z} < h_{X,Y}l_{X,Z} < \theta_{Y,X} < l_{X,Y}h_{X,Z}  < h_{X,Y}h_{X,Z} \\
    & l_{X,Y}l_{X,Z}   < \theta_{Y,X} < \{l_{X,Y}h_{X,Z}, \; h_{X,Y}l_{X,Z}\} < h_{X,Y}h_{X,Z} \\
    & \theta_{Y,X} < l_{X,Y}l_{X,Z}   < \{l_{X,Y}h_{X,Z}, \; h_{X,Y}l_{X,Z}\} < h_{X,Y}h_{X,Z}
\end{align*}
A DSGRN parameter is a collection of inequalities: a choice of one order parameter and one logic parameter for each node in the network. We will call the inequalities for node $i$ a \textit{factor parameter for $i$.}
We remark that the number of DSGRN parameters scales poorly with the number of edges in a network. As the number of in-edges to a node grows, the size of the partial order grows exponentially. As the number of out-edges from a node grows, the number of order parameters grows factorially. Additionally, the complexity of inserting multiple thresholds into the partial order causes a large increase in the number of logic parameters.

\subsection{Parameter graph and remainder parameter}\label{sec:pg}

An important element of this work is the transition between WT and mutant dynamical phenotypes in parameter space. The following two sections explain the details of this transition.

The collection of all possible factor parameters for $i$ can be represented as an undirected graph, called the factor graph for node $i$. These factor graphs have nodes representing each factor parameter and the edges between these nodes represent a single change in an inequality.  

Each factor graph for a node $i$ can be written as the product of $M$ logic graphs, where the nodes of a logic graph are the collection of all logic parameters for $i$ and an edge exists between two nodes when there is a single change in a logic parameter inequality (for a fixed order parameter). The number of logic graphs $M$ is the number of order parameters for node $i$. Connections between logic graphs exist whenever there is a single change in the order parameter while the logic parameter remains the same, provided that the two swapped thresholds in the order parameter have no intervening logic value. This is most easily seen by an example.

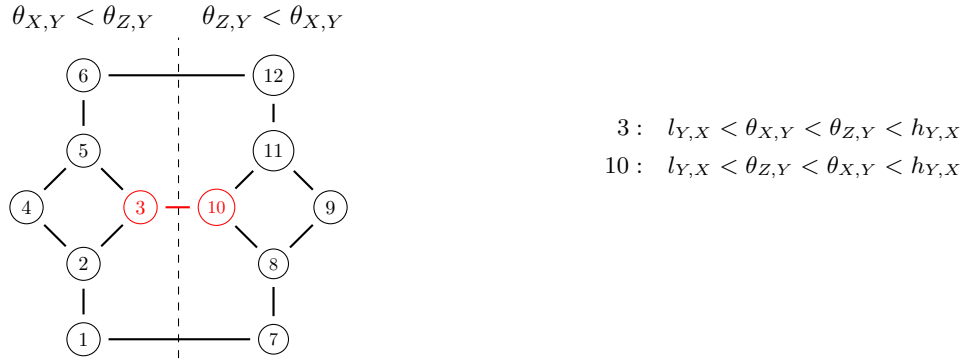
\begin{figure}[!h]
\begin{minipage}{0.5\linewidth}
\begin{center}
\begin{tikzpicture}[rotate around x=0, rotate around y=0, rotate around z=0, main node/.style={circle, draw, scale = 0.75, inner sep=1pt, minimum size=0pt}, node distance=2cm,z={(60:-0.5cm)}]

\draw[dashed] (1.25,-.25)--(1.25,4);

\node[main node] (n1) at (0,0) [circle,inner sep=3pt]{1};
\node[main node] (n2) at (0,1)[circle,inner sep=3pt]{2};
\node[main node] (n3) at (-.75,1.75) [circle,inner sep=3pt]{4};
\node[main node] (n4) at (.75,1.75)[circle,red,inner sep=3pt]{3};
\node[main node] (n5) at (0,2.5) [circle,inner sep=3pt]{5};
\node[main node] (n6) at (0,3.5)[circle,inner sep=3pt]{6};
\node[main node] (n7) at (2.5,0) [circle,scale= .9,inner sep=3pt]{7};
\node[main node] (n8) at (2.5,1)[circle,scale=.9,inner sep=3pt]{8};
\node[main node] (n9) at (1.75,1.75) [circle,red,scale=.9,inner sep=3pt]{10};
\node [main node](n10) at (3.25,1.75)[circle,inner sep=3pt]{9};
\node[main node] (n11) at (2.5,2.5) [circle,inner sep=3pt]{11};
\node[main node] (n12) at (2.5,3.5)[circle,inner sep=3pt]{12};

\node (n13) at (0,4.25) [ellipse,inner sep=3pt]{$\theta_{X,Y}<\theta_{Z,Y}$};
\node (n14) at (2.5,4.25)[ellipse,inner sep=3pt]{$\theta_{Z,Y}<\theta_{X,Y}$};

\path[>=angle 90,thick]
(n1) edge[shorten <= 3pt, shorten >= 3pt] node[] {} (n2)
(n2) edge[shorten <= 3pt, shorten >= 3pt,] node[] {} (n3)
(n2) edge[shorten <= 3pt, shorten >= 3pt] node[] {} (n4)
(n3) edge[shorten <= 3pt, shorten >= 3pt] node[] {} (n5)
(n4) edge[shorten <= 3pt, shorten >= 3pt] node[] {} (n5)
(n5) edge[shorten <= 3pt, shorten >= 3pt] node[] {} (n6)
(n6) edge[ shorten <= 3pt, shorten >= 3pt] node[] {} (n12)
(n4) edge[red, shorten <= 3pt, shorten >= 3pt] node[] {} (n9)
(n1) edge[shorten <= 3pt, shorten >= 3pt] node[] {} (n7)
(n7) edge[shorten <= 3pt, shorten >= 3pt] node[] {} (n8)
(n8) edge[ shorten <= 3pt, shorten >= 3pt] node[] {} (n9)
(n8) edge[shorten <= 3pt, shorten >= 3pt] node[] {} (n10)
(n10) edge[ shorten <= 3pt, shorten >= 3pt] node[] {} (n11)
(n11) edge[shorten <= 3pt, shorten >= 3pt] node[] {} (n12)
(n9) edge[shorten <= 3pt, shorten >= 3pt] node[] {} (n11);
\end{tikzpicture}
\end{center}
\end{minipage}%
\hfill
\begin{minipage}{0.5\linewidth}
\begin{small}
\begin{align*}
3:\quad & l_{Y,X} < \theta_{X,Y}  < \theta_{Z,Y} <  h_{Y,X} \\
10:\quad & l_{Y,X} < \theta_{Z,Y} < \theta_{X,Y} <  h_{Y,X} \\
\end{align*}
\end{small}
\end{minipage}%
\caption{ The factor graph for the $Y$ node in Fig~\ref{fig:example-networks} (Right), shown as the product of two logic graphs, one for each possible threshold ordering for $Y$. The threshold ordering can be seen above the corresponding logic graph of the factor graph.}
\label{fig:Y_logic_graph}
\end{figure}

Let us now construct the factor graph for $Y$ from the three-node network in Fig~\ref{fig:example-networks} (Right). Node $Y$ has a single in-edge and two out-edges, meaning that two thresholds are inserted between the low and high production rates of $Y$: $0 < l_{Y,X} < h_{Y,X}$. This results in six logic parameters for $Y$:
\begin{align*}
    1:& \;\;\; l_{Y,X} < h_{Y,X} < \theta_1 < \theta_2 \\
    2:& \;\;\; l_{Y,X} < \theta_1 < h_{Y,X} < \theta_2 \\
    3:& \;\;\; l_{Y,X} < \theta_1 < \theta_2 < h_{Y,X} \\
    4:& \;\;\; \theta_1 < l_{Y,X} < h_{Y,X} < \theta_2 \\
    5:& \;\;\; \theta_1 < l_{Y,X} < \theta_2  < h_{Y,X}\\
    6:& \;\;\; \theta_1 < \theta_2 < l_{Y,X} < h_{Y,X}
\end{align*}
where $\theta_1 < \theta_2$ is some ordering of the two thresholds of $Y$,  $\theta_{X,Y}$ and $\theta_{Z,Y}$. These six inequalities are the nodes of the logic graph of $Y$. Due to the two order parameters of $Y$, there are two copies of this logic graph, one for each ordering, as shown in Fig~\ref{fig:Y_logic_graph}. The numbering of the logic parameters above is associated to the logic graph left of the dashed line, which is associated to order parameter $\theta_{X,Y}<\theta_{Z,Y}$. The isomorphic logic graph for order parameter $\theta_{Z,Y}<\theta_{X,Y}$ is on the right. Edges only exist between the two logic graphs when there is a single inequality change in the order parameter between thresholds that are adjacent in the logic parameter. An example of this type of edge can be seen in red in Fig \ref{fig:Y_logic_graph}.

For a GRN $N$ with nodes $i=1,\dots,n$, the product of the factor graphs $FG(i)$ is the DSGRN parameter graph $PG(N)$: 
\[\prod\limits_{i=1}^n FG(i) = PG(N).\]
The parameter graph contains all possible DSGRN parameters as nodes and encodes adjacency between real-valued parameter regions as edges~\cite{cummins:2016}.

\begin{figure}
\begin{tabular}{l@{}c@{}c@{}}
\begin{tabular}{ccc}
{\huge $X$} & \hspace{0.9in}{\huge $Y$} & \hspace{0.7in}{\huge $Z$} 
% \node (n1) at (0,0) [scale=1.75]{X};
% \node (n2) at (3,0) [scale=1.75]{Y};
% \node (n3) at (6,0) [scale=1.75]{Z};
% \node (n4) at (0,0) {corresponding Z inequality};
\end{tabular}
\\
\\
\begin{varwidth}{\linewidth}
\begin{center}
\begin{tikzpicture}
\draw[thick] (1,1.25)--(1.5,.75);
\draw[thick] (1,.75)--(1.5,1.25);
\draw[thick] (5,1.25)--(5.5,.75);
\draw[thick] (5,.75)--(5.5,1.25);

\node (n1) at (2.5,0) [red,scale=0.6,circle,fill,inner sep=3pt]{};
\node (n2) at (2.5,.5)[scale=0.6,circle,fill,inner sep=3pt]{};
\node (n3) at (2,1) [scale=0.6,circle,fill,inner sep=3pt]{};
\node (n4) at (3,1)[scale=0.6,circle,fill,inner sep=3pt]{};
\node (n5) at (2.5,1.5) [scale=0.6,circle,fill,inner sep=3pt]{};
\node (n6) at (2.5,2)[scale=0.6,circle,fill,inner sep=3pt]{};
\node (n7) at (4,0) [scale=0.6,circle,fill,inner sep=3pt]{};
\node (n8) at (4,.5)[scale=0.6,circle,fill,inner sep=3pt]{};
\node (n9) at (3.5,1) [scale=0.6,circle,fill,inner sep=3pt]{};
\node (n10) at (4.5,1)[scale=0.6,circle,fill,inner sep=3pt]{};
\node (n11) at (4,1.5) [scale=0.6,circle,fill,inner sep=3pt]{};
\node (n12) at (4,2)[scale=0.6,circle,fill,inner sep=3pt]{};

\node (n13) at (0,0) [red,scale=0.6,circle,fill,inner sep=3pt]{};
\node (n14) at (0,.5)[scale=0.6,circle,fill,inner sep=3pt]{};
\node (n15) at (-.5,1) [scale=0.6,circle,fill,inner sep=3pt]{};
\node (n16) at (.5,1)[scale=0.6,circle,fill,inner sep=3pt]{};
\node (n17) at (0,1.5) [scale=0.6,circle,fill,inner sep=3pt]{};
\node (n18) at (0,2)[scale=0.6,circle,fill,inner sep=3pt]{};

\node (n19) at (6,0)[scale=0.6,circle,fill,inner sep=3pt]{};
\node (n20) at (6,1) [scale=0.6,circle,fill,inner sep=3pt]{};
\node (n21) at (6,2)[red,scale=0.6,circle,fill,inner sep=3pt]{};

\path[>=angle 90,thick]
(n1) edge[shorten <= 3pt, shorten >= 3pt] node[] {} (n2)
(n2) edge[shorten <= 3pt, shorten >= 3pt,] node[] {} (n3)
(n2) edge[shorten <= 3pt, shorten >= 3pt] node[] {} (n4)
(n3) edge[shorten <= 3pt, shorten >= 3pt] node[] {} (n5)
(n4) edge[shorten <= 3pt, shorten >= 3pt] node[] {} (n5)
(n5) edge[shorten <= 3pt, shorten >= 3pt] node[] {} (n6)
(n6) edge[ shorten <= 3pt, shorten >= 3pt] node[] {} (n12)
(n4) edge[shorten <= 3pt, shorten >= 3pt] node[] {} (n9)
(n1) edge[shorten <= 3pt, shorten >= 3pt] node[] {} (n7)
(n7) edge[shorten <= 3pt, shorten >= 3pt] node[] {} (n8)
(n8) edge[ shorten <= 3pt, shorten >= 3pt] node[] {} (n9)
(n8) edge[shorten <= 3pt, shorten >= 3pt] node[] {} (n10)
(n10) edge[ shorten <= 3pt, shorten >= 3pt] node[] {} (n11)
(n11) edge[shorten <= 3pt, shorten >= 3pt] node[] {} (n12)
(n9) edge[shorten <= 3pt, shorten >= 3pt] node[] {} (n11)

(n13) edge[ shorten <= 3pt, shorten >= 3pt] node[] {} (n14)
(n14) edge[shorten <= 3pt, shorten >= 3pt] node[] {} (n15)
(n14) edge[shorten <= 3pt, shorten >= 3pt] node[] {} (n16)
(n16) edge[shorten <= 3pt, shorten >= 3pt] node[] {} (n17)
(n15) edge[ shorten <= 3pt, shorten >= 3pt] node[] {} (n17)
(n17) edge[shorten <= 3pt, shorten >= 3pt] node[] {} (n18)

(n19) edge[ shorten <= 3pt, shorten >= 3pt] node[] {} (n20)
(n20) edge[shorten <= 3pt, shorten >= 3pt] node[] {} (n21);
\end{tikzpicture}
\end{center}
\end{varwidth}%

\begin{minipage}{0.50\linewidth}
\begin{align*}
Z: \theta_{X,Z} < l_{Z,Y}<h_{Z,Y}
\end{align*}
\end{minipage}%
\\
\\
\begin{varwidth}{\linewidth}
\begin{center}
\begin{tikzpicture}
\draw[thick] (1,1.25)--(1.5,.75);
\draw[thick] (1,.75)--(1.5,1.25);
\draw[thick] (5,1.25)--(5.5,.75);
\draw[thick] (5,.75)--(5.5,1.25);

\node (n1) at (2.5,0) [red,scale=0.6,circle,fill,inner sep=3pt]{};
\node (n2) at (2.5,.5)[scale=0.6,circle,fill,inner sep=3pt]{};
\node (n3) at (2,1) [scale=0.6,circle,fill,inner sep=3pt]{};
\node (n4) at (3,1)[scale=0.6,circle,fill,inner sep=3pt]{};
\node (n5) at (2.5,1.5) [scale=0.6,circle,fill,inner sep=3pt]{};
\node (n6) at (2.5,2)[scale=0.6,circle,fill,inner sep=3pt]{};
\node (n7) at (4,0) [scale=0.6,circle,fill,inner sep=3pt]{};
\node (n8) at (4,.5)[scale=0.6,circle,fill,inner sep=3pt]{};
\node (n9) at (3.5,1) [scale=0.6,circle,fill,inner sep=3pt]{};
\node (n10) at (4.5,1)[scale=0.6,circle,fill,inner sep=3pt]{};
\node (n11) at (4,1.5) [scale=0.6,circle,fill,inner sep=3pt]{};
\node (n12) at (4,2)[scale=0.6,circle,fill,inner sep=3pt]{};

\node (n13) at (0,0) [red,scale=0.6,circle,fill,inner sep=3pt]{};
\node (n14) at (0,.5)[scale=0.6,circle,fill,inner sep=3pt]{};
\node (n15) at (-.5,1) [scale=0.6,circle,fill,inner sep=3pt]{};
\node (n16) at (.5,1)[scale=0.6,circle,fill,inner sep=3pt]{};
\node (n17) at (0,1.5) [scale=0.6,circle,fill,inner sep=3pt]{};
\node (n18) at (0,2)[scale=0.6,circle,fill,inner sep=3pt]{};

\node (n19) at (6,0)[red,scale=0.6,circle,fill,inner sep=3pt]{};
\node (n20) at (6,1) [scale=0.6,circle,fill,inner sep=3pt]{};
\node (n21) at (6,2)[scale=0.6,circle,fill,inner sep=3pt]{};

\path[>=angle 90,thick]
(n1) edge[shorten <= 3pt, shorten >= 3pt] node[] {} (n2)
(n2) edge[shorten <= 3pt, shorten >= 3pt,] node[] {} (n3)
(n2) edge[shorten <= 3pt, shorten >= 3pt] node[] {} (n4)
(n3) edge[shorten <= 3pt, shorten >= 3pt] node[] {} (n5)
(n4) edge[shorten <= 3pt, shorten >= 3pt] node[] {} (n5)
(n5) edge[shorten <= 3pt, shorten >= 3pt] node[] {} (n6)
(n6) edge[ shorten <= 3pt, shorten >= 3pt] node[] {} (n12)
(n4) edge[shorten <= 3pt, shorten >= 3pt] node[] {} (n9)
(n1) edge[shorten <= 3pt, shorten >= 3pt] node[] {} (n7)
(n7) edge[shorten <= 3pt, shorten >= 3pt] node[] {} (n8)
(n8) edge[ shorten <= 3pt, shorten >= 3pt] node[] {} (n9)
(n8) edge[shorten <= 3pt, shorten >= 3pt] node[] {} (n10)
(n10) edge[ shorten <= 3pt, shorten >= 3pt] node[] {} (n11)
(n11) edge[shorten <= 3pt, shorten >= 3pt] node[] {} (n12)
(n9) edge[shorten <= 3pt, shorten >= 3pt] node[] {} (n11)

(n13) edge[ shorten <= 3pt, shorten >= 3pt] node[] {} (n14)
(n14) edge[shorten <= 3pt, shorten >= 3pt] node[] {} (n15)
(n14) edge[shorten <= 3pt, shorten >= 3pt] node[] {} (n16)
(n16) edge[shorten <= 3pt, shorten >= 3pt] node[] {} (n17)
(n15) edge[ shorten <= 3pt, shorten >= 3pt] node[] {} (n17)
(n17) edge[shorten <= 3pt, shorten >= 3pt] node[] {} (n18)

(n19) edge[ shorten <= 3pt, shorten >= 3pt] node[] {} (n20)
(n20) edge[shorten <= 3pt, shorten >= 3pt] node[] {} (n21);
\end{tikzpicture}
\end{center}
\end{varwidth}%

\begin{minipage}{0.50\linewidth}
\begin{align*}
Z: l_{Z,Y}<h_{Z,Y}<\theta_{X,Z}
\end{align*}
\end{minipage}%

\end{tabular}
\caption{The DSGRN parameter graph of the three-node network in Fig~\ref{fig:example-networks} (Right) is the product of three factor graphs corresponding to the nodes $X$, $Y$, and $Z$.  The example $X$, $Y$, $Z$ triples of red factor parameters (top panel and bottom panel) correspond to two different DSGRN parameters for the network. The fixed $X$, $Y$ factor parameters together form a remainder parameter for $Z$. The inequalities for the two different factor parameters for $Z$ are shown to the left.}
\label{fig:remainder_param_3n}
\end{figure}
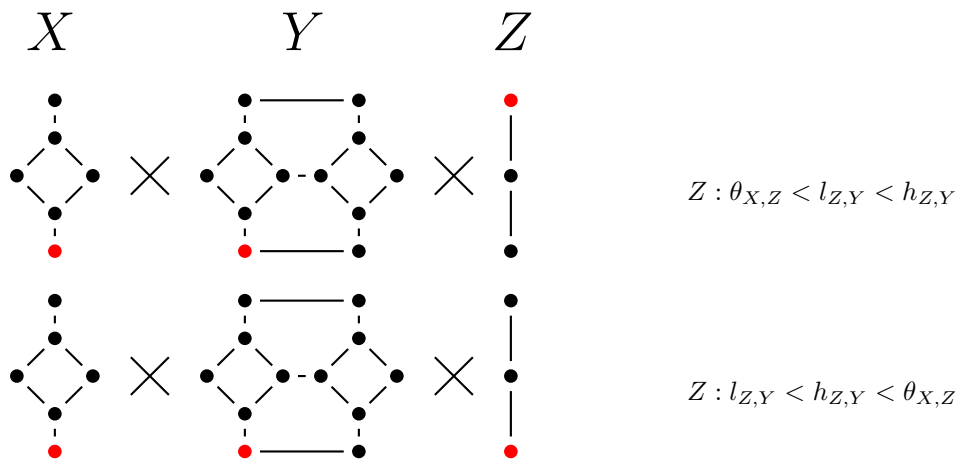

As an example, the DSGRN parameter graph for the three-node network in Fig~\ref{fig:example-networks} (Right) is shown in Fig~\ref{fig:remainder_param_3n}. Each choice of a node from the $X$ factor graph, a node from the $Y$ factor graph, and a node from the $Z$ factor graph is a DSGRN parameter node. Two such DSGRN parameter nodes are shown in red in the top and bottom panels respectively in Fig~\ref{fig:remainder_param_3n}. By examining the factor graphs, we see that the DSGRN parameter graph for the three-node network has a size of $6\times12\times3 = 216$. 

% An important idea for our research is that of the remainder parameter. 
Suppose a single factor parameter for a node $x_i$ is fixed, for example, representing a mutant phenotype such as a knock-out. The combination of the factor parameters for all of the remaining nodes form a \textit{remainder parameter}. We demonstrate the idea of a remainder parameter using Fig~\ref{fig:remainder_param_3n} by  allowing the parameter node of $Z$ to vary while those of $X$ and $Y$ remain fixed. In this case we have the remainder parameter composed of the $X$ and $Y$ factor parameters marked in red as $Z$ varies from its highest parameter to its lowest.  The high and low $Z$ factor parameters could represent an up-regulation or knock-out of a gene in a biological scenario. In the first case, $Z$ is continuously expressed at its highest level while in the second, it is expressed at its lowest level. The lowest expression level may as well be taken to be zero, as in a knock-out, since $Z$ is never expressed highly enough to regulate its downstream target.

% When a pair of DSGRN parameters differs only in a logic parameter of $x_i$ but has the same remainder parameter, we say that the pair of DSGRN parameters is \textit{$x_i$-divergent}. The pair of DSGRN parameters shown in red in Fig~\ref{fig:remainder_param_3n} (top vs bottom panel) are $Z$-divergent. By examining dynamical phenotypes, we can establish how the network dynamics differ if only the node $Z$ is disturbed. This gives us the ability to investigate how an up-regulation or a knock-down of a specific gene can affect the dynamics of a GRN. The goal of our research is investigating the effect that Clb2 activity  has on the mini pulse generator and using Clb2-divergent DSGRN parameter pairs we can directly track and determine this effect.

\subsection{Application to the mini wavepool}\label{sec:methods:clb2}

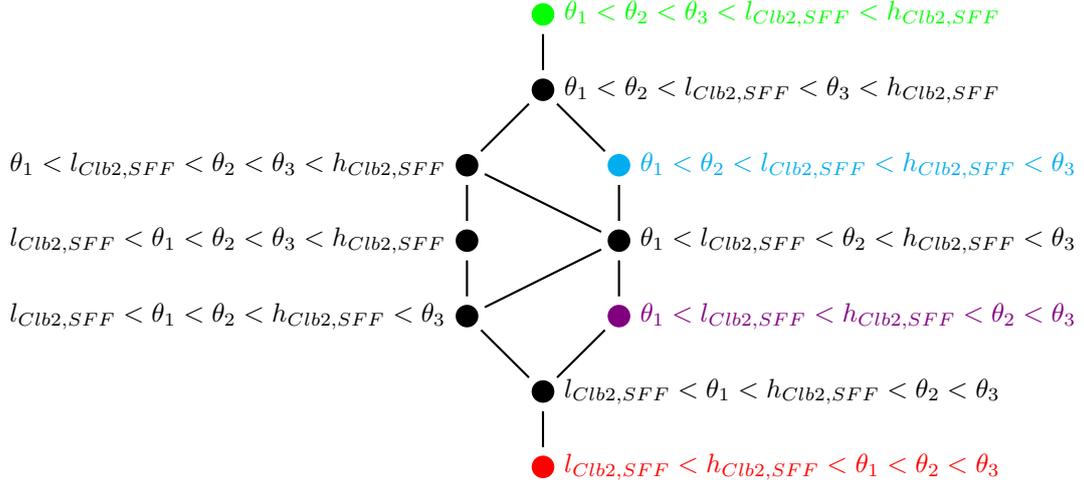
\begin{figure}[h!]
\centering
\begin{tikzpicture}
\node (n1) at (0,5) [label=left:$\theta_1<l_{Clb2,SFF}<\theta_2<\theta_3<h_{Clb2,SFF}$,circle,fill,inner sep=3pt]{};
\node (n2) at (0,4)[label=left:$l_{Clb2,SFF}<\theta_1<\theta_2<\theta_3<h_{Clb2,SFF}$,circle,fill,inner sep=3pt]{};
\node (n3) at (0,3) [label=left:$l_{Clb2,SFF}<\theta_1<\theta_2<h_{Clb2,SFF}<\theta_3$,circle,fill,inner sep=3pt]{};
\node[color=green] (n4) at (1,7)[label=right:\textcolor{green}{$\theta_1<\theta_2<\theta_3<l_{Clb2,SFF}<h_{Clb2,SFF}$},circle,fill,inner sep=3pt]{};
\node (n5) at (1,6) [label=right:$\theta_1<\theta_2<l_{Clb2,SFF}<\theta_3<h_{Clb2,SFF}$,circle,fill,inner sep=3pt]{};
\node[color=cyan] (n6) at (2,5)[label=right:\textcolor{cyan}{$\theta_1<\theta_2<l_{Clb2,SFF}<h_{Clb2,SFF}<\theta_3$},circle,fill,inner sep=3pt]{};
\node (n7) at (2,4) [label=right:$\theta_1<l_{Clb2,SFF}<\theta_2<h_{Clb2,SFF}<\theta_3$,circle,fill,inner sep=3pt]{};
\node[color=violet] (n8) at (2,3)[label=right:\textcolor{violet}{$\theta_1<l_{Clb2,SFF}<h_{Clb2,SFF}<\theta_2<\theta_3$},circle,fill,inner sep=3pt]{};
\node (n9) at (1,2) [label=right:$l_{Clb2,SFF}<\theta_1<h_{Clb2,SFF}<\theta_2<\theta_3$,circle,fill,inner sep=3pt]{};
\node[color=red] (n10) at (1,1)[label=right:\textcolor{red}{$l_{Clb2,SFF}<h_{Clb2,SFF}<\theta_1<\theta_2<\theta_3$},circle,fill,inner sep=3pt]{};

\path[>=angle 90,thick]
(n1) edge[shorten <= 3pt, shorten >= 3pt] node[] {} (n2)
(n2) edge[shorten <= 3pt, shorten >= 3pt,] node[] {} (n3)
(n3) edge[shorten <= 3pt, shorten >= 3pt] node[] {} (n9)
(n1) edge[shorten <= 3pt, shorten >= 3pt] node[] {} (n7)
(n3) edge[shorten <= 3pt, shorten >= 3pt] node[] {} (n7)
(n5) edge[shorten <= 3pt, shorten >= 3pt] node[] {} (n1)
(n4) edge[ shorten <= 3pt, shorten >= 3pt] node[] {} (n5)
(n5) edge[shorten <= 3pt, shorten >= 3pt] node[] {} (n6)
(n6) edge[shorten <= 3pt, shorten >= 3pt] node[] {} (n7)
(n7) edge[shorten <= 3pt, shorten >= 3pt] node[] {} (n8)
(n8) edge[ shorten <= 3pt, shorten >= 3pt] node[] {} (n9)
(n9) edge[shorten <= 3pt, shorten >= 3pt] node[] {} (n10)

;                       
\end{tikzpicture}
\caption{One of the six logic parameter graphs corresponding to the Clb2 node from the mini wavepool in Fig~\ref{fig:wavepool} (Right), where the thresholds $\theta_1, \theta_2, \theta_3$ are associated to the nodes for SBF, SFF, and Swi5 via some fixed mapping. The Clb2 ON logic parameter is represented in green at the top of the factor graph. The Clb2 OFF logic parameter is represented in red at the bottom of the factor graph. The Clb2 INT-H logic parameter in blue is two steps up from the Clb2 INT-L logic parameter in violet. The WT logic parameters are the remaining black inequalities. The checkpoint phenotypes are not restricted to any particular Clb2 logic parameter.}
\label{fig:Clb2-factor-graph}
\end{figure}

Recall from Fig~\ref{fig:wavepool} (Right) that Clb2 has one in-edge and three out-edges. This means it has $3! = 6$ order parameters and ten logic parameters (see Fig~\ref{fig:Clb2-factor-graph}). Four of the ten logic parameters are taken to be representative of various Clb2 mutants (see the caption of Fig~\ref{fig:Clb2-factor-graph}). We propose that the WT phenotype is associated to one of the logic parameters in black. Notice that each of the logic parameters that denote the state of a Clb2 mutant (inequalities in color) have both the low and high values $l_{Clb2,SFF}$ and $h_{Clb2,SFF}$ together between thresholds. In other words, the model of Clb2 activity implies that even if the molecular concentration isn't perfectly constant, there will never be a sufficiently large change in concentration to trigger a change in regulation at downstream genes. In contrast, the logic parameters for Clb2 WT ensures that changes in Clb2 concentration will impact at least one downstream target. We do not require regulatory activity at all downstream nodes, because the only information provided by the data is that Clb2 has noticeable oscillations. The DSGRN parameters at which we see consistency with the experimental data can give us guidance in assessing which Clb2 interactions may be important in the WT phenotype.

Dynamical phenotypes II and III involve transitions in the Clb2 logic graph from WT to mutant parameters that exhibit consistency with the corresponding mutant time series when the remainder parameter (referred to as the mini pulse generator parameter) is constant.
%In the context of the mini wavepool, the remainder parameter is referred to as the mini pulse generator parameter.
As an example, consider the Clb2 ON transcriptional phenotype. A match according to dynamical phenotype II is a single mini pulse generator parameter that exhibits (1) mutant cycling consistent with the \textit{cdc20$\Delta$} dataset at the phenotype-permissible green node at the top of the Clb2 logic graph as well as (2) WT cycling at any one of the six black nodes in the Clb2 logic graph.

\subsection{State transition graphs and Morse graphs}\label{sec:stg}

Each DSGRN parameter has a corresponding state transition graph (STG) that graphically represents the dynamics of the network at that DSGRN parameter. For ODE systems, the dynamics of a network are described by time-dependent trajectories in phase space in which all gene products associated to the network are changing concentration. Phase space is the $N$-dimensional real-valued and positive space where each coordinate represents a node in the GRN. DSGRN does not concern itself with these trajectories but instead looks at directional flow across thresholds and considers paths through an STG. 
 
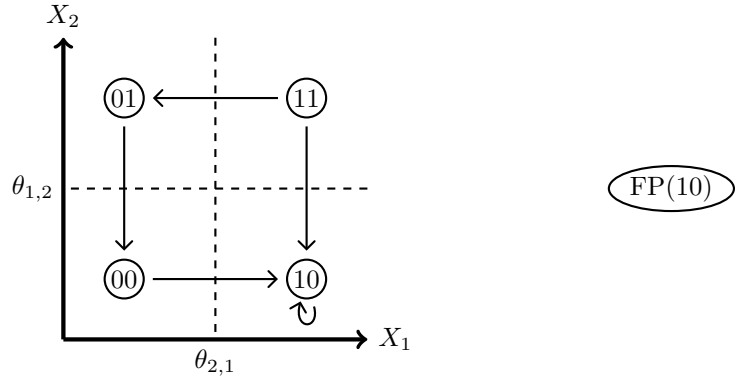
\begin{figure}[!h]
\centering
\begin{tikzpicture}[main node/.style={circle, draw, thick, inner sep=1pt, minimum size=0pt},scale=0.8]
\draw[dashed, thick] (5,2.5)--(0,2.5) node [left]{$\theta_{1,2}$};
\draw[dashed, thick] (2.5,5)--(2.5,0) node [below]{$\theta_{2,1}$};
\draw[->,ultra thick] (0,0)--(5,0) node[right]{$X_1$};
\draw[->,ultra thick] (0,0)--(0,5) node[above]{$X_2$};

\node[main node] (n1) at (1,1) {00};
\node[main node] (n2) at (1 , 4){01};
\node[main node] (n3) at (4,4) {11};
\node[main node] (n4) at (4 , 1){10};

\path[->,>=angle 90,thick]
(n1) edge[shorten <= 3pt, shorten >= 3pt] node[] {} (n4)
(n3) edge[shorten <= 3pt, shorten >= 3pt] node[] {} (n2)
(n3) edge[shorten <= 3pt, shorten >= 3pt] node[] {} (n4)
(n2) edge[shorten <= 3pt, shorten >= 3pt] node[] {} (n1)
(n4) edge[shorten <= 3pt, shorten >= 3pt, loop below] node[] {} (n4)
;

\node[main node,shape=ellipse] () at (10,2.5) {FP(10)};

\end{tikzpicture}

\caption{(Left) The  STG superposed over the rectangular domains dividing phase space for the toggle switch with the DSGRN parameter in~\eqref{eq:toggle-switch-param}. (Right) The corresponding MG.}
\label{fig:stg}
\end{figure}

 An example phase space for the toggle switch is shown in Fig~\ref{fig:stg}. We discretize phase space by dividing it up into rectangular boxes called domains using the collection of thresholds for the switching system; these are the dotted lines in Fig~\ref{fig:stg}, showing the division of the positive plane into four domains. Each domain corresponds to a level of $X_1$ and $X_2$ where $1$ indicates above-threshold values and $0$ indicates below-threshold values. As an example the domain labeled (01) corresponds to high $X_2$ and low $X_1$. We say that $(01)$ is the \textit{state} corresponding to the upper left domain. 
 
 We can define flow across the boundaries of the domains by specifying a DSGRN parameter. 
 In Fig~\ref{fig:stg}, we use
\begin{align}
 DSG & RN  \quad parameter: \quad \label{eq:toggle-switch-param}\\
X_1:\quad & l_{1,2} < \theta_{2,1} < h_{1,2} \nonumber\\
X_2: \quad & l_{2,1} < h_{2,1} < \theta_{1,2}\nonumber
 \end{align}
 to determine the arrows. For example, the top arrow pointing left and the right arrow pointing down are determined by considering the starting domain, where $X_1$ and $X_2$ have high concentrations. In this case, the concentrations of $X_1$ and $X_2$ are both decreasing since the nodes are effective in repressing each other. However, note that $X_2$ is \textit{always} decreasing. This happens because of the DSGRN parameter for $X_2$. Regardless of the presence or absence of $X_1$, $X_2$ will never increase across its threshold $\theta_{1,2}$ since its high concentration is still below $\theta_{1,2}$. There is also a DSGRN parameter where $X_2$ (or symmetrically $X_1$) is always increasing, indicating that one or both nodes are ineffective repressors.

The STG for the toggle switch is superposed on phase space in Fig~\ref{fig:stg}. The nodes of the STG correspond to the states of the domains and the directed edges between states generally correspond to the flow across the intervening boundary. The one exception is that self-edges are added whenever it is not possible to leave that domain because all of the flow points inward; see the lower right domain in Fig~\ref{fig:stg}.

The size of the STG grows with the size of the GRN, since the addition of a node to the GRN adds another dimension to phase space and the addition of an edge, or threshold value, adds another set of domains to an existing dimension of phase space. Both of these mechanisms increase the number of nodes in the STG, which rapidly becomes large and difficult to interpret. It is therefore useful for clarity to examine a summary of the STG called a Morse graph (MG).
We build the MG from the recurrent components of the STG. A recurrent component, or Morse set, is a maximal set of nodes in the STG that contain a path from any domain $u$ to any other domain $v$ for all $u,v$ within the recurrent component. Each of these Morse sets are represented in the MG as Morse nodes, where an edge between Morse nodes indicates that there is a path between some domain $u$ in Morse set 1 to some domain $v$ in Morse set 2.

Each Morse node in the MG has an annotation indicating the dynamics in the associated Morse set. The annotations consist of the labels full cycle (FC), partial cycle (PC), and fixed point (FP). The FC annotation indicates a Morse set in which there is a looped path, or cycle, in the STG that crosses at least one threshold value for each $x_i$. The PC annotation indicates a Morse set in which the looped path crosses thresholds for only a subset of the $x_i$. The FP annotation indicates a Morse set consisting of a single state with a self-edge. The Morse graph for the toggle switch is the FP with domain coordinates (10), denoted FP(10); see Fig~\ref{fig:stg}. The SAC and DRC FPs for the mini wavepool network that are shown qualitatively in Table~\ref{tab:qualitative_fps} are represented by the domain coordinates shown in Table~\ref{fig:sac_fps}.

\begin{table}[]
\centering
\begin{tabular}{|lll|}
\hline
\multicolumn{3}{|l|}{Checkpoint FPs: FP(Swi4, Nrm1/Yox1, Ndd1, Swi5/Ace2, Clb2)}                                                                                                                                                                                                \\ \hline
\multicolumn{1}{|l|}{SAC FPs:} & \multicolumn{2}{l|}{\begin{tabular}[c]{@{}l@{}}All proxies: \\  FP(0, 0, 2, 1, 1) \\  FP(0, 0, 2, 1, 2) \\  FP(0, 0, 2, 1, 3)\end{tabular}}                                                                                                     \\ \hline
\multicolumn{1}{|l|}{DRC FPs:} & \begin{tabular}[c]{@{}l@{}}Yox1 proxies:\\  FP(0, 0, 2, 1, 1) \\  FP(0, 0, 2, 1, 2) \\  FP(0, 0, 2, 1, 3)\end{tabular} & \begin{tabular}[c]{@{}l@{}}Nrm1 Proxies:\\  FP(0, 1, 2, 1, 1) \\  FP(0, 1, 2, 1, 2) \\  FP(0, 1, 2, 1, 3)\end{tabular} \\ \hline
\end{tabular}
\caption{The order of the nodes for the checkpoint phenotype FPs is Swi4, which has activity levels in the range 0-3, Nrm1/Yox1, which has levels 0-1, Ndd1, which has levels 0-2, Clb2, which has levels 0-3, and Swi5/Ace2, which has levels 0-1. See Section~\ref{sec:dsgrn} for an explanation of the differing state discretizations for each of the nodes.} 
\label{fig:sac_fps}
\end{table}

The Morse graph also encodes the \textit{stability} of a Morse set. Stability in the sense of dynamical systems roughly means that trajectories close to a stable manifold in phase space will approach that manifold asymptotically over time. Stability of a Morse node in the Morse graph, whether a fixed point or a cycle, is identified with having no possible exit from a Morse set once a path enters the Morse set; i.e., a Morse node is stable if and only if it has no out-edges in the MG. Morse nodes with out-edges are unstable, since there are potential exit paths from the associated Morse set. 

In the mini wavepool model, wild-type cycling is evaluated within stable FCs, mutant cycling is evaluated in stable PCs where Clb2 is in a fixed state and the mini pulse generator nodes are oscillating, and checkpoints are evaluated as FPs.

As an example, consider again the single recurrent component in Fig~\ref{fig:stg}, labeled FP(10). By inspection of the STG we can see that the domain (10) has no out edges to other domains therefore it is a stable FP. One can see that the key feature of the STG, namely the single equilibrium, is more readily identifiable from the MG rather than the STG. In general, the simplicity of the MG makes it much easier to interpret the dynamical behavior of a GRN at a particular DSGRN parameter.

In order to understand the implications of the different dynamic phenomena described above, we will analyze by hand the example 3-node network seen on the right in Fig \ref{fig:example-networks}, which exhibits all possible annotations. We will be looking at three different DSGRN parameters that give rise to three different Morse graphs. Since this network has three nodes, phase space will have three dimensions where the states are represented in the order ($XYZ$); for example, (010), (110), etc. The $X$ and $Z$ dimensions are divided into two domains each since they have only one out-edge apiece. The $Y$ dimension contains three domains since there are two out-edges. This means that instead of two states, 0 and 1, $Y$ has three states, 0, 1, and 2, representing the three possible positions with respect to two distinct thresholds.

The first example DSGRN parameter we will investigate gives rise to an MG containing a stable FC. The corresponding STG and inequalities can be seen in Fig \ref{fig:fc_stg}. The nodes have been arranged in the same spatial order as the rectangular domains in phase space. The red arrows and nodes indicate the path that contains the FC. It can be seen through inspection of the STG that starting anywhere within the STG leads to this stable cycle. This cycle is considered an FC because the set of domains that are passed through encompass all three directions in phase space.

\begin{figure}[!htbp]
\centering
\begin{minipage}{0.75\linewidth}
\begin{tikzpicture}[rotate around x=0, rotate around y=0, rotate around z=0, main node/.style={circle, draw, thick, inner sep=1pt, minimum size=0pt}, scale=2.4, node distance=1cm,z={(60:-0.5cm)}]
\tikzstyle{every loop}=[looseness=14]

\node[main node, fill=white] (000) at (0,0,0) {$000$} ;
\node[main node, fill=white] (001) at (0,1,0) {$001$} ;
\node[main node, fill=white] (010) at (1,0,0) {$010$} ;
\node[main node, fill=white] (100) at (0,0,1) {$100$} ;
\node[main node,red, fill=white] (110) at (1,0,1) {$110$} ;
\node[main node, fill=white] (101) at (0,1,1) {$101$} ;
\node[main node,red, fill=white] (011) at (1,1,0) {$011$} ;
\node[main node,red, fill=white] (111) at (1,1,1) {$111$} ;
\node[main node,red, fill=white] (020) at (2,0,0) {$020$} ;
\node[main node,red, fill=white] (120) at (2,0,1) {$120$} ;
\node[main node,red, fill=white] (021) at (2,1,0) {$021$} ;
\node[main node, fill=white] (121) at (2,1,1) {$121$} ;

\draw[very thick,-,shorten >= 0pt,shorten <= 0pt]
(101) edge[->] (001)
(000) edge[dashed, ->] (001)
(000) edge[dashed,->] (100)
(100) edge[->] (101)
(101) edge[line width=0.1cm, white] (111)
(010) edge[dashed, ->] (011)
(010) edge[dashed,->] (110)
(001) edge[->] (011)
(000) edge[dashed,->] (010)
(100) edge[->] (110)
(110) edge[line width=0.1cm, white] (111)
(110) edge[->, red] (111)

(111) edge[->,red] (011)
(121) edge[line width=0.1cm, white] (120)
(101) edge[->] (111)

(111) edge[line width=0.1cm, white] (121)
(121) edge[->] (021)

(121) edge[->] (120)
(021) edge[->,red] (020)
(020) edge[->,red] (120)
(011) edge[->,red] (021)
(121) edge[->] (111)
(120) edge[->, red] (110)
(010) edge[dashed, ->] (020);

\node[main node, red, shape=ellipse,scale=1.5] () at (3,0.5) {$FC$};

\end{tikzpicture}
\end{minipage}%
\smallskip

\begin{align*}
 DSG & RN  \quad parameter \quad FC:: \quad \\
X:\quad & l_{X,Y}l_{X,Z} < h_{X,Y}l_{X,Z} < \theta_{Y,X}< l_{X,Y}h_{X,Z} < h_{X,Y}h_{X,Z} \\
Y: \quad & \theta_{X,Y}<l_{Y,X}<\theta_{Z,Y}<h_{Y,X}\\
Z: \quad & l_{Z,Y}<\theta_{X,Z}<h_{Z,Y}
\end{align*}

\caption{The corresponding STG (Left), MG (Right), and set of inequalities (Bottom) for an example DSGRN parameter for the 3-node network in Fig~\ref{fig:example-networks} exhibiting a stable full cycle shown by the red nodes and edges. The dashed lines are for visual effect only. They represent edges on the ``back'' of the cubes. }
\label{fig:fc_stg}
\end{figure}

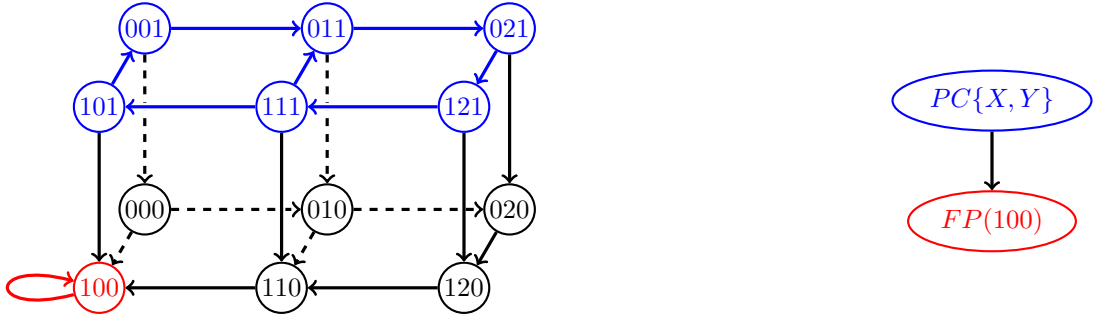
\begin{figure}[!htbp]
\centering
\begin{minipage}{.75\linewidth}
\begin{tikzpicture}[rotate around x=0, rotate around y=0, rotate around z=0, main node/.style={circle, draw, thick, inner sep=1pt, minimum size=0pt}, scale=2.4, node distance=1cm,z={(60:-0.5cm)}]
\tikzstyle{every loop}=[looseness=14]

\node[main node,red, fill=white] (000) at (0,0,0) {$000$} ;
\node[main node,blue, fill=white] (001) at (0,1,0) {$001$} ;
\node[main node,red, fill=white] (010) at (1,0,0) {$010$} ;
\node[main node,red, fill=white] (100) at (0,0,1) {$100$} ;
\node[main node,red, fill=white] (110) at (1,0,1) {$110$} ;
\node[main node,blue, fill=white] (101) at (0,1,1) {$101$} ;
\node[main node,blue, fill=white] (011) at (1,1,0) {$011$} ;
\node[main node,blue, fill=white] (111) at (1,1,1) {$111$} ;
\node[main node,red, fill=white] (020) at (2,0,0) {$020$} ;
\node[main node,red, fill=white] (120) at (2,0,1) {$120$} ;
\node[main node,blue, fill=white] (021) at (2,1,0) {$021$} ;
\node[main node,blue, fill=white] (121) at (2,1,1) {$121$} ;

\draw[very thick,-,shorten >= 0pt,shorten <= 0pt]
(101) edge[->,blue] (001)
(001) edge[dashed,->] (000)
(100) edge[dashed,->,red] (000)
(101) edge[->] (100)
(011) edge[dashed,->] (010)
(011) edge[line width=0.1cm, white] (111)
(110) edge[red,dashed,->] (010)
(001) edge[->,blue] (011)
(000) edge[dashed,->,red] (010)
(110) edge[->,red] (100)
(111) edge[line width=0.1cm, white] (110)
(111) edge[->] (110)

(021) edge[line width=0.1cm, white] (121)
(111) edge[->,blue] (011)
(121) edge[line width=0.1cm, white] (120)
(111) edge[->,blue] (101)

(021) edge[->,blue] (121)
(121) edge[->] (120)
(021) edge[->] (020)
(020) edge[->,red] (120)
(011) edge[->,blue] (021)
(121) edge[->,blue] (111)
(120) edge[->, red] (110)
(010) edge[dashed,->,red] (020);
\end{tikzpicture}
\end{minipage}%
\hfill
\begin{minipage}{0.25\linewidth}
\begin{tikzpicture}[main node/.style={ellipse, draw, thick, inner sep=3pt, minimum size=0pt},scale=0.8]
\node[main node, red] (n1) at (0,0) {$PC\{X,Y\}$};
\node[main node, blue] (n2) at (0,2) {$PC\{X,Y\}$};

\draw[very thick,-,shorten >= 0pt,shorten <= 0pt]
(n2) edge[->] (n1);

\end{tikzpicture}
\end{minipage}%

\begin{align*}
 DSG & RN  \quad parameter \quad PC: \quad \\
X:\quad & l_{X,Y}l_{X,Z} < l_{X,Y}h_{X,Z} < \theta_{Y,X} < h_{X,Y}l_{X,Z} < h_{X,Y}h_{X,Z} \\
Y: \quad & l_{Y,X}<\theta_{Z,Y}<\theta_{X,Y}<h_{Y,X}\\
Z: \quad & l_{Z,Y}<h_{Z,Y}<\theta_{X,Z}
\end{align*}

\caption{The corresponding STG (Left), MG (Right), and set of inequalities (Bottom) for an example DSGRN parameter for the 3-node network in Fig~\ref{fig:example-networks} exhibiting  a stable partial cycle in the variables $X$ and $Y$ (red) and an unstable partial cycle in $X$ and $Y$ (blue).} 
\label{fig:PC_stg}
\end{figure}

The second example DSGRN parameter contains both a stable and an unstable PC in the variables $X$ and $Y$. The STG and DSGRN parameter inequalities are shown in Fig \ref{fig:PC_stg}. Again the stable cycle is represented by the red arrows and nodes while the unstable cycle is represented in the blue arrows and nodes. Through inspection of the STG it can be seen why the red cycle is stable and the blue cycle is unstable. At any node within the blue cycle, one could keep traveling on the blue cycle or drop down to the red cycle. Once moving in the red cycle one can not make it back to the blue cycle which dictates whether a cycle is stable or unstable. Unlike the FC parameter we only see oscillations in a subset of gene product concentrations in the GRN, seen because the cycle only moves through two of the three dimensions. In particular, the stable cycle does not exist outside of the domain where Z is 0.

The last DSGRN parameter we will look at in this network exhibits an FP. The corresponding STG and inequalities can be seen in Fig \ref{fig:fp_stg}. The FP for this DSGRN parameter can be seen as the single red node with a self-edge within the STG, FP(100). Similar to the last example there exists an unstable PC(X,Y) that is represented in the blue nodes and arrows. Through inspection it can be seen that the (100) domain has no out-edges and therefore must be an FP. Also notice that once the unstable cycle is left one can never return to that cycle and will eventually end up at the FP(100).

\begin{figure}[h!]
\centering
\begin{minipage}{.25\linewidth}
\begin{tikzpicture}[rotate around x=0, rotate around y=0, rotate around z=0, main node/.style={circle, draw, thick, inner sep=1pt, minimum size=0pt}, scale=2.4, node distance=1cm,z={(60:-0.5cm)}]
\tikzstyle{every loop}=[looseness=14]

\node[main node, fill=white] (000) at (0,0,0) {$000$} ;
\node[main node,blue, fill=white] (001) at (0,1,0) {$001$} ;
\node[main node, fill=white] (010) at (1,0,0) {$010$} ;
\node[main node,red, fill=white] (100) at (0,0,1) {$100$} ;
\node[main node, fill=white] (110) at (1,0,1) {$110$} ;
\node[main node,blue, fill=white] (101) at (0,1,1) {$101$} ;
\node[main node,blue, fill=white] (011) at (1,1,0) {$011$} ;
\node[main node,blue, fill=white] (111) at (1,1,1) {$111$} ;
\node[main node, fill=white] (020) at (2,0,0) {$020$} ;
\node[main node, fill=white] (120) at (2,0,1) {$120$} ;
\node[main node,blue, fill=white] (021) at (2,1,0) {$021$} ;
\node[main node,blue, fill=white] (121) at (2,1,1) {$121$} ;

\draw[very thick,-,shorten >= 0pt,shorten <= 0pt]
(101) edge[->,blue] (001)
(001) edge[dashed,->] (000)
(000) edge[dashed, ->] (100)
(100) edge[->,red, loop left] (100)
(101) edge[line width=0.1cm, white] (111)
(011) edge[dashed,->] (010)
(010) edge[dashed,->] (110)
(001) edge[->,blue] (011)
(000) edge[dashed,->] (010)
(110) edge[->] (100)
(111) edge[line width=0.1cm, white] (110)
(111) edge[->] (110)

(120) edge[line width=0.1cm, white] (121)
(111) edge[->,blue] (011)
(111) edge[line width=0.1cm, white] (121)
(111) edge[->,blue] (101)

(021) edge[->,blue] (121)
(101) edge[->] (100)
(121) edge[->] (120)
(021) edge[->] (020)
(020) edge[->] (120)
(011) edge[->,blue] (021)
(121) edge[->,blue] (111)
(120) edge[->] (110)
(010) edge[dashed,->] (020);
\end{tikzpicture}
\end{minipage}%
\hfill
\begin{minipage}{0.25\linewidth}
\begin{tikzpicture}[main node/.style={ellipse, draw, thick, inner sep=3pt, minimum size=0pt},scale=0.8]
\node[main node,red] (00) at (0,0) {$FP(100)$};
\node[main node,blue] (02) at (0,2) {$PC\{X,Y\}$};

\draw[very thick,-,shorten >= 0pt,shorten <= 0pt]
(n2) edge[->] (n1);
\end{tikzpicture}
\end{minipage}%

\begin{align*}
 DSG & RN  \quad parameter \quad FP: \quad \\
X:\quad & l_{X,Y}l_{X,Z} < \theta_{Y,X} < h_{X,Y}l_{X,Z}, \; l_{X,Y}h_{X,Z} < h_{X,Y}h_{X,Z} \\
Y: \quad & l_{Y,X}<\theta_{Z,Y}<\theta_{X,Y}<h_{Y,X}\\
Z: \quad & l_{Z,Y}<h_{Z,Y}<\theta_{X,Z}
\end{align*}

\caption{The corresponding STG and set of inequalities for an example DSGRN parameter for the 3-node network in Fig~\ref{fig:example-networks} exhibiting an FP.}
\label{fig:fp_stg}
\end{figure}

\subsection{Time Series Discretization}\label{sec:timeseries}

To evaluate DSGRN network model consistency with a dataset exhibiting oscillations, we need to extract the sequence of maxima and minima (together called extrema) from a time series dataset. We briefly discuss the methodology from~\cite{berry:2020}, in which for each extremum, a time interval is assigned representing experimental uncertainty in the timing of extrema. 
The idea is that at a specified noise level, an extremum can occur anywhere within the assigned time interval, perhaps due to sparse temporal sampling and/or measurement error. Since the time intervals assigned to different extrema may overlap, the ordering of some extrema with respect to others is indeterminate. However, within any single time series the order of extrema is known.

\begin{figure}[!h]
\centering
\includegraphics[width=.6\textwidth]{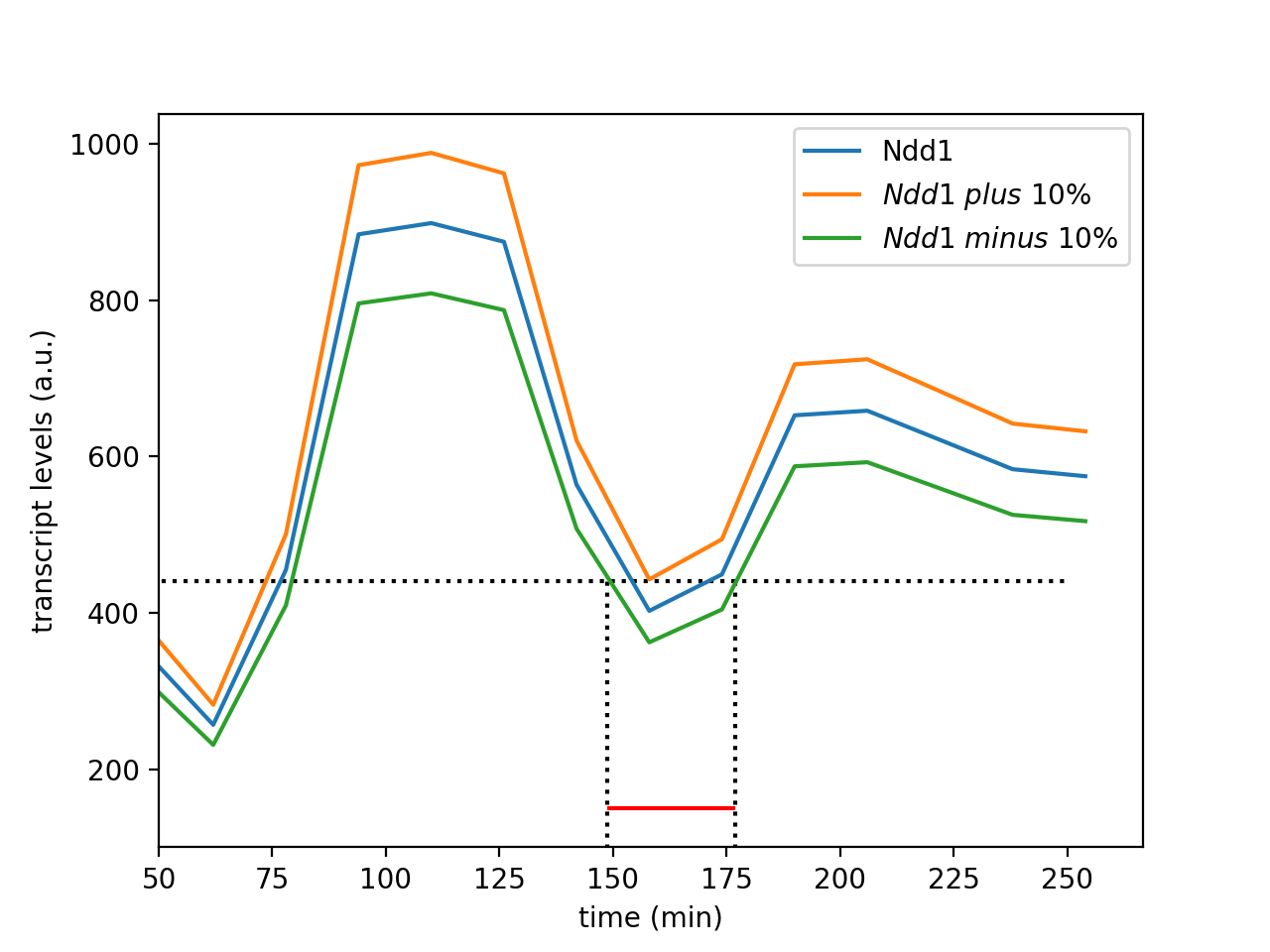}
\caption{The time interval corresponding to a 10\% noise level at the second minimum within the wild-type Ndd1 time series. The red line indicates the size of the time interval. \Bree{fix font for Ndd1}}
\label{fig:time_interval_epsilon}
\end{figure}

Fig~\ref{fig:time_interval_epsilon} schematically shows the process for assigning a time interval for an extremum. The example shown is the gene expression level for Ndd1 in the WT dataset. The blue line is the collected data (with linear interpolation) and the orange and green lines are the original curve $\pm$10\% of the difference between the global maximum and global minimum. This is called a 10\% noise level of the Ndd1 data. The purpose of choosing a noise level is to smooth out small, spurious extrema and also account for imprecision in the timing of large extrema.

The example computation in Fig~\ref{fig:time_interval_epsilon} assigns a time interval to the second minimum in the time series which occurs at 158 minutes, called an extremal interval, by locating a region of the graph around the minimum that is below the $+$10\% curve and above the $-$10\% curve. It was proven in~\cite{cummins2:2018} that any perturbation of the original data that is bounded within $\pm \epsilon$ curves for some noise level $\epsilon$ is guaranteed to have a minimum located within the extremal interval. That is, the existence of the minimum is robust within the extremal interval. The extremal interval for a maximum is recovered in a similar manner. 

This procedure is repeated for every extremum in every time series, so that the end result is a collection of time intervals. The intervals across time series are compared to see if they are overlapping. When two intervals overlap, it is possible for the extrema associated to the intervals to occur in either order in time. In other words, the ordering of the extrema is indeterminate and the extrema are called incomparable. However, when intervals are non-overlapping, then the order between the two extrema is known, and they are called comparable.

A \textit{pattern diagram} is a discrete, graphical representation of all the comparability relationships between extremal intervals. This is known mathematically as the Hasse diagram of a partially ordered set.
In the most straightforward case, all extrema of a dataset are comparable and form a linear sequence of events in time. This is true for any single time series. However, for more than one time series this never happens, since the start of every time series co-occurs and the concentrations for each node are either at a maximum or a minimum.

\begin{figure}[!h]
    \centering
    \includegraphics[scale=0.35]{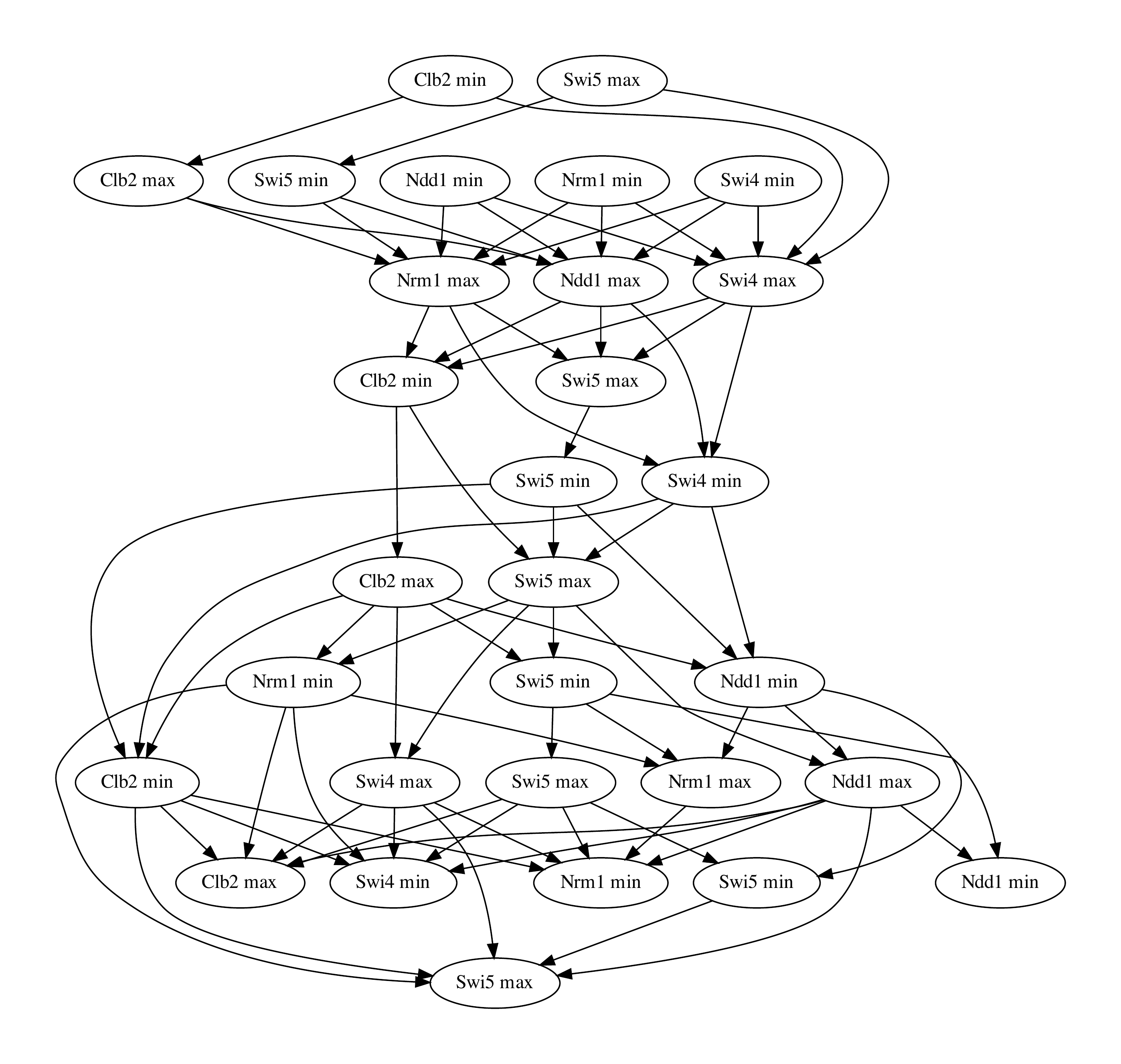}
    \caption{An example pattern diagram for the WT dataset at a noise level of 10\%.}
    \label{fig:ts_poset_g}%
\end{figure}

Fig~\ref{fig:ts_poset_g} shows the pattern diagram for the WT data  (Fig~\ref{fig:ts_poset}b) for the nodes in the mini wavepool. The time intervals for each extremum were computed at 10\% noise, as were all time series discretizations in this manuscript. Each node in the pattern diagram is a type of extremum, with nodes at the top occurring at the beginning of the time series. A directed arrow between a source node and target node means that the source extremum is known to occur earlier in time than the target extremum. If there is no path between two nodes, then the associated extrema are incomparable.

\subsection{Pattern Matching}\label{sec:patternmatch}

Pattern diagrams are the basis for testing DSGRN network model consistency with oscillating time series data. While the technical details differ~\cite{cummins2:2018}, a straightforward demonstration of pattern matching can be done by labeling the edges of a DSGRN state transition graph with the extrema labels that occur in the nodes of the pattern diagram. We begin our example of pattern matching by discussing how to use the pattern diagram, and then proceed to the labeling of the STG.

The key idea is that of a \textit{linear extension} of the pattern diagram. A linear extension is any sequence of all the nodes in a pattern diagram that does not contradict any arrows in the pattern diagram. For example, in Fig~\ref{fig:ts_poset_g}, the top two nodes Clb2 min and Swi5 max may occur in either order in a linear extension because they are incomparable. However, the arrow from Clb2 min to Clb2 max enforces that in any linear extension, Clb2 min occurs before Clb2 max. We claim that any of the linear extensions arising from a pattern diagram might be the ``true'' sequence of extrema in the biological system given the measurement error and sampling density of the data.

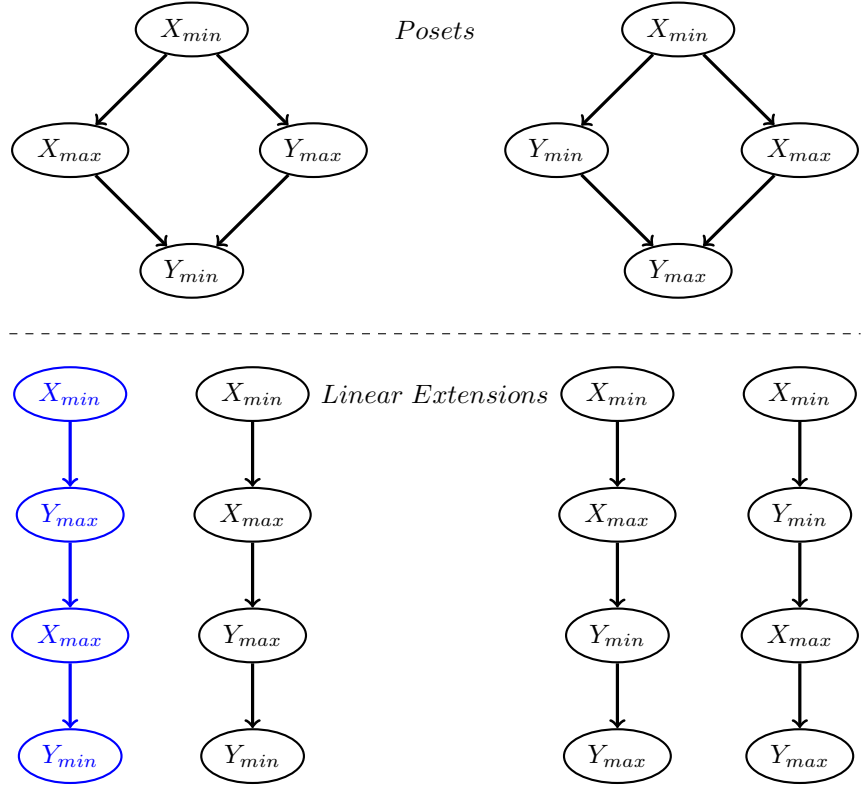
\begin{figure}[!h]
\centering

\begin{tikzpicture}[main node/.style={ellipse,draw, thick, inner sep=3pt, minimum size=0pt}, scale=0.8]

\node[main node, fill=white] (-45) at (-4,5) {$X_{min}$} ;
\node[main node, fill=white] (-23) at (-2,3) {$Y_{max}$} ;
\node[main node, fill=white] (-63) at (-6,3) {$X_{max}$} ;
\node[main node, fill=white] (-41) at (-4,1) {$Y_{min}$} ;

\node[main node, fill=white] (45) at (4,5) {$X_{min}$} ;
\node[main node, fill=white] (23) at (2,3) {$Y_{min}$} ;
\node[main node, fill=white] (63) at (6,3) {$X_{max}$} ;
\node[main node, fill=white] (41) at (4,1) {$Y_{max}$} ;

\node (n1) at (0,5)[ellipse,inner sep=3pt]{$Posets$};

\draw[very thick,-,shorten >= 0pt,shorten <= 0pt]
(-45) edge[->] (-23)
(-45) edge[->] (-63)
(-23) edge[->] (-41)
(-63) edge[->] (-41)
(45) edge[->] (23)
(45) edge[->] (63)
(23) edge[->] (41)
(63) edge[->] (41);

\end{tikzpicture}%

\hfill

\begin{tikzpicture}[main node/.style={ellipse,draw, thick, inner sep=3pt, minimum size=0pt}, scale=0.8]

\draw[dashed] (-7,7)--(7,7);
\node[main node,blue, fill=white] (-66) at (-6,6) {$X_{min}$} ;
\node[main node,blue, fill=white] (-64) at (-6,4) {$Y_{max}$} ;
\node[main node,blue, fill=white] (-62) at (-6,2) {$X_{max}$} ;
\node[main node,blue, fill=white] (-60) at (-6,0) {$Y_{min}$} ;

\node[main node, fill=white] (-36) at (-3,6) {$X_{min}$} ;
\node[main node, fill=white] (-34) at (-3,4) {$X_{max}$} ;
\node[main node, fill=white] (-32) at (-3,2) {$Y_{max}$} ;
\node[main node, fill=white] (-30) at (-3,0) {$Y_{min}$} ;

\node[main node, fill=white] (66) at (6,6) {$X_{min}$} ;
\node[main node, fill=white] (64) at (6,4) {$Y_{min}$} ;
\node[main node, fill=white] (62) at (6,2) {$X_{max}$} ;
\node[main node, fill=white] (60) at (6,0) {$Y_{max}$} ;

\node[main node, fill=white] (36) at (3,6) {$X_{min}$} ;
\node[main node, fill=white] (34) at (3,4) {$X_{max}$} ;
\node[main node, fill=white] (32) at (3,2) {$Y_{min}$} ;
\node[main node, fill=white] (30) at (3,0) {$Y_{max}$} ;

\node (n2) at (0,6) [ellipse,inner sep=3pt]{$Linear \; Extensions$};

\draw[very thick,-,shorten >= 0pt,shorten <= 0pt]
(-66) edge[->,blue] (-64)
(-64) edge[->,blue] (-62)
(-62) edge[->,blue] (-60)
(-36) edge[->] (-34)
(-34) edge[->] (-32)
(-32) edge[->] (-30)
(66) edge[->] (64)
(64) edge[->] (62)
(62) edge[->] (60)
(36) edge[->] (34)
(34) edge[->] (32)
(32) edge[->] (30);
\end{tikzpicture}%

\caption{Two example pattern diagrams corresponding to two hypothetical datasets (top row) and their corresponding collections of linear extensions (bottom row). The linear extension in blue is a pattern match to the unstable partial cycle in Fig~\ref{fig:PC_stg}, as seen by the labeling of the STG in Fig~\ref{fig:search_graph}. }
\label{fig:poset_extensions}
\end{figure}

Continuing our example of the three-node network in Fig \ref{fig:example-networks},we have constructed two example pattern diagrams that could have arisen from two different datasets of nodes $X$ and $Y$ (top row of Fig~\ref{fig:poset_extensions}), with the idea of checking if a pattern match exists between either of these datasets and the unstable partial cycle in Fig~\ref{fig:PC_stg} (blue nodes). Via inspection of the two example pattern diagrams it can be seen which extrema are comparable and incomparable. In the first pattern diagram on the left, the node $X_{min}$ is comparable with any other node in the pattern diagram, but $X_{max}$ and $Y_{max}$ are not comparable with each other. For the pattern diagram on the right, we see a shift in the ordering of the extrema where $X_{max}$ and $Y_{min}$ are incomparable. Each pattern diagram has two linear extensions (bottom row of Fig~\ref{fig:poset_extensions}). To construct the linear extensions for the left pattern diagram a decision is made that defines the ordering of the two incomparable extrema $X_{max}$ and $Y_{max}$. Similarly for the poset on the right, the decision is between the order of $Y_{min}$ and $X_{max}$.

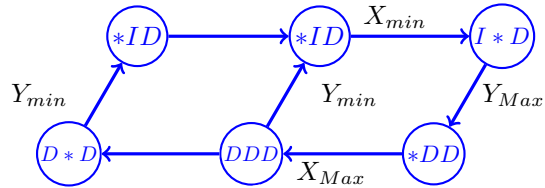
\begin{figure}[h!]
    \centering
\begin{tikzpicture}[rotate around x=0, rotate around y=0, rotate around z=0, main node/.style={circle, draw, thick, inner sep=1pt, minimum size=0pt}, scale=2.4, node distance=1cm,z={(60:-0.5cm)}]

\node[main node,blue, fill=white] (001) at (0,1,0) {$*ID$} ;
\node[main node,blue, fill=white,scale = .8] (101) at (0,1,1.5) {$D*D$} ;
\node[main node,blue, fill=white] (011) at (1,1,0) {$*ID$} ;
\node[main node,blue, fill=white,scale = .8] (111) at (1,1,1.5) {$DDD$} ;
\node[main node,blue, fill=white,scale = .9] (021) at (2,1,0) {$I*D$} ;
\node[main node,blue, fill=white,scale = .9] (121) at (2,1,1.5) {$*DD$} ;

\node (n1) at (-.35,1,.75) {$Y_{min}$};
% \node(n2) at (.5, 1, 1.75) {$Z_{Max}$};
\node (n3) at (1.5, 1, 1.75) {$X_{Max}$};
\node (n4) at (2.25,1,.75) {$Y_{Max}$};
% \node(n5) at (.35, 1, -.25) {$Z_{Max}$};
\node (n6) at (1.35, 1, -.25) {$X_{min}$};
\node (n7) at (1.35, 1, .75) {$Y_{min}$};

\draw[very thick,-,shorten >= 0pt,shorten <= 0pt]
(101) edge[->,blue] (001)
(101) edge[line width=0.1cm, white] (111)
(001) edge[->,blue] (011)
(111) edge[line width=0.1cm, white] (110)
(111) edge[->,blue] (011)
(111) edge[line width=0.1cm, white] (121)
(111) edge[->,blue] (101)
(021) edge[->,blue] (121)
(011) edge[->,blue] (021)
(121) edge[->,blue] (111);

\end{tikzpicture}%
\caption{The part of the STG in Fig \ref{fig:PC_stg} corresponding to the unstable partial cycle that is labeled according to whether node concentrations are increasing or decreasing within the corresponding domains in phase space.}
\label{fig:search_graph}
\end{figure}

The nodes of an STG can be labeled with information about whether a node concentration is increasing ($I$), decreasing ($D$), or both ($*$) within the corresponding domain in phase space (for details, see~\cite{cummins2:2018}). From this information, edge labels in the STG can be deduced that describe which extremum could occur between any two neighboring domains. 
To continue our example, Fig \ref{fig:search_graph} contains the unstable PC from Fig \ref{fig:PC_stg} along with the $I$, $D$, and $*$ labels  on each node in the order $(XYZ)$. Consider the edge going from the node labeled $(DDD)$ to the node labeled $(*ID)$. In the first domain, $Y$ is decreasing and in the second domain, it is increasing. Therefore $Y$ had to undergo a minimum on the edge between the two. Continuing to the node labeled $(I\!*\!D)$, we note that we are moving along the $Y$ axis. $Y$ does not regulate itself, and therefore we know that starting from node $(*ID)$ means $Y$ must continue to increase, although it may decrease starting from another initial condition. Therefore, no extremum may occur in $Y$ on the arrow from $(*ID)$ to $(I\!*\!D)$. However, on the next edge to $(*DD)$, $Y$ transitions from increasing to decreasing, therefore must have achieved a maximum. Similar arguments allow the assignment of the other edge labels.

 To continue with pattern matching, we compare the sequence of extrema in the STG to the linear extensions in Fig~\ref{fig:poset_extensions}. It can be seen that within the STG $X_{min}$ is followed by $Y_{max}$ then $X_{max}$ and finally $Y_{min}$. This is exactly the left-most linear extension in Fig~\ref{fig:poset_extensions}, high-lighted in blue. On the other hand, looking at the rest of the linear extensions in Fig \ref{fig:poset_extensions}, it can be seen that none of them match the path through the STG. We conclude that the three-node network model in Fig~\ref{fig:example-networks} is consistent with the hypothetical dataset that generated the pattern diagram on the left, but is not consistent with the pattern diagram on the right.

 This pattern matching procedure is conducted analogously for the WT and mutant datasets against the mini wavepool model in Section~\ref{sec:results}.

 \subsection{Hill models} \label{sec:hill}

To model an activating regulation $x \to y$, we used a Hill function of the form:
\begin{align*}
H^+_{y,x}(x) = (h_{y,x}-l_{y,x}) \dfrac{x^n}{\theta_{y,x}^n + x^n} + l_{y,x}
\end{align*}
where $l$, $h$, and $\theta$ are the switching system parameters introduced in Section~\ref{sec:dsgrn}. Similarly, a repressing regulation $x \dashv y$ is given by:
\begin{align*}
H^-_{y,x}(x) = (h_{y,x}-l_{y,x}) \dfrac{\theta_{y,x}^n}{\theta_{y,x}^n + x^n} + l_{y,x}
\end{align*}

We created a Hill model of the mini wavepool network in Fig~\ref{fig:wavepool} of the following form:
\begin{align*}
\dot S &= -S + H^-(N)H^-(C)H^+(S)H^+(W) \\
\dot N & = -N + H^+(S) \\
\dot D &= -D + H^+(S)H^+(C) \\
\dot W &= -W + H^-(C)H^+(D) \\
\dot C &= -C + H^+(N)
\end{align*}
where $S$ (Swi4), $N$ (Nrm1), $D$ (Ndd1), $W$ (Swi5), and $C$ (Clb2) represent expression levels of the corresponding mRNA. The subscripts on the Hill functions have been suppressed for clarity.

The parameters for the simulations in Figs~\ref{fig:wt_hill_model},~\ref{fig:phenoII_hill_model}, and~\ref{fig:hill_model_sac} are given as supplementary data files. The Hill exponent $n$ was taken to be 10 in all simulations.

\section*{Data Availability} The figures and results within this paper were generated using the Github repository at \url{https://github.com/julianfox8/2022-DSGRN-Phenotypes-Yeast.git}. This repository includes the jointly normalized data used in this manuscript and the Jupyter notebooks and code used to jointly normalize the datasets, but does not include the original data. That is available by request from the Haase lab at Duke University.

\section*{Acknowledgements}
JF was partially supported by a grant from the Undergraduate Scholars Program at Montana State University. BC and SBD were partially supported by NSF TRIPODS+X grant DMS-1839299 and NIH 5R01GM126555-01. RCM was partially supported by NIH 5R01GM126555-01. MG was partially supported by the National Science Foundation under awards DMS-1839294 and HDR TRIPODS award CCF-1934924, National Institutes of Health award R01 GM126555, FAPESP grant 2019/06249-7, and by CNPq grant 309073/2019-7.

\bibliographystyle{plain}
\bibliography{bibliography.bib}

\end{document}